\documentclass[11pt]{article}

\usepackage[T1]{fontenc}
\usepackage[dvips]{color}
\usepackage{amsmath,amssymb,epsfig}
\usepackage{comment}
\newtheorem{defi}{Definition}[section]
\newtheorem{prop}[defi]{Proposition}
\newtheorem{theo}[defi]{Theorem}
\newtheorem{conj}[defi]{Conjecture}
\newtheorem{lemm}[defi]{Lemma}
\newtheorem{coro}[defi]{Corollary}
\newtheorem{rema}[defi]{Remark}
\newtheorem{exem}[defi]{Example}
\newtheorem{exems}[defi]{Examples}

\newcommand{\bdefi}{\begin{defi}}
\newcommand{\edefi}{\end{defi}}
\newcommand{\bprop}{\begin{prop}}
\newcommand{\eprop}{\end{prop}}
\newcommand{\btheo}{\begin{theo}}
\newcommand{\etheo}{\end{theo}}
\newcommand{\blemm}{\begin{lemm}}
\newcommand{\brema}{\begin{rema}}
\newcommand{\erema}{\end{rema}}
\newcommand{\bexer}{\begin{exem}}
\newcommand{\eexer}{\end{exem}}
\newcommand{\bexems}{\begin{exems}}
\newcommand{\eexems}{\end{exems}}
\newcommand{\bconj}{\begin{conj}}
\newcommand{\econj}{\end{conj}}
\newcommand{\elemm}{\end{lemm}}
\newcommand{\bcoro}{\begin{coro}}
\newcommand{\ecoro}{\end{coro}}
\newcommand{\dem}{\noindent{\bf Proof. }}
\newcommand{\rem}{\noindent{\bf Remark. }}

\usepackage{mathrsfs}
\renewcommand\mathcal{\mathscr}

\newcommand{\A}{{\cal A}}
\newcommand{\R}{{\cal R}}

\renewcommand{\L}{{\cal L}}

\newcommand{\N}{{\cal N}}

\newcommand{\D}{{\cal D}}
\newcommand{\E}{{\cal E}}

\renewcommand{\H}{{\cal H}}
\renewcommand{\O}{{\cal O}}
\newcommand{\C}{{\cal C}}

\renewcommand{\S}{{\cal S}}

\newcommand{\maths}[1]{{\mathbb #1}}  

\newcommand{\RR}{\maths{R}}
\newcommand{\NN}{\maths{N}}
\newcommand{\CC}{\maths{C}}
\newcommand{\QQ}{\maths{Q}}

\newcommand{\HH}{\maths{H}}

\newcommand{\ZZ}{\maths{Z}}
\newcommand{\PP}{\maths{P}}

\newcommand{\ra}{\rightarrow}

\newcommand{\ov}[1]{{\overline #1}} 
\newcommand{\wt}[1]{{\widetilde{#1}}}
\newcommand{\wh}[1]{{\widehat{#1}}}

\newcommand{\ga}{\gamma}
\newcommand{\Ga}{\Gamma}

\newcommand{\bs}{\backslash}

\newcommand{\cqfd}{\hfill$\Box$}

\newcommand{\PSLZ}{\mbox{${\rm{ PSL}}_{2}(\ZZ)$}}

\newcommand{\hdr}{{\HH}^2_\RR}
\newcommand{\htr}{{\HH}^3_\RR}
\newcommand{\hnr}{{\HH}^n_\RR}
\newcommand{\hnc}{{\HH}^n_\CC}

\newcommand{\spd}[1]{\operatorname{Sp}(#1)}

\newcommand{\cp}{{\operatorname{\mathfrak c\mathfrak r%
\mathfrak p}}}
\newcommand{\ftp}{{\operatorname{\mathfrak f\mathfrak t%
\mathfrak p}}}
\newcommand{\diam}{{\operatorname{diam}}}
\newcommand{\Lki}{{\operatorname{Lk_\infty}}}
\renewcommand{\Re}{{\operatorname{Re}}}

\newcommand{\arsinh}{{\operatorname{argsinh}}}
\newcommand{\arcosh}{{\operatorname{argcosh}}}

\newcommand{\hdc}{{\HH}^2_\CC}
\renewcommand{\S}{\textsection}
\newcommand{\Qq}{\QQ_{\rm quad}}
\newcommand{\Kq}{K_{\rm quad}}

\newcounter{fig}

\title{Spiraling spectra of geodesic lines
\\  in negatively curved manifolds}
\author{Jouni Parkkonen \and Fr\'ed\'eric Paulin}
\date{}

\begin{document} 
\maketitle

\begin{abstract} 
\noindent   
Given a negatively curved geodesic metric space $M$, we study the
asymptotic penetration behaviour of geodesic lines of $M$ in
small neighbourhoods of closed geodesics and of other compact 
convex subsets of $M$. We define a {\it spiraling spectrum} which
gives precise information on the asymptotic spiraling lengths of
geodesic lines around these objects. We prove analogs of the theorems
of Dirichlet, Hall and Cusick in this context. As a consequence, we obtain
Diophantine approximation results of elements of $\RR,\CC$ or the
Heisenberg group by  quadratic irrational ones.  
\end{abstract}

%%%%%%%%%%%%%%%%%%%%%%%%
%% Modified 6.6.08/Jouni
%%%%%%%%%%%%%%%%%%%%%%%%
\section{Introduction} 
\label{sec:intro} 

Let $M$ be a  finite volume connected complete Riemannian manifold with
dimension $n$ at least $2$ and sectional curvature at most $-1$. Let
$e$ be an end of $M$, and let $C$ be a closed geodesic in $M$. One of
the aims of this paper is to study the asymptotic spiraling behaviour
of the (locally) geodesic lines in $M$ starting from $e$ around the closed
geodesic $C$.

Just for the sake of normalization, fix a Margulis neighbourhood $N$
of the cusp $e$ in $M$ (see for instance \cite{BK}). Let ${\rm
  Lk}_N(M)$ be the set of  geodesic lines starting from $e$
that first meet $\partial N$ at time $0$, and do not converge to a
cusp of $M$. Let $d_N$ be the Hamenst\"adt distance on ${\rm Lk}_N(M)$
(see \cite{HPMZ}), which is a natural distance inducing the
compact-open topology on ${\rm Lk}_N(M)$, and which coincides with the
induced Riemannian distance on the first intersection points with
$\partial N$ if $N$ has constant curvature.

Let ${\rm Lk}_{N,C}(M)$ be the (countable, dense) set of elements
$\rho$ in ${\rm Lk}_N(M)$ that spiral indefinitely around $C$, that is
such that $\lim_{t\ra+\infty}\;d(\rho(t),C)=0$. For every $r$ in ${\rm
  Lk}_{N,C}(M)$, let $D(r)$ be the shortest length of a path between
$\partial N$ and $C$ which is homotopic (while its endpoints stay in
$\partial N$ and $C$ respectively), for any $t$ big enough, to the
path obtained by following $r$ from $r(0)$ to $r(t)$, and then a
shortest geodesic between $r(t)$ and its closest point on $C$. This
number $D(r)$ naturally measures the wandering of $r$ in $M$ before
$r$ seriously starts to spiral indefinitely around $C$. See the end of
Section \ref{sec:framework} for explicit computations when $M$ is
locally symmetric.

We define the {\it spiraling constant} around $C$ of $\xi\in {\rm
  Lk}_N(M)$ by
$$
c(\xi)=\liminf_{r\in{\rm Lk}_{N,C}(M)\;,\;
D(r)\ra+\infty} \;\;
e^{D(r)}d_N(\xi, r)\;,
$$
which measures how well $\xi$ is approximated by geodesic lines
spiraling indefinitely  around $C$, and, when small, says that,
asymptotically, $\xi$
 has long periods of time during which it spirals around
$C$. We define the {\it spiraling spectrum} around $C$ in $M$ by
$$
\operatorname{Sp}_{N,C}(M)= \big\{c(\xi)\;:\;
\xi\in {\rm Lk}_{N}(M)-{\rm Lk}_{N,C}(M)\big\}\;.
$$

Here is a sample of our results.

\btheo[Dirichlet-type theorem] \label{theo:mainintroun}
The spiraling spectrum $\operatorname{Sp}_{N,C}(M)$ is a
boun\-ded subset of $[0,+\infty[$.
\etheo

\btheo[Cusick-type and Hall-type theorem] \label{theo:mainintrodeux}
If $M$ has constant curvature, then the spiraling spectrum
$\operatorname{Sp}_{N,C}(M)$ is closed. If, in addition, the dimension
of $M$ is at least $3$, then the spectrum contains an interval $[0,c]$
for some $c>0$.  \etheo

When $C$ is replaced by a cusp (and spiraling a long time around $C$
is replaced by having a long excursion in a fixed cusp neighbourhood),
the analogous results are motivated by Diophantine approximation
results: see for instance \cite{For,Coh,Pat,Sul,Ser,Haa,CF,Vul,Dal},
as well as below in this introduction, and Remark
\ref{rem:motivation}.  In that context, the boundedness of the
spectrum was proved in \cite[Theorem 1.1]{HPMZ}, the closedness of the
spectrum was shown in \cite{Mau}, and the existence of a Hall ray was
proved in \cite[Theorem 1.6]{PP2}.

\medskip Although our arithmetic applications are going to be in the
setting defined at the beginning of this introduction, our results are
true in much more general situations (see the beginning of Section
\ref{sec:framework}).  In particular, $M$ does not need to have a cusp
(and for instance could be compact): we may replace $e$ by a point
$x_0$ in $M$, and then consider the geodesic rays starting from
$x_0$. Or $M$ could be allowed to have a compact totally geodesic
boundary, and we may replace $e$ by a connected component
$\partial_0M$ of $\partial M$, considering the geodesic rays starting
from a point of $\partial_0M$ perpendicularly to
$\partial_0M$. Furthermore, $C$ can be replaced by a connected
embedded totally geodesic submanifold of positive nonmaximal
dimension, or by the convex core of a precisely invariant
quasifuchsian subgroup (see for instance \cite{MT} for
definitions). The theorems \ref{theo:mainintroun} and
\ref{theo:mainintrodeux} remain valid under certain more general
hypotheses on $M$ (see the theorems \ref{theo:Dirichlet},
\ref{theo:closedspec} and Corollary \ref{coro:corohalgene} for
statements). Section \ref{sec:hallray}, where we prove the existence
of Hall rays in spiraling spectra, relies on \cite{PP2}. In Section
\ref{subsec:upperbounds}, we also give upper bounds on the spiraling
spectra in several classical examples.

\medskip To conclude this introduction, we give Diophantine
approximation results which follow from the above theorems in
Riemannian geometry.  Recall that for $x\in\RR-\QQ$, the {\it
approximation constant} of $x$ by rational numbers is 
$$
c(x)=
\;\liminf_{p,q\in\ZZ,\; q\ra+\infty}\;\; q^2\Big|x-\frac pq\Big|,
$$
and that the {\it
Lagrange spectrum} is $\operatorname{Sp}_\QQ= \{c(\xi)\;:\;
\xi\in \RR-\QQ\}$. Numerous properties of the Lagrange spectrum are
known (see for instance \cite{CF}). In particular,
$\operatorname{Sp}_\QQ$ is bounded (Dirichlet 1842), has maximum
$\frac{1}{\sqrt{5}}$ (Korkine-Zolotareff 1873, Hurwitz 1891), is
closed (Cusick 1975), contains a {\it Hall ray}, that is a maximal non
trivial interval $[0,\mu]$ (Hall 1947), with $\mu=491 993 569/(2 221
564 096+283 748\sqrt{468})$ (Freiman 1975). Also, recall Khintchine's
result \cite{Khi} saying that almost every real number is badly
approximable by rational numbers.  The following result, which is a
quite particular case of the results of Section \ref{sec:appdioapp},
gives analogous Diophantine approximation results of real numbers by
(families of)  quadratic irrational elements.

For every real  quadratic irrational number $\alpha$ over $\QQ$, let
$\alpha^\sigma$ be its Galois conjugate. Let $\alpha_0$ be a fixed
real quadratic irrational number over $\QQ$. Let $\E_{\alpha_0}={\rm
PSL}_2(\ZZ)\cdot\{\alpha_0,\alpha_0^\sigma\}$ be its (countable, dense
in $\RR$) orbit for the action by homographies and anti-homographies
of ${\rm PSL}_2(\ZZ)$ on $\RR\cup\{\infty\}$.  For instance, if $\phi$
is the Golden Ratio $\frac{1+\sqrt{5}}{2}$, then $\E_{\phi}$ is the
set of real numbers whose continued fraction expansion ends with an
infinite string of $1$'s.

For every $x\in \RR-(\QQ\cup \E_{\alpha_0})$, define the {\it
approximation constant} of $x$ by elements of $\E_{\alpha_0}$, as
%$$
%c_{\alpha_0}(x)=
%\;\liminf_{\alpha\in\E_{\alpha_0}\;:\;h(\alpha)\ra +\infty}\;\;
%h(\alpha)\;|x-\alpha|\;,
%$$
%$$
%c_{\alpha_0}(x)=
%\;\liminf_{\alpha\in\E_{\alpha_0}\;:\;|\alpha-\alpha^\sigma|\ra 0}\;\;
%\frac{2}{|\alpha-\alpha^\sigma|}\;|x-\alpha|\;,
%$$
$$
c_{\alpha_0}(x)=
\;\liminf_{\alpha\in\E_{\alpha_0}\;:\;|\alpha-\alpha^\sigma|\ra 0}\;\;
2\;\frac{|x-\alpha|}{|\alpha-\alpha^\sigma|}\;,
$$
%where $h(\alpha)=2/|\alpha-\alpha^\sigma|$, 
and the corresponding {\it approximation spectrum}, by
$$
\operatorname{Sp}_{\alpha_0} =\{c_{\alpha_0}(x)\;:\;x\in \RR-(\QQ\cup
\E_{\alpha_0})\}\;.
$$

\btheo \label{theo:goldratintro} %
Let $\alpha_0$ be a real quadratic irrational number over $\QQ$.  
Then $\operatorname{Sp}_{\alpha_0}$ is a  closed bounded subset of
$[0,+\infty[\,$.

Furthermore, let $\psi:\;]0,+\infty[\;\ra\;]0,+\infty[$ be a map such
    that $t\mapsto \log(\psi(e^{-t}))$ is Lipschitz. If $\int_{0}^1
    \psi(t)/t^2\;dt$ diverges (resp.~converges), then for Lebesgue
    almost all $x\in\RR$,
$$
\;\liminf_{\alpha\in\E_{\alpha_0}\;:\;|\alpha-\alpha^\sigma|\ra 0}\;\;
\frac{|x-\alpha|}{\psi(|\alpha-\alpha^\sigma|)}=
0\;\;({\rm resp.} = +\infty)\;.
$$  
\etheo

In this particular case,
% the existence of such a nontrivial interval
%$[0,\kappa]$ uses the work \cite{SS}, and 
the last statement can be derived from \cite{BV} or \cite{DMPV}.  
%If $[0,\kappa_{\alpha_0}]$ is the maximal interval containing $0$
%contained in $\operatorname{Sp}_{\alpha_0}$, we do not know the value
%of $\kappa_{\alpha_0}$
%%%%%%
%\note{we should try to get an explicit lower bound on $\kappa$, and 
%it will not follow directly from \cite{PP2}}
%%%%%% 
%(an analog of Freiman's constant). 
Except for the following result, we do not know the exact value of the
maximum $K_{\alpha_0}$ of $\operatorname{Sp}_{\alpha_0}$ (an analog of
Hurwitz's constant). We prove an upper bound $K_{\alpha_0}\leq
(1+\sqrt{2})\sqrt{3}\approx 4.19$ for any $\alpha_0$, see Section
\ref{subsec:upperbounds}.

\bprop For the Golden Ratio $\phi=\frac{1+\sqrt{5}}{2}$, we have
$K_\phi=1-1/\sqrt{5}\approx 0.55 $, and $K_\phi$ is not isolated in
$\operatorname{Sp}_{\phi}$.  \eprop

There are many papers on the Diophantine approximation of real numbers
by algebraic numbers. After pioneering work by Mahler, Koksma, Roth
and Wirsing, the following Dirichlet type theorem has been proved by
Davenport and Schmidt. Let $\Qq$ be the set of real 
quadratic irrational numbers over $\QQ$, and denote by $H(\alpha)$ the naive
height of an algebraic number $\alpha$ (the maximal absolute value of
the coefficients of its minimal polynomial over $\ZZ$).  For every
nonquadratic irrational real number $x$, Davenport and Schmidt
\cite{DS1} proved that
$$
\liminf_{\alpha\in \QQ\cup\Qq\;:\;H(\alpha)\ra+\infty} \;\;
H(\alpha)^{3}\;|x-\alpha|\;\;<\;+\infty\;.
$$
Sprind\v{z}uk \cite{Spr} proved that this result is
generically optimal: For every $\epsilon>0$, for Lebesgue
almost every $x$ in $\RR$,
$$
\liminf_{\alpha\in\QQ\cup\Qq\;:\;H(\alpha)\ra+\infty} \;\;
H(\alpha)^{3+\epsilon}\;|x-\alpha|\;\;=+\infty\;.
$$
We refer to \cite{Bug} and its impressive bibliography for further
references. But note that none of the works that we know of is
approximating by elements in the orbit under integral homographies of
a given algebraic number; almost all of them are approximating using
(a simple function of) the naive height as a complexity, but none
using our complexity $h(\alpha)=2/|\alpha-\alpha^\sigma|$. This
complexity (see \cite[Lem.~5.2]{BPP} for an algebraic interpretation)
behaves very differently from the naive height $H(\alpha)$, even in
such an orbit, see Section \ref{subsec:appdioappRC}.

In Section \ref{sec:appdioapp}, expanding Theorem
\ref{theo:goldratintro}, we will give arithmetic applications
analogous to the results of Dirichlet, Cusick, and Khintchine for the
Diophantine approximation of points of $\RR$ (resp.~$\CC$, the
Heisenberg group ${\rm Heis}_{2n-1}(\RR)$) by classes of quadratic
irrational elements over $\QQ$ (resp.~quadratic irrational elements
over imaginary quadratic extensions of $\QQ$, elements whose
coefficients are rational or quadratic over an imaginary quadratic
extension of $\QQ$), and to Hall's result in $\CC$ and ${\rm
  Heis}_{2n-1}(\RR)$.

\medskip\noindent {\small {\it Acknowledgements : } The first author
  acknowledges the financial support of the Ecole Normale
  Sup\'erieure and Yrj\"o, Vilho ja Kalle V\"ais\"al\"an s\"a\"ati\"o,
  and the second author that of the University of 
  Georgia at Athens, for visits during which part of this paper was
  investigated.  Both authors warmly thank Sa'ar Hersonsky for many
  hours of discussions; furthermore, the likelyhood of the measure
  theoretic statements in the theorems \ref{theo:goldratintro},
  \ref{theo:rc} and \ref{theo:applicapproxheis} was explicited during
  his work with the second author in \cite{HPpre}. We thank Y.~Bugeaud
  and K.~Belabas for very helpful comments concerning the arithmetic
  applications, and Y.~Benoist.}

\section{Preliminaries} 
\label{sec:background}

Throughout the paper, $(X,d)$ will be a proper CAT($-1$) geodesic
metric space, and $\partial_\infty X$ its boundary at infinity. We use
\cite{BH} as a general reference for this section. Unless otherwise
stated, balls and horoballs are closed. If $\epsilon>0$ and $A$ is a
subset of $X$, we denote by $\N_\epsilon A$ the (closed)
$\epsilon$-neighbourhood of $A$ in $X$.

Let $\Gamma$ be a discrete group of isometries of $X$. The limit set
of $\Gamma$ is denoted by $\Lambda\Gamma$. The conical limit set of
$\Gamma$ is denoted by $\Lambda_c\Gamma$. When $\Lambda\Gamma$
contains at least two points, the convex hull of $\Lambda\Gamma$ is
denoted by $\C\Gamma$. The group $\Ga$ is  {\em
  convex-cocompact} if $\Lambda\Ga$ contains at least two points, and
if the action of $\Ga$ on $\C \Ga$ has compact quotient.

We will say that a subgroup $H$ of a group $G$ is {\it almost
  malnormal} if, for every $g$ in $G-H$, the subgroup $gHg^{-1}\cap H$
is finite.  We refer for instance to \cite[Prop.~2.6]{HPpre} for a
proof of the following well known result.

\bprop\label{prop:equivmalnormal}
%Let $X$ be a proper ${\rm CAT}(-1)$ geodesic metric space.  
Let $\Ga_0$ be a convex-cocompact subgroup of 
%a discrete group 
$\Ga$.
%of isometries of $X$. 
The following assertions are equivalent.
\begin{enumerate}
\item[(1)] 
  $\Ga_0$ is almost malnormal in $\Ga$;
\item[(2)] the limit set of $\Ga_0$ is precisely invariant under
  $\Ga_0$, that is for every $\ga\in\Ga-\Ga_0$, the set
  $\Lambda\Ga_0\cap\ga\Lambda\Ga_0$ is empty;
\item[(3)] 
  $\C\Ga_0\cap \ga \C\Ga_0$ is compact for every $\ga\in\Ga-\Ga_0$;
\item[(4)] for every $\epsilon>0$, there exists
  $\kappa=\kappa(\epsilon)>0$ such that ${\rm diam}\big(\N_\epsilon
  \C\Ga_0\cap \ga\N_\epsilon \C\Ga_0\big)\leq \kappa$ for every
  $\ga\in\Ga-\Ga_0$. \hfill\cqfd
\end{enumerate}
\eprop

For every $\xi$ in $\partial_\infty X$, the {\em Busemann function at
  $\xi$} is the map $\beta_\xi$ from $X\times X$ to $\RR$ defined by
$$
\beta_\xi(x,y)=\lim_{t\ra+\infty} d(x,\xi_t)-d(y,\xi_t)\;,
$$
for any geodesic ray $t\mapsto \xi_t$ ending at $\xi$. 

Let $C$ be a nonempty closed convex subset of $X$. We denote by
$\partial_\infty C$ its set of points at infinity, and by $\partial C$
its boundary in $X$. The {\it closest point map} of $C$ is the
map $\pi_C:(X\cup\partial_\infty X)\ra (C\cup\partial_\infty C)$ which
associates to a point $x\in X$ its closest point in $C$ in the usual
sense, which fixes all points of $\partial_\infty C$, and which
associates to a point $\xi\in\partial_\infty X-\partial_\infty C$ the
point of $C$ which minimizes the map $x\mapsto \beta_\xi(x,x_0)$ for
any $x_0$ in $X$.  This map is continuous.

As in \cite{HPpre}, we define the {\em distance-like map}
$d_C:\big(\partial_\infty X-\partial_\infty C\big)^2\to[0,+\infty[$
associated to $C$ as follows: For $\xi,\xi'\in\partial_\infty
X-\partial_\infty C$, let $x=\pi_C(\xi),x'=\pi_C(\xi')$ be their
closest points in $C$. Let $\xi_t,\xi'_t:[0,+\infty[\;\to X$ be the
geodesic rays starting at $x,x'$ and converging to $\xi,\xi'$ as
$t\to\infty$. Let
\begin{equation}\label{eq:distancelike}
d_{C}(\xi,\xi')=\lim_{t\to+\infty}e^{\frac{1}{2}d(\xi_t,\xi'_t)-t}.
\end{equation}
The distance-like map is invariant under the diagonal action of the
isometries of $X$ preserving $C$, and generalizes the visual and
Hamenst\"adt distances: If $C$ consists of a single point $x$, then
$d_C$ is the {\em visual distance} on $\partial_\infty X$ based at $x$
(see for instance \cite{Bou}), and we denote it by $d_x$. If $C$ is a
ball, then the distance-like map $d_C$ is a positive constant multiple
of the visual distance based at the center of $C$.  If $C$ is a
horoball with point at infinity $\xi_0$, then $d_C$ is the {\em
  Hamenst\"adt distance} on $\partial_\infty X-\{\xi_0\}$, and we also
denote it by $d_{\xi_0,\partial C}$ to put the emphasis on $\xi_0$.

Although $d_C$ is not always an actual distance on $\partial_\infty
X-\partial_\infty C$, it follows from \cite{HPpre} that for every
$\epsilon>0$, there exists $\eta>0$ such that for every
$\xi,\xi',\xi''$ in $\partial_\infty X-\partial_\infty C$, if
$d_{C}(\xi,\xi')<\eta$ and $d_{C}(\xi',\xi'')< \eta$, then
$d_{C}(\xi,\xi'')< \epsilon$.  Indeed, if $d_{C}(\xi,\xi')$ and
$d_{C}(\xi',\xi'')$ are small, then by \cite[Lemma 2.3 (4)]{HPpre}, the
geodesic lines $]\xi,\xi'[$ and $]\xi',\xi''[$ are far from $C$. By
hyberbolicity of $X$, the geodesic line $]\xi,\xi'''[$ is also far
from $C$. Hence $\pi_C(\xi)$ and $\pi_C(\xi'')$ are close. Therefore,
by \cite[Lemma 2.3 (3)]{HPpre}, the value of $d_{C}(\xi,\xi'')$ is
small. In particular, the family of subsets
$$
\Big\{W_n=\big\{(\xi,\xi')\in \big(\partial_\infty X-\partial_\infty
C\big)^2\;:\;d_{C}(\xi,\xi')< \frac{1}{n+1}\big\} \Big\}_{n\in\NN}
$$
is a countable separating system of entourages of a metrisable
uniform structure on $\partial_\infty X-\partial_\infty C$ (see
\cite[TG II.1]{Bou}), whose induced topology is the usual one, by
\cite[Lemma 2.3 (1)]{HPpre}, and which is invariant by the diagonal
action of the isometries of $X$ preserving $C$.

\medskip The {\em crossratio} of four pairwise distinct points
$a,b,c,d\in \partial_\infty X$ is
\begin{equation}\label{eq:crossdefn}
[a,b,c,d]=  
\frac{1}{2}\lim _{t\rightarrow+\infty}
                d(a_t,c_t)-d(c_t,b_t)+d(b_t,d_t)-d(d_t,a_t),
\end{equation}
where $a_t,b_t,c_t,d_t$ are any geodesic rays converging to $a,b,c,d$,
respectively. For the existence of the limit and the properties of the
crossratios, see \cite{Ota} where the order convention is different,
and \cite{Bou} whose crossratio is the exponential of ours; we will be
using the same expression as in \cite{HPCMH,PP2}. If $x_0\in X$, then
$$
[a,b,c,d]=  
\log\;\frac{d_{x_0}(c,a)}{d_{x_0}(c,b)}\;
               \frac{d_{x_0}(d,b)}{d_{x_0}(d,a)}\;.
$$
If $H$ is a horosphere with center $\xi\in\partial_\infty X$, then
for $a,b,c,d\in\partial_\infty X-\{\xi\}$,
$$
[a,b,c,d]=  
\log\;\frac{d_{H}(c,a)}{d_{H}(c,b)}\;
               \frac{d_{H}(d,b)}{d_{H}(d,a)}\;.
$$
If $\xi$ and $a$ coincide, the above expression simplifies as follows,
see \cite[Section 3.1]{PP2}: 
\begin{equation}\label{eq:formcrosham}
[\xi,b,c,d]  =   \log\;\frac{d_{H}(d,b)}{d_{H}(c,b)}\;.
\end{equation}

\bigskip 
Let $\xi\in X\cup\partial_\infty X$. We say that a geodesic line
$\rho:\;]-\infty,+\infty[\;\ra X$ (resp.~geodesic ray
$\rho:[\iota_0,+\infty[\; \ra X$) {\em starts} from $\xi$ if
$\xi=\rho(-\infty)$ (resp.~$\xi=\rho(\iota_0)$). We denote by $T^1_\xi
X$ the space of geodesic rays or lines starting from $\xi$, endowed
with the compact-open topology.

In \cite{PP2}, the penetration of geodesic rays and lines in
neighbourhoods of convex subsets of $X$ was studied by means of
penetration maps.  We now recall the definitions of three classes of
such maps $\ell_C,\cp_L,\ftp_L:T^1_{\xi} X\to[0,+\infty]$, where $C$
is the (closed) $\epsilon$-neighbourhood of a closed convex subset in
$X$ for some $\epsilon>0$, and $L$ is a geodesic line in $X$, with
endpoints $L_1,L_2$.
\begin{enumerate}
\item[(1)] The {\it penetration length map} $\ell_C$ associates to
  every $\rho$ in $T^1_{\xi} X$ the length of the intersection of $C$
  and of the image of $\rho$. (This intersection is connected by
  convexity; there was the assumption in \cite{PP2} that $\xi\notin
  C\cup\partial_\infty C$, which is not necessary here.)
\item[(2)] The {\it fellow-traveller penetration map} $\ftp_L$ is
  defined by
  $$
  {\ftp}_L:\rho\mapsto d\big(\pi_L(\xi),\pi_L(\rho(+\infty))\big)\;,
  $$ with the convention that this distance is $+\infty$ if
  $\pi_L(\xi)$ or $\pi_L(\rho(+\infty))$ is in $\partial_\infty X$
  (there was the assumption in \cite{PP2} that $\xi\notin \N_\epsilon
  L\cup\partial_\infty L$ where $\epsilon>0$ was arbitrary but fixed,
  which is not necessary here).
\item[(3)] When $\xi\in\partial_\infty X$, the {\it crossratio
  penetration map} ${\cp}_{L}$ is defined by
  $$
  {\cp}_L:\rho\mapsto \max\big\{0,[\xi,L_1,\rho(+\infty),L_2],
  [\xi,L_2,\rho(+\infty),L_1]\big\}\;,
  $$ if $\xi,\rho(+\infty)\notin\{L_1,L_2\}$, and ${\cp}_L(\rho)=
  +\infty$ otherwise (there was the assumption in \cite{PP2} that
  $\xi\notin \partial_\infty L$, which is not necessary here).
\end{enumerate}

It is shown in Section 3.1 of \cite{PP2} (and it is easy to see that
the result is still true if $\xi$ belongs respectively to
$C\cup\partial_\infty C$, $\N_\epsilon L\cup\partial_\infty L$,
$\partial_\infty L$) that the above maps are continuous and that, for
every $\epsilon>0$, we have, with the convention that $x-y=0$ if
$x=y=+\infty$, the following inequalities
\begin{equation}\label{eq:normftpl}
||{\ftp}_L-\ell_{\N_\epsilon L}||_\infty=\sup_{\rho\in T^1_{\xi}X}\;
|{\ftp}_L(\rho)-\ell_{\N_\epsilon L}(\rho)|\leq
2\,c'_1(\epsilon)+2\epsilon\;,
\end{equation} 
where $c'_1(\epsilon)= 2\,\arsinh(\coth\epsilon)$, and when
$\xi\in\partial_\infty X$,
\begin{equation}\label{eq:normcrpftp}
||{\cp}_L-{\ftp}_L||_\infty=\sup_{\rho\in T^1_{\xi}X}\;
|{\cp}_L(\rho)-{\ftp}_L(\rho)|\leq 4\log(1+\sqrt{2})\;.
\end{equation}

In constant curvature, the crossratio penetration map has the
following geometric interpretation (see for instance \cite[\S
7.23--7.24]{Bea}, and \cite[\S V.3]{Fenchel}, for related formulas).
Recall that the {\it complex distance} $\ell+i\theta$ between two
oriented geodesic lines $\ga$ and $L$ (in this order) in $\hnr$ is
defined as follows. It is $0+i0$ if they are simultaneously asymptotic
at $+\infty$ or at $-\infty$, and $0+i\pi$ if the terminal point at
infinity of one is the original point at infinity of the
other. Otherwise, if $[p,q]$ is the common perpendicular arc (with
$p=q$ the common intersection point of $\ga$ and $L$ if they
intersect), where $p\in L$, then $\ell=d(p,q)$ and $\theta$ is the
angle at $p$ between the parallel transport of $\ga$ along $[p,q]$ and
$L$.

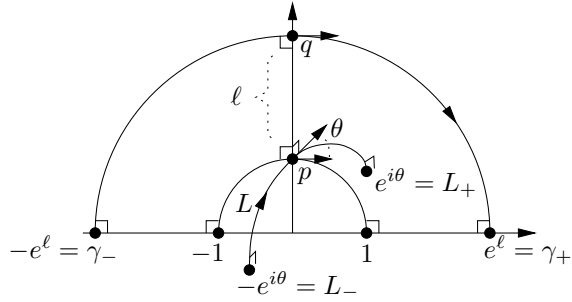
\begin{center}
\input{fig_crossratiogeom.pstex_t}
\end{center}

\blemm \label{lem:calccrossratiogeom} Let $\ga$ and $L$ be oriented
geodesic lines in $\hnr$ with pairwise distinct endpoints
$\ga_-,\ga_+$ and $L_-,L_+$, respectively, and with complex distance
$\ell+i\theta$. Then
$$
[\ga_-,L_-,\ga_+,L_+]=-\log\;\frac{\cosh\ell+\cos\theta}2\;.
$$
In particular,
$$
\cp_L(\ga)=\max\big\{0,-\log\;\frac{\cosh\ell\pm\cos\theta}2\big\}\;.
$$
\elemm

\dem Using isometries, we may assume that $\ga$ and $L$ are both
contained in the upper halfspace model of $\HH^3_\RR$, that the common
perpendicular segment $[p,q]$ is on the vertical axis, with $p$ at
Euclidean height one and $q$ above $p$, and that $\ga_+$ is a positive
real number. By an easy computation, we then have
$$
[\ga_-,L_-,\ga_+,L_+]=
\log\;\frac{|\ga_+-\ga_-|\;|L_+-L_-|}{|\ga_+-L_-|\;|L_+-\ga_-|}
=\log\;\frac{4e^\ell}{|e^\ell+e^{i\theta}|^2}=
-\log\;\frac{\cosh\ell+\cos\theta}2\;.
$$
The result follows.
\cqfd

\medskip
Note that in $\hdr$, we have $\ell>0$ if and only if $\theta=0$
or $\theta=\pi$.

\section{The approximation and spiraling spectra} 
\label{sec:framework} 

In this section, we set up the general framework for our approximation
results. We begin by the definition of the quadruples of data that we
study.

\medskip
%\bdefi  [Quadruple of data $\D$]\label{defi:quadata}
\noindent{\bf The definition of $\D$.}  {\em Let $\Ga$ be a discrete
  group of isometries of a proper ${\rm CAT}(-1)$ geodesic metric
  space $X$.  Let $\Ga_0$ be an almost malnormal convex-cocompact
  subgroup of infinite index in $\Ga$ and let $C_0=\C\Gamma_0$ be the
  convex hull of $\Ga_0$.  Let $C_\infty$ be a nonempty closed convex
  subset of $X$, with stabilizer $\Ga_\infty$ in $\Ga$.  Assume that
  $C_\infty$ does not contain $\C\Ga$, that
  $\Ga_\infty\backslash \partial C_\infty$ is compact, and that the
  intersection of $C_\infty$ and $\ga C_0$ is nonempty for only
  finitely many classes $[\ga]$ in $\Ga_\infty\backslash \Ga/\Ga_0$.
  We will denote the quadruple of data $(X,\Ga,\Ga_0,C_\infty)$ by
  $\D$.}
%\edefi

\medskip
Since $\Lambda\Ga_0$ has at least two points and since $\Ga_0$ has
infinite index in $\Ga$, note that $\Ga$ is nonelementary. By
Proposition \ref{prop:equivmalnormal} (3) and since $C_0$ is
noncompact, the subgroup $\Ga_0$ is the stabilizer of $C_0$ in $\Ga$.
By the discreteness of $\Ga$ and the cocompactness of $\Ga_0$ on $C_0$, a
compact subset of $X$ intersects only finitely many $\ga C_0$ for
$\ga\in \Ga/\Ga_0$.

\medskip Recall that the {\em distance} between two subsets
$A,B$ of $X$ is $d(A,B)=\inf_{a\in A,b\in B}d(a,b)$. For every
$r=[\ga]$ in $\Ga_\infty\backslash \Ga/\Ga_0$, define
$$
D(r)=d(C_\infty,\ga C_0)\;,
$$
which does not depend on the choice of the representative $\ga$ of
$r$. By the cocompactness of the action of $\Ga_\infty$ on $\partial
C_\infty$ and the fact that only finitely many translates of $C_0$
meet a given compact subset, the intersection $C_\infty\cap \ga C_0$
is empty if and only if $D(r)>0$. By convexity, this condition implies
that $\partial_\infty C_\infty\cap \ga\partial_\infty C_0$ is also
empty.  For the same reasons, the following result holds, see
\cite[Lemma 4.1]{HPpre} for a proof in the case $C_\infty=C_0$.

\blemm \label{lem:depthgrows}
For every $T\geq 0$, there are only finitely many elements $r$ in 
$\Ga_\infty\backslash \Ga/\Ga_0$ such that $D(r)\leq T$. \cqfd
\elemm

\blemm \label{lem:infinitedoublecoset} 
The set of double cosets $\Ga_\infty\backslash \Ga/\Ga_0$ is infinite.
\elemm

\dem We first claim that there exists a hyperbolic element $\ga$ in
$\Ga$ whose attractive fixed point $\ga_+$ does not belong to
$\partial_\infty C_\infty$.  The limit set $\Lambda\Ga$ is the closure
of the set of attractive fixed points of elements of $\Ga$, since
$\Ga$ is nonelementary.  Thus, if no such $\ga_+$ exists, $\Lambda\Ga$
is contained in $\partial_\infty C_\infty$, which contradicts the
hypotheses on $\D$ by the convexity of $C_\infty$.

By Lemma \ref{prop:equivmalnormal} (2), and since $\Ga\neq \Ga_0$, we
have $\bigcap_{\ga'\in\Ga}\ga'\Lambda\Ga_0=\emptyset$, hence there
exists $\ga'\in \Ga$ such that the repulsive fixed point $\ga_-$ of
$\ga$ does not belong to $\ga'\Lambda\Ga_0$. Hence the sequence of
closed subsets $(\ga^n\ga'C_0)_{n\in\NN}$ of the compact space $X\cup
\partial_\infty X$ converges to the singleton $\{\ga_+\}$ as $n$ goes
to $+\infty$. This implies that $D([\ga^n\ga'])= d(C_\infty, \ga^n
\ga'C_0)$ converges to $+\infty$ as $n$ goes to $+\infty$. In
particular, the set $\{D(r): r\in \Ga_\infty\backslash \Ga/\Ga_0\}$ is
infinite, and the result follows.  \cqfd

\medskip 
The {\em link} of $\D$ (which depends only on $X,\Ga$ and $C_\infty$)
is
$$
\Lki=
%\Lki(\Ga,C_\infty)=
\Ga_\infty \bs (\Lambda_c\Ga-\partial_\infty C_\infty).
$$ 
The quotient space $\Ga_\infty \bs (\Lambda\Ga-\partial_\infty
C_\infty)$, which contains $\Lki$, is compact, since the closest point
map from $\partial_\infty X -\partial_\infty C_\infty$ to $\partial
C_\infty$ is continuous and $\Ga_\infty$-equivariant. Furthermore,
$\Lki$ is dense in $\Ga_\infty \bs (\Lambda\Ga-\partial_\infty
C_\infty)$.  For every $r=[\ga]$ in $\Ga_\infty\bs\Ga/\Ga_0$ such that
$D(r)>0$, let
$$
\Lambda_r=\pi_\infty(\ga\partial_\infty C_0)
$$
be the image by the canonical projection
$\pi_\infty:\Lambda_c\Ga-\partial_\infty C_\infty\ra {\rm Lk}_\infty$
of $\ga\partial_\infty C_0=\ga\Lambda\Ga_0$.  Note that
$\ga\partial_\infty C_0$ is indeed contained in $\Lambda_c\Ga$ since
$\Ga_0$ is convex-cocompact, and that $\partial_\infty C_\infty$ is
disjoint from $\ga\partial_\infty C_0$ if $D(r)>0$, as seen before.
Furthermore, the sets $\Lambda_r$ are compact subsets of ${\rm
  Lk}_\infty$, that are pairwise disjoint by Lemma
\ref{prop:equivmalnormal} (2), and the union
\begin{equation}\label{eq:deflkio}
{\rm Lk}_{\infty,0}=
\bigsqcup_{r\in\Ga_\infty\backslash\Ga/\Ga_0\,,\;D(r)>0}\Lambda_r
\end{equation}
is dense in ${\rm Lk}_\infty$. In this paper, we study how
dense ${\rm Lk}_{\infty,0}$ is in ${\rm Lk}_\infty$.

Let $\wt d_\infty:(\partial_\infty X-\partial_\infty C_\infty)^2
\ra[0,+\infty[$ be the distance-like map associated to $C_\infty$, and let
$d_\infty$ be its quotient map on ${\rm Lk}_\infty$, which defines, as
in Section \ref{sec:background}, a metrisable uniform structure on
$\Ga_\infty\bs(\partial_\infty X-\partial_\infty C_\infty)$, inducing
the quotient topology.  We endow the 
%(infinite) 
double coset space
$\Ga_\infty\bs\Ga/\Ga_0$ with the Fr\'echet filter of the complements
of the finite subsets, and denote by ${\displaystyle \liminf_r}\;
f(r)$ the lower limit of a real valued map $f$ along this filter.  The
{\it approximation constant} of $\xi\in{\rm Lk}_\infty-{\rm
  Lk}_{\infty,0}$ is
\begin{equation}\label{eq:apprconst}
c(\xi)=\liminf_r \;e^{D(r)}d_\infty(\xi,\Lambda_r)\;,
\end{equation}
and the subset of $[0,+\infty]$ defined by
$$
%{\rm Sp}_{C_\infty}(\Ga_0)=
\spd{\D} \;=
\big\{c(\xi)\;:\;\xi\in {\rm Lk}_\infty-{\rm Lk}_{\infty,0}\big\}
$$
is called the {\it approximation spectrum} of points of 
${\rm Lk}_\infty$ by points of ${\rm Lk}_{\infty,0}$.
We define the {\it Hurwitz constant} of $\D$ as
$$
K_\D=\sup \,\spd{\D}\;\in[0,+\infty]\;.
$$

\brema\label{rem:motivation} 
{\rm Let us give some background and motivations for the
terminology introduced in this paper:
In the definition of the quadruple of data $\D$, let us specialise to the
situation when $X$ is a Riemannian manifold with pinched negative
curvature and $\Ga$ is geometrically finite. If we change the assumptions
on $\Ga_0$ and $\Ga_\infty$ such that  
$\Ga_0=\Ga_\infty$ is the stabilizer of a parabolic fixed
point $\xi_\infty$ of $\Ga$ and $C_\infty$ is the maximal precisely
invariant horoball centered at $\xi_\infty$, then we recover the
framework of Diophantine approximation in negatively curved manifolds
that was developped in \cite{HPMZ,HPsurv,HPETDS,PP2,PP3}.  In this
situation, $\Ga_0$ is not convex-cocompact and  $C_0=C_\infty$, and the 
new quadruple does not have the properties we require of the quadruples of
data in this paper. However, if we take 
 ${\rm
  Lk}_\infty =\Ga_\infty\bs\Lambda_c\Ga$,
$\Lambda_{[\ga]}=\Ga_\infty\ga(+\infty)$, ${\rm
  Lk}_{\infty,0}=\pi_\infty(\Ga\cdot\infty)$, all the constructions in
Section \ref{sec:framework} are still valid.

In particular, let $X$ be the upper halfplane model of the real
hyperbolic plane $\hdr$, let $\Ga=\PSLZ$, let $C_\infty$ be the
horoball in $X$ of points having Euclidean height at least $1$, let
$\Ga_0=\Ga_\infty$ be the cyclic group generated by $z\mapsto z+1$,
and let $\D=(X,\Ga,\Ga_0,\Ga_\infty)$. Then (see \cite[section
2.3]{HPMZ}, \cite{PP3}) ${\rm Lk}_\infty =(\RR-\QQ)/\ZZ$; for every
$r=[\ga]\in\Ga_\infty\bs(\Ga-\Ga_\infty)/\Ga_\infty$, we have
$D(r)=2\log q$ if $\ga\infty=p/q$ with $p\in\ZZ$ and $q\in\NN-\{0\}$
relatively prime; for every $\xi\in\RR-\QQ$, the approximation
constant $c(\xi \!\mod \ZZ)$ is the classical approximation constant
of the irrational number $\xi$ by rational numbers; the approximation
spectrum $\spd{\D}$ is the classical Lagrange spectrum, and the
Hurwitz constant $K_\D$ is the classical Hurwitz constant
$\frac{1}{\sqrt{5}}$ (see the introduction for the definition of these
 objects).
}
\erema

The end of this section is devoted to the study of geometric
examples. 

\medskip Let $M$ be a nonelementary complete connected Riemannian
manifold with sectional curvature at most $-1$, and dimension at least
$2$. Let $A_0$ be a closed geodesic in $M$, not necessarily simple
(for more general $A_0$'s, as for instance in the introduction, we
refer to the general setup). Let $A_\infty$ be a closed codimension
$0$ submanifold of $M$ with smooth connected compact locally convex
boundary, disjoint from $A_0$.

Recall that a locally geodesic ray $\rho$ in $M$ is {\it recurrent}
if, as a map from $[0,+\infty[$ to $M$, it is not proper, that is if
there exist a compact subset $K$ of $M$ and a sequence
$(t_n)_{n\in\NN}$ in $[0,+\infty[$ converging to $+\infty$ such that
$\rho(t_n)\in K$ for every $n$. We say that a locally geodesic ray
$\rho$ in $M$ {\it spirals around} $A_0$ if $d(\rho(t),A_0)$ converges
to $0$ as $t$ goes to $+\infty$.

Let ${\rm Lk}_{A_\infty}(M)$ be the set of recurrent locally geodesic
rays starting perpendicularly from $\partial A_\infty$ and exiting
$A_\infty$, and let ${\rm Lk}_{A_\infty,A_0}(M)$ be the subset of
elements of ${\rm Lk}_{A_\infty}(M)$ that spiral around $A_0$.  Recall
that a geodesic line in a complete simply connected manifold that
crosses a horosphere perpendicularly starts from (up to time reversal)
the point at infinity of this horosphere. Hence when $A_\infty$ is a
small Margulis neighbourhood $N_\infty$ of a cusp $e_\infty$ with
compact boundary, it is equivalent to require that a geodesic ray
exits perpendicularly from $A_\infty$ or that the negative subray of
the geodesic line containing it is a minimizing geodesic ray starting
from the boundary of $N_\infty$ and converging to $e_\infty$.

For every $\rho,\rho'$ in ${\rm Lk}_{A_\infty}(M)$, and every $t\in
[0+\infty[$, let $\ell_t$ be the shortest length of a path homotopic
(relative to the endpoints) to the path obtained by following (the
inverse of) $\rho$ from $\rho(t)$ to $\rho(0)$, then following a
shortest path contained in $A_\infty$ between $\rho(0)$ and
$\rho'(0)$, then following $\rho'$ from $\rho'(0)$ to $\rho'(t)$;
define
\begin{equation}\label{eq:defdinftydown}
\delta_{A_\infty}(\rho,\rho')=
\lim_{t\ra+\infty} e^{\frac{1}{2}\ell_t - t}\;.
\end{equation}
(We will show below that the limit does exist).  For every $\bar r$ in
${\rm Lk}_{A_\infty,A_0}(M)$, let $D(\bar r)$ be the shortest length
of a geodesic segment $[a,b]$ with $a$ in $\partial A_\infty$, $b$ in
$A_0$ such that there exists a (locally) geodesic ray $\rho$ starting
from $b$, contained in $A_0$, such that the path obtained by following
$[a,b]$ from $a$ to $b$ and then $\rho$ is properly homotopic to $\bar
r$ while its origin remains in $\partial A_\infty$.  The {\it
  spiraling constant} around $A_0$ of an element $\bar\xi$ of ${\rm
  Lk}_{A_\infty}(M)$ is
$$
c(\bar\xi)=\liminf_{\bar r\in{\rm Lk}_{A_\infty,A_0}(M)\;,\;
D(\bar r)\ra+\infty} \;\;
e^{D(\bar r)}\delta_{A_\infty}(\bar\xi,\bar r)\;,
$$
and the subset of $[0,+\infty]$ defined by 
$$
\operatorname{Sp}_{A_\infty,A_0}(M)= \big\{c(\bar\xi)\;:\;
\bar\xi\in {\rm Lk}_{A_\infty}(M)-{\rm Lk}_{A_\infty,A_0}(M)\big\}
$$
is called the {\it spiraling spectrum} of geodesic rays in ${\rm Lk}
_{A_\infty}(M)$ around $A_0$.  These notions coincide with the
similarly named ones in the introduction if $M$ has finite volume and
$A_\infty$ is the chosen Margulis neighbourhood of the cusp $e$.

\medskip 
To see the connection with the framework outlined at the beginning of
this section, we may define a quadruple of data
$$
\D_{M,A_0,A_\infty}=(X,\Ga,\Ga_0,C_\infty)\,
$$
%%%%%%
% satisfying Definition \ref{defi:quadata}
%%%%%%%% 
as follows. If $\wt M\ra M$ is a universal Riemannian covering of $M$
with covering group $\Ga$, let $X=\C\Ga$ be the convex hull of $\Ga$,
let $\Ga_0$ be the stabilizer in $\Ga$ of a fixed lift $C_0$ of $A_0$
to $\wt M$, and let $C_\infty$ be the intersection with $X$ of a fixed
connected component $\wt A_\infty$ of the preimage of $A_\infty$ in
$\wt M$.

Note that the image in $M$ of a geodesic ray $\rho$ in $\wt M$ is
recurrent if and only if the endpoint at infinity of $\rho$ is a
conical limit point of $\Ga$. Consider the map $\wt \Phi$ from
$\Lambda_c\Ga-\partial_\infty C_\infty$ to ${\rm Lk}_{A_\infty}(M)$,
which associates to an element $\xi$ of $\Lambda_c\Ga-\partial_\infty
C_\infty$ the image in $M$ of the geodesic ray in $X$ starting from
the closest point on $\wt A_\infty$ to $\xi$ and converging to
$\xi$. By taking the quotient by $\Ga_\infty$, this map induces a
homeo\-morphism $\Phi:{\rm Lk}_{\infty}\ra{\rm Lk}_{A_\infty}(M)$,
which maps ${\rm Lk}_{\infty,0}$ to ${\rm Lk}_{A_\infty,A_0}(M)$.  By
construction, the map $\Phi$ preserves the maps $d_\infty$ and
$\delta_{A_\infty}$ (which proves along the way that the limit in
\eqref{eq:defdinftydown} exists).

For every $\bar r$ in ${\rm Lk}_{A_\infty,A_0}(M)$, by definition of
${\rm Lk}_{\infty,0}$ (see Equation \eqref{eq:deflkio}), there exists
a unique $r$ in $\Ga_\infty\backslash\Ga/\Ga_0$ such that
$\Phi^{-1}(\bar r)$ belongs to $\Lambda_r$. The map $\bar r\mapsto r$
from ${\rm Lk}_{A_\infty,A_0}(M)$ to $\Ga_\infty\backslash\Ga/\Ga_0$
satisfies $D(r)=D(\bar r)$ if $D(r)>0$; the complementary subset in
$\Ga_\infty\backslash\Ga/\Ga_0$ of its image is finite, since there
are only finitely many $r\in \Ga_\infty\backslash\Ga/\Ga_0$ such that
$D(r)\leq 0$; every point in its image has at most two preimages,
since $\Lambda_r$ has at most two points.

Hence, for every $\xi$ in ${\rm Lk}_{\infty}$, we have by construction
$c(\xi)=c(\Phi(\xi))$.  Therefore, as $\Phi$ is surjective,
\begin{equation}\label{eq:divspec}
\spd{\D_{M,A_0,A_\infty}}=
\operatorname{Sp}_{A_\infty,A_0}(M)\;,
\end{equation}
and we conclude that  to obtain results on the spiraling spectrum, it
is sufficient to prove results on the approximation spectrum.

\bigskip

\noindent{\bf Example 1: Spiraling around a closed geodesic in a real
  hyperbolic manifold}

\medskip We will use the upper halfspace model of the real hyperbolic
$n$-space $\HH_\RR^n$, with constant sectional curvature $-1$, so that
$\partial_\infty \HH^n_\RR= \RR^{n-1}\cup\{\infty\}$.  Let $\Ga$ be a
nonelementary discrete subgroup of isometries of $\HH_\RR^n$. Assume
that $\infty$ is a parabolic fixed point of $\Ga$, with stabilizer
$\Ga_\infty$, and that the interior of the horoball $\H_1$ of points
of Euclidean height at least one is {\it precisely invariant} under
$\Ga_\infty$ (that is for all $\ga\in\Ga-\Ga_\infty$, the horoballs
$\H_1$ and $\ga \H_1$ have disjoint interior). When $\Ga$ is
torsion-free and has finite covolume, $\H_1$ covers a Margulis
neighbourhood of a cusp of $\Ga\backslash \HH^n_\RR$. Define $C_\infty
= \H_1\cap \C\Gamma$, and assume that $\Ga_\infty \backslash \partial
C_\infty$ is compact (for instance if $\Ga$ is geometrically finite,
and in particular if $\Ga$ has finite covolume).  Let $\ga_0$ be a
hyperbolic element of $\Ga$, with translation axis $C_0$.  Let $\Ga_0$
be the stabilizer of $C_0$ in $\Ga$, which contains the cyclic group
generated by $\ga_0$ as a subgroup of finite index.  The quadruple
$\D= (\C\Ga,\Ga,\Ga_0,C_\infty)$ satisfies the hypotheses of the
beginning of the section.

\medskip
\noindent\begin{minipage}{6cm}
%\begin{center}
%\!\!\!\!\!\!\!\!\!\!
~~~\input{fig_calculDr.pstex_t}
%\end{center}
\end{minipage}
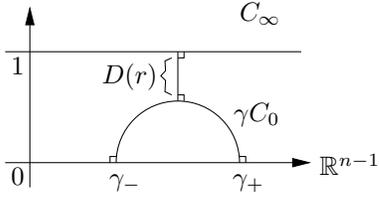
\begin{minipage}{8.9cm}
  ~~~ For every $\ga$ in $\Ga$, let $\ga_\pm$ be the fixed points of
  the hyperbolic element $\ga\ga_0\ga^{-1}$. If $r$ is the double
  coset of $\ga$ in $\Ga_\infty\backslash \Ga/\Ga_0$, and if
  $D(r)=d(C_\infty,\ga C_0)>0$, then by an easy computation in
  hyperbolic geometry
\begin{equation}\label{eq:calcDr}
D(r)=-\log \frac{1}{2}\|\ga_+-\ga_-\|\;,
\end{equation}
where $\|\cdot \|$ is the standard Euclidean norm on $\RR^{n-1}$.
\end{minipage}

\medskip Let $\R_{\Ga_0}$ be the set of fixed points of the conjugates
of $\ga_0$, endowed with its Fr\'echet filter.  For every $\alpha\in
\R_{\Ga_0}$, let $\alpha^*$ be the other endpoint of the translation
axis of a conjugate of $\ga_0$ containing $\alpha$ at infinity.  The
distance-like map $d_{C_\infty}$ coincides with the Hamenst\"adt
distance $d_{\infty,\partial \H_1}$ on the limit set $\Lambda\Gamma$.
In $\RR^{n-1}=\partial_\infty\hnr-\{\infty\}$, the Hamenst\"adt
distance $d_{\infty,\partial \H_1}$ coincides with the Euclidean
metric (see for instance \cite{HPMZ}).  For every $\xi$ in
$\Lambda_c\Ga-\R_{\Ga_0}$ and $\alpha\in \R_{\Ga_0}$, define
$\ell_\alpha+i\theta_\alpha$ to be the complex distance between the
oriented geodesic lines from $\infty$ to $\xi$ and from $\alpha^*$ to
$\alpha$.  Note that $\ell_{\alpha^*}=\ell_{\alpha}$ and
$\theta_{\alpha^*}= \theta_\alpha+\pi$. Then, we have by Equation
\eqref{eq:calcDr}, by Equation \eqref{eq:formcrosham} and by Lemma
\ref{lem:calccrossratiogeom}, respectively,
\begin{equation}\label{eq:calccstapproxconstcurv}
  c(\Ga_\infty\xi)=\liminf_{\alpha\in \R_{\Ga_0}}
\;2\,\frac{\|\xi-\alpha\|}{\|\alpha-\alpha^*\|}
=\liminf_{\alpha\in \R_{\Ga_0}} \;2\,e^{-[\infty,\alpha,\xi,\alpha^*]}
=\liminf_{\alpha\in \R_{\Ga_0}}\;\big(\cosh\ell_\alpha-\cos\theta_\alpha\big)\;.
\end{equation}
Furthermore, by definition, $\spd{\D}=\big\{\;c(\xi)\;:\; \xi\in
\Ga_\infty\backslash(\Lambda_c\Ga-\R_{\Ga_0})\big\}$.

\bigskip
\noindent{\bf Example 2: Spiraling around a closed geodesic in a
  complex hyperbolic manifold}

\medskip 
Let $n\geq 2$. The elements of $\CC^{n-1}$ are identified with their
coordinate column vectors and for every $w,w'$ in $\CC^{n-1}$, we
denote by $w^*w'$ their standard Hermitian product, where $w^*$ is the
conjugate transpose of $w$, and $|w|^2=w^*w$.

Let $\HH^n_\CC$ be the Siegel domain model of the complex hyperbolic
$n$-space. Its underlying manifold is
$$
\HH^n_\CC=\big\{(w_0,w)\in
\CC\times\CC^{n-1}\;:\; 2\,{\rm Re}\;w_0 -|w|^2>0\big\}\;.
$$
The complex hyperbolic distance $d_{\hnc}$ is defined by the
Riemannian metric
$$
ds^2=\frac{4}{(2\,{\rm Re}\;w_0 -|w|^2)^2}
\big((dw_0-dw^*\;w)(\overline{dw_0}-w^*\;dw)+
(2\,{\rm Re}\;w_0 -|w|^2)\;dw^*\;dw\big)
$$
(see for instance \cite[Sect.~4.1]{Gol}).  The complex hyperbolic
space has constant holomorphic sectional curvature $-1$, hence its
real sectional curvatures are bounded between $-1$ and $-\frac14$.
When we want to consider $\hnc$ as a CAT($-1$) space, we will use the
distance $d'_{\hnc}=\frac 12d_{\hnc}$.  The boundary at infinity of
$\hnc$ is
$$
\partial_\infty\HH^n_\CC=\{(w_0,w)\in
\CC\times\CC^{n-1}\;:\; 2\,{\rm Re}\;w_0 -|w|^2=0\}
\cup\{\infty\}\;.
$$
The horoballs centered at $\infty$ in $\HH^n_\CC$ are the subspaces
$$
\H_s=\{(w_0,w)\in
\CC\times\CC^{n-1}\;:\; 2\,{\rm Re}\;w_0 -|w|^2\geq s\}\;,
$$
for $s>0$. The submanifold $\{(w_0,w)\in\HH^n_\CC\;:\; w=0\}$, with
the induced Riemannian metric, is the right halfplane model of the
real hyperbolic plane with constant curvature $-1$, and it is totally
geodesic in $\HH^n_\CC$.  Hence the map $c_0:\RR\ra \HH^n_\CC$ defined
by $c_0:t\mapsto (e^{-t},0)$ is a unit speed geodesic line for
$d_{\hnc}$, starting from $\infty$, ending at $(0,0)\in \partial_\infty 
\HH^n_\CC$ and meeting the horosphere $\partial \H_2$ at time $t=0$.
In particular, the distance between two horospheres centered at
$\infty$ is
\begin{equation}\label{eq:calcdisthoro}
d_{\hnc}(\partial\H_s,\partial\H_{s'})=|\log(s'/s)|.
\end{equation}

Let $G_{\H_2}$ be the group of isometries of $\hnc$ preserving
(globally) $\H_2$. The {\it Cygan distance} $d_{\rm Cyg}$ (see for
instance \cite[page 160]{Gol}) is the unique distance on
$\partial_\infty\HH^n_\CC-\{\infty\}$ invariant under $G_{\H_2}$ such
that
$$
d_{\rm Cyg}((w_0,w),(0,0))= \sqrt{2|w_0|}.
$$ 
Similarly, we introduced in \cite[Lem.~6.1]{PP2} the {\it modified
  Cygan distance} $d'_{\rm Cyg}$, as the unique distance on
$\partial_\infty\HH^n_\CC -\{\infty\}$ invariant under $G_{\H_2}$ such
that
$$
d'_{\rm
  Cyg}((w_0,w),(0,0))= \sqrt{2|w_0|+|w|^2}.
$$ 
%Note that $d'_{\rm Cyg}\geq d_{\rm Cyg}$.

\medskip Let $\Ga$ be a discrete subgroup of isometries of $\HH_\CC^n$
with finite covolume. Assume that $\infty$ is a parabolic fixed point,
whose stabilizer in $\Ga$ we denote by $\Ga_\infty$, such that the
horoball $\H_2$ is precisely invariant under $\Ga_\infty$. Let $\ga_0$
be a hyperbolic element of $\Ga$, with translation axis $C_0$.  Let
$\Ga_0$ be the stabilizer of $C_0$ in $\Ga$.  The quadruple $\D=
(\HH^n_\CC,\Ga,\Ga_0,\H_2)$ satisfies the hypotheses of the beginning
of the section.

In the following result, we compute the associated map
$D:\Ga_\infty\bs\Ga/\Ga_0\ra\RR$ where $D([\ga])=d'_{\hnc}(\H_2,\ga
C_0)$.

\blemm\label{lem:calcDcomphyp}
 If $[\ga]\in\Ga_\infty\bs\Ga/\Ga_0$ and $D([\ga])>0$ then, with
$\ga_\pm$ the fixed points of the hyperbolic element
$\ga\ga_0\ga^{-1}$, we have
$$
D([\ga])=-\log\frac 12 \frac{d_{\rm Cyg}(\ga_-,\ga_+)^2}{d'_{\rm
    Cyg}(\ga_-,\ga_+)}\;. 
$$
\elemm

\dem Identify $\HH_\CC^n\cup\partial_\infty\HH_\CC^n$ with its image
in the projective space $\PP_n(\CC)$ by $(w_0,w)\mapsto [w_0:w:1]$ and
$\infty\mapsto [1:0:0]$. Note that $\ga C_0$ is the geodesic line
between $\ga_-$ and $\ga_+$. By invariance under $G_{\H_2}$, we may
assume that $\ga_-=(0,0)$ and $\ga_+= (w_0,w) \in \partial_\infty
\HH_\CC^n$ such that $w_0\ne 0$. The projective action of
$g=\left(\begin{array}{ccc} 1 & 0 & 0 \\ \frac{w}{w_0} & {\rm Id} & 0 \\
    \frac{1}{w_0} & \frac{w^*}{\overline{w_0}} & 1 \end{array} \right)
$ is an isometry of (the image of) the Siegel domain, fixing the point
$(0,0)$ of $\partial_\infty\HH_\CC^n$, and mapping $\infty$ to
$(w_0,w)$, hence sending the geodesic line between $(0,0)$ and
$\infty$ to the one between $(0,0)$ and $(w_0,w)$. Therefore the map
$\gamma_{w_0,w}:\RR\to\hnc$ defined by
$$
\gamma_{w_0,w}(t)=
     %\big[e^{-t}:e^{-t}\frac{w}{w_0}:\frac{w_0+e^{-t}}{w_0}\big]=
     %[w_0:w:1+w_0e^{-t}]=
\Big(\frac{w_0}{1+w_0e^{t}}\,,\frac{w}{1+w_0e^{t}}\Big)\;
%vaihdoin merkkia! 
$$ 
is a geodesic line with endpoints $(w_0,w)$ and $(0,0)$ in
$\partial_\infty\hnc$. The point $\gamma_{w_0,w}(t)$ belongs to the
horosphere $\partial\H_{s(t)}$, where
$$
s(t)=  
2\Re\big(\frac{w_0}{1+w_0e^{t}}\big)-\big|\frac{w}{1+w_0e^{t}}\big|^2
=  \frac{2\Re(w_0(1+\ov{w_0}e^t))-|w|^2}{|1+e^tw_0|^2}=
 \frac{2e^t|w_0|^2}{|1+e^tw_0|^2}.
$$
If $w_0=re^{i\phi}$ (in polar coordinates) and if $T=e^t$, then
$$
s(t)=\frac{2 Tr^2}{T^2r^2+2T r\cos\phi+1}\;.
$$
The map $t\mapsto s(t)$ reaches its maximum at $T=1/r$, that is at
$t=-\log|w_0|$, and its maximum value is
$$
s=\frac{|w_0|}{1+\frac{\Re\; w_0}{|w_0|}}=
%\frac{|w_0|^2}{|w_0|+\Re\; w_0}= 
\frac{|w_0|^2}{|w_0|+|w|^2/2}=
%\frac{2|w_0|^2}{2|w_0|+|w|^2}=
\frac 12 \frac{d_{\rm Cyg}((0,0),(w_0,w))^4}{d'_{\rm Cyg}((0,0),(w_0,w))^2}
\;.
$$
The result then follows from Equation \eqref{eq:calcdisthoro}, since
$\H_2$ and $\ga C_0$ are disjoint if and only if $s<2$. \cqfd

\medskip In $\hnc$ with its CAT($-1$) distance $d'_{\hnc}$, the
distance-like map $d_{\H_2}$ coincides (as seen in Section
\ref{sec:background}) with the Hamenst\"adt distance
$d_{\infty,\,\partial \H_2}$. Recall (see \cite[Prop.~3.12]{HPsurv})
that
\begin{equation}\label{eq:egalhamcyg}
d_{\infty,\,\partial \H_2}=\frac{1}{\sqrt{2}}\,d_{Cyg}\;.
\end{equation}
Let $\R_{\Ga_0}$ be the set of fixed points of the conjugates of
$\ga_0$, endowed with its Fr\'echet filter.  For every $\alpha\in
\R_{\Ga_0}$, let $\alpha^*$ be the other endpoint of the translation
axis of a conjugate of $\ga_0$ containing $\alpha$ at infinity.  We
therefore have
\begin{equation}\label{eq:spdCygan}
\spd{\D}=\big\{c(\Ga_\infty\xi)=\;\liminf_{\alpha\in \R_{\Ga_0}} 
\;\sqrt{2}\,\frac{d'_{\rm Cyg}(\alpha,\alpha^*)\;d_{\rm
    Cyg}(\xi,\alpha)}{d_{\rm 
    Cyg}(\alpha,\alpha^*)^2}\;:\;
\xi\in \Lambda_c\Ga-\R_{\Ga_0}\big\}\;.
\end{equation}

%%%%%%%%%%%%%%%%
%% Modified 27.12.09/Frédéric
%%%%%%%%%%%%%%%%

\section{The basic properties of the approximation spectra}
\label{sec:bounded}

Let $\D=(X,\Ga,\Ga_0,C_\infty)$ be a quadruple of data as defined in
Section \ref{sec:framework}.  In this section, we study the upper
bound of the approximation spectrum $\operatorname{Sp}(\D)\subset \RR$
of $\D$, and we give a closedness result for $\operatorname{Sp}(\D)$.

\subsection{The nontriviality of the approximation spectra}
\label{subsec:nontrivial}

A map $f:[0,+\infty[\;\ra\;]0,+\infty[$ is called {\it slowly varying}
if it is measurable and if there exist constants $B>0$ and $A\geq 1$
such that for every $x,y$ in $\RR_+$, if $|x-y|\leq B$, then $f(y)\leq
A\,f(x)$. Recall that this implies that $f$ is locally bounded, hence
it is locally integrable; also, if $\log f$ is Lipschitz, then $f$ is
slowly varying. 

Let $\epsilon$ be a positive real number, and let $f,g:[0,+\infty[\,
\to\,]0,+\infty[$.  A geodesic ray or line $\rho$ in $X$ will be
called {\it $(\epsilon,g)$-Liouville} (with respect to $\D$) if there
exist a sequence $(t_n)_{n\in \NN}$ of positive times converging to
$+\infty$ and a sequence $(\ga_n)_{n\in \NN}$ of elements of $\Ga$
such that $\rho(t)$ belongs to $\N_\epsilon (\ga_n C_0)$ for every $t$
in $[t_n,t_n+ g(t_n)]$.  A geodesic ray or line $\rho$ in $X$ such
that $\rho(+\infty)\notin \partial_\infty C_\infty$ will be called
{\it $f$-well approximated} (with respect to $\D$) if there exist
infinitely many $\ga$ in $\Ga/\Ga_0$ such that
$$
\wt d_\infty\big(\rho(+\infty),\ga\Lambda\Ga_0\big)\leq 
f\big(D([\ga])\big)\;e^{-D([\ga])}\;.
$$

The following result is proved in \cite[Lemma 5.2]{HPpre} (when
$C_\infty=C_0$, but the proof is the same).

\blemm \label{lem:euivliouborcan} %
Let $f:[0,+\infty[\;\ra\;]0,1[$ be slowly varying, and let $g:t\mapsto
-\log f(t)$.  Let $\epsilon>0$.  There exists $c=c(\epsilon,f)>0$ such
that for every geodesic ray or line $\rho$ in $X$ such that
$\rho(+\infty) \notin \partial_\infty C_\infty\cup \bigcup_{\ga\in\Ga}
\;\ga\; \partial_\infty C_0$, if $\rho$ is $(\epsilon,g)$-Liouville,
then $\rho$ is $(cf)$-well approximated, and conversely, if $\rho$ is
$(\frac{1}{c}f)$-well approximated, then $\rho$ is
$(\epsilon,g)$-Liouville.  
\cqfd \elemm

Our first result says in particular that $\{0\}\varsubsetneq
\spd{\D}$.  We refer to Section \ref{sec:hallray} for much stronger
results for particular cases of $\D$.

\bprop \label{prop:hurwitzpositiv} 
The approximation spectrum of $\D$ contains $0$ as a nonisolated
point, and hence the Hurwitz constant of $\D$ is positive.  
\eprop

The following consequence, amongst other similar ones, follows from
Equation \eqref{eq:divspec}.

\bcoro Let $M$ be a nonelementary complete connected Riemannian
manifold with sectional curvature at most $-1$ and dimension at least
$2$. Let $A_0$ be a closed geodesic in $M$, and let $A_\infty$ be a
closed codimension $0$ submanifold of $M$ with smooth connected
compact locally convex boundary, disjoint from $A_0$. Then the
spiraling spectrum $\operatorname{Sp}_{A_\infty,A_0}(M)$ around $A_0$
contains $0$ as a nonisolated point.  
\cqfd \ecoro

\noindent{\bf Proof of Proposition \ref{prop:hurwitzpositiv}. } Let us
first prove that there exists an element $\ga$ in $\Ga-\Ga_0$ such
that $d(C_0,\ga C_0)$ and $d(C_\infty,\ga C_0)$ are both as big as we
need.

By the lemmae \ref{lem:depthgrows} and \ref{lem:infinitedoublecoset},
there exists a nontrivial double class $[\ga_0]\in \Ga_\infty\bs \Ga/
\Ga_0$ such that $d(C_\infty,\ga_0C_0)$ is big. Since $\Ga_0$ contains
a hyperbolic element, there exists a hyperbolic element $\ga_1$ in
$\Ga$ whose attractive fixed point $(\ga_1)_+$ belongs to $\ga_0 
\partial_\infty C_0$, and in particular is not in $\partial_\infty
C_\infty$.  As $\ga_0\notin \Ga_0$ and $\Ga_0$ is almost malnormal, we
have $(\ga_1)_+\notin \partial_\infty C_0$ by Lemma
\ref{prop:equivmalnormal} (2).  Since $(\ga_1)_+\notin \partial_\infty
C_0\cup\partial_\infty C_\infty$, if $n$ is big enough, then
$\ga=\ga_1^n$ is an element in $\Ga-\Ga_0$ such that $d(C_0,\ga C_0)$
and $d(C_\infty,\ga C_0)$ are both big enough.

\medskip Now, let $\ga$ be as above.  Let $[p,q]$ be the shortest
segment between $C_0$ and $\ga C_0$, with $p\in C_0$. Let $\alpha$ be
a hyperbolic element in $\Ga_0$ with big translation length, and
$\beta=\ga\alpha\ga^{-1}$. Let $(k_n)_{n\in\NN}$ be a sequence of
positive integers. In particular, $L=d(p,q)$,
$L'_n=d(p,\alpha^{k_n} p)$ and $L''_n= d(q,\beta^{k_n} q)$ are big
(independently of $(k_n)_{n\in\NN}$). 

For every $n$ in $\NN$, define
$\ga_n=\beta^{k_1}\alpha^{k_1}\beta^{k_2}\alpha^{k_2}\dots
\beta^{k_n}\alpha^{k_n}$, so that $\ga_0={\rm id}$ and
$\ga_1=\beta^{k_1}\alpha^{k_1}$. Consider the piecewise geodesic ray
which is geodesic between the consecutive points
$$
p, q, \beta^{k_1} q, \beta^{k_1} p, \ga_1 p, \dots, \ga_n p, \ga_n q,
\ga_n\beta^{k_{n+1}} q, \ga_n\beta^{k_{n+1}} p, \ga_{n+1} p, \dots\;.
$$
Then if $A=\min\{L,L'_n,L''_n\;:\;n\in\NN\}$ is big enough, as the
comparison angles at the above points between the incoming and
outgoing segments are at least $\pi/2$ by convexity, this piecewise
geodesic ray is quasi-geodesic.

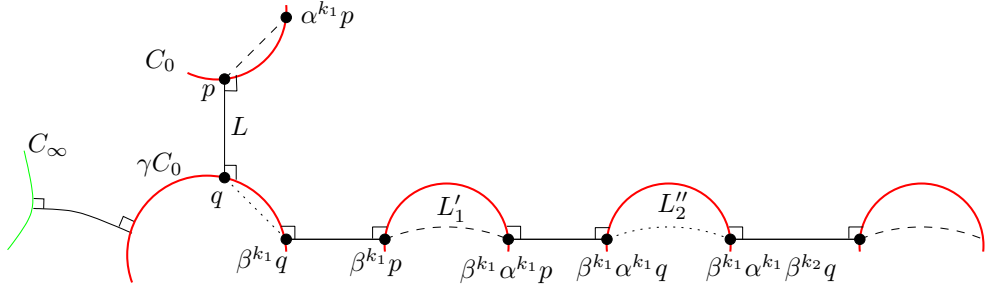
\begin{center}
\input{fig_quasigeod.pstex_t}
\end{center}

Hence, it stays at bounded distance (depending only on $A$) from a
geodesic ray $\rho^*$ starting from $p$. Note that, by convexity, the
segments $[\ga_n q, \ga_n\beta^{k_{n+1}} q]$ and  $[\ga_n\beta^{k_{n+1}}
p, \ga_{n+1} p]$ (which are long if the $L''_n$, $L'_n$ are big) are
contained in images under $\Ga$ of $C_0$. The point at infinity $\xi$
of $\rho^*$ is in particular a conical limit point (since $\Ga_0$ is
convex-cocompact, there are points in one orbit under $\Ga$ that
accumulate to $\xi$ while staying at bounded distance from $\rho^*$).
Up to taking the translation length of $\alpha$, and hence $A$,
big enough, the point $\xi$ belongs neither to $\partial_\infty
C_\infty$, nor to any $\ga'\partial_\infty C_0$ for $\ga'\in\Ga$
(otherwise, two copies of $C_0$ would be close for a too long time,
contradicting Lemma \ref{prop:equivmalnormal} (4)). Hence the
approximation constant $c(\Ga_\infty\xi)$ is well defined. 

In order to apply Lemma \ref{lem:euivliouborcan}, we fix $\epsilon>0$.
If the sequence $(k_n)_{n\in\NN}$ tends to $+\infty$, then the
geodesic ray $\rho^*$ spends longer and longer time in the images by
$\Ga$ of the $\epsilon$-neighbourhood of $C_0$.  Thus,
$c(\Ga_\infty\xi)$ is equal to $0$, by Lemma \ref{lem:euivliouborcan}.

\medskip To prove that $0$ is not isolated, take the sequence
$(k_n)_{n\in\NN}$ to be constant, with $k_1$ big compared with
$\kappa(\epsilon)$ (which has been defined in Lemma
\ref{prop:equivmalnormal} (4)), $L$ and the bounded distance between
$\rho^*$ and the above quasi-geodesic. In particular, $\rho^*$ is
$(\epsilon,g)$-Liouville for $g$ a constant map, having a big value if
$k_1$ is big. By Lemma \ref{prop:equivmalnormal} (4), since $\rho^*$
spends intervals of time of only bounded length outside
$\Ga\N_\epsilon C_0$, the geodesic ray $\rho^*$ is not
$(\epsilon,g')$-Liouville for $g'>g$ a big enough constant map.

By Lemma \ref{lem:euivliouborcan}, this implies that the approximation
constant of (the image modulo $\Ga_\infty$ of) $\xi$ is positive, and
small if $k_1$ is big. \cqfd

\medskip \rem Let us notice here that the approximation constants are
generically equal to $0$, hence that the nonvanishing of an
approximation constant is a quite rare behaviour. We will make this
explicit only in a particular case.

Assume that $X$ is a Riemannian manifold and $C_0$ a geodesic line.
For every $v\in \Ga\backslash T^1 X$, let $\xi_v\in
\Ga_\infty\backslash\partial_\infty X$ be the (orbit under $\Ga_\infty$
of the) endpoint of a geodesic line in $X$ whose tangent vector at the
origin maps to $v$ by the quotient by $\Ga$. (Several choices are
possible, but they will give the same approximation constant.) Let
$\mu$ be a (finite, positive, Borel) measure on $\Ga\backslash T^1 X$
invariant and ergodic under the quotient geodesic flow
$(\phi_t)_{t\in\RR}$. Assume that the support of $\mu$ contains the
orbit under $\Ga$ of the lift of $C_0$ to $ T^1 X$ by its unit tangent
vector, and that the (measurable) subset of unit vectors $v$ such that
$\xi_v\in \Ga_\infty\backslash\big((\partial_\infty X- \Lambda_c\Ga)
\cup \partial_\infty C_\infty\big)$ has measure $0$. For instance,
this is true if $\Ga$ has finite covolume, $\mu$ is the Liouville
measure and $C_\infty$ is a precisely invariant horoball, or if $\Ga$
is cocompact and $\mu$ is the maximal entropy measure, and $C_\infty$
is the translation axis of a hyperbolic element. The ergodicity
assumption implies that $\{\phi_t v\} _{t\in\RR^+}$ is dense in the
support of $\mu$ for almost every $v$. Recall that if two unit tangent
vectors are very close, then the geodesic lines they define are close
for a long time. Hence for $\mu$-almost every $v$, we have $\xi_v\in
\Ga_\infty \backslash (\Lambda_c\Ga-\partial_\infty C_\infty)$ and
$c(\xi_v)=0$ by Lemma \ref{lem:euivliouborcan}.

\subsection{The boundedness of the approximation spectra}
\label{subsec:bounded}

If $\Ga$ is geometrically finite (see for instance \cite{Bow}), then
there exists a $\Ga$-equivariant family $\H$ of horoballs centered at
the parabolic fixed points of $\Ga$, with pairwise disjoint interiors.
There are many possible choices for such an $\H$ (though only one
maximal one if $\Ga$ has only one orbit of parabolic fixed points). In
the computations of Section \ref{subsec:upperbounds}, we will choose
natural ones. We call $X_0=\C\Ga-\bigcup\H$ the {\em thick part} of
$\C\Ga$.  Clearly, $X_0$ is $\Ga$-invariant, and $\Ga$ acts
isometrically on it. We call $\Ga\backslash X_0$ the {\em thick convex
  core} of $\Ga\backslash X$.

\medskip 
The next result gives a sufficient condition for the Hurwitz constant
of $\D$ to be finite.  In particular, this condition is satisfied when
$X$ is a Riemannian manifold and $\Ga$ has finite covolume.  Recall
that the Hurwitz constant of $\D$ is
$$
K_\D=\sup \,\spd{\D}\;\in[0,+\infty]\;.
$$

\btheo\label{theo:Dirichlet} If $\Ga$ is geometrically finite, then
$\spd{\D}$ is bounded, hence $0<K_\D<\infty$.  \etheo

The following consequence, amongst other similar ones, follows from
Equation \eqref{eq:divspec}.

\bcoro Let $M$ be a geometrically finite complete connected Riemannian
manifold with sectional curvature at most $-1$ and dimension at least
$2$. Let $A_0$ be a closed geodesic in $M$, and let $A_\infty$ be a
closed codimension $0$ submanifold of $M$ with smooth connected
compact locally convex boundary, disjoint from $A_0$. Then the
spiraling spectrum $\operatorname{Sp}_{A_\infty,A_0}(M)$ around $A_0$
is bounded.  \cqfd \ecoro

\noindent{\bf Proof of Theorem \ref{theo:Dirichlet}. }  Let $\H$ be as
above. Since $\Lambda\Ga_0$ contains at least two points, $C_0$ is not
contained in any element of $\H$, hence $C_0$ intersects $X_0$. Since
$\Ga$ is geometrically finite, the diameter $\Delta$ of the quotient
metric space $\Ga\bs X_0$ is finite.  For every $\xi\in \Lambda_c\Ga
-(\partial_\infty C_\infty\cup\bigcup_{\ga\in \Ga} \ga\partial_\infty
C_0)$, let $\rho_\xi$ be a geodesic ray starting from the closest
point to $\xi$ on $C_\infty$ and converging to $\xi$. As $\xi$ is a
conical limit point, there exists a sequence of positive times
$(t_n)_{n\in\NN}$ converging to $+\infty$ such that $\rho_\xi(t_n)\in
X_0$ for every $n$. Hence, there exists a sequence of elements
$(\ga_n)_{n\in\NN}$ such that for every $n$ in $\NN$,
$$
d(\rho_\xi(t_n),\ga_nC_0)\leq \Delta\;.
$$

For $n$ big enough, the distance between the convex subsets $C_\infty$
and $\ga_n C_0$ is big. Indeed, if $d(C_\infty,\ga_{n_k} C_0)$ is
bounded for some subsequence $(n_k)_{k\in\NN}$ tending to $+\infty$,
then by Lemma \ref{lem:depthgrows} and up to extracting a subsequence,
the double cosets $[\ga_{n_k}]\in\Ga_\infty\bs\Ga/\Ga_0$ are
constant. Since $\ga_{n_k}C_0$ contains a point whose closest point on
$C_\infty$ is at bounded distance (at most $\Delta$) from the point
$\rho_\xi(0)$, up to extracting a subsequence and up to multiplying
$\ga_{n_k}$ on the right by an element of $\Ga_0$, we may assume that
$(\ga_{n_k})_{k\in\NN}$ is constant. Since $\rho_\xi$ converges to
$\xi$, the construction of $(\ga_n)_{n\in\NN}$ implies that $\xi$
belongs to the closed subset $\ga_{n_0}\partial_\infty C_0$, a
contradiction.

\smallskip

\noindent
\begin{minipage}{8.6cm}
  If $\rho_\xi$ does not meet $\ga_nC_0$, we denote by $[p_n,q_n]$ the
  shortest segment between $\rho_\xi$ and $\ga_n C_0$, with $p_n$ in
  $\rho_\xi$ (which exists since $\xi\notin \bigcup_{n\in\NN} \ga_n
  \partial_\infty C_0$). Otherwise, define $p_n=q_n$ to be the first
  intersection point of $\rho_\xi$ with $\ga_nC_0$. In particular,
  $d(p_n,q_n)\leq \Delta$.  Let $h_n$ be the point of $\ga_n C_0$ the
  closest to $C_\infty$, and $h'_n$ its closest point on $C_\infty$.
  Since $\Ga_0$ is convex-cocompact, there exists a constant $c_1\geq
  0$ (depending only on $\Ga_0$) such that the distance between $q_n$
  and a geodesic ray starting from $h_n$ and staying inside $\ga_nC_0$
  is at most $c_1$, for every $n$ in $\NN$.  Denote by $\xi_n$ the
  point at infinity of this geodesic ray, which does not belong to
  $\partial_\infty C_\infty$. Let $\rho_n$ be the geodesic ray
  starting from the closest point to $\xi_n$ in $C_\infty$ and
  converging to $\xi_n$.
\end{minipage}
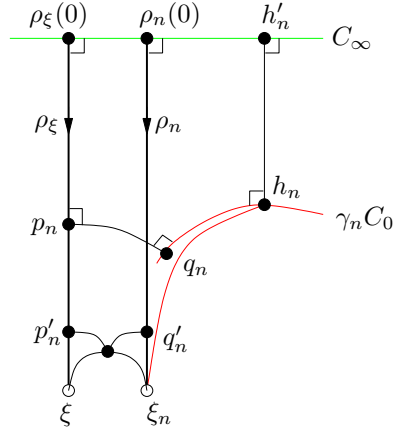
\begin{minipage}{6.3cm}
\begin{center}
\input{fig_dirichlet_F.pstex_t}
\end{center}
\end{minipage}

\medskip By the distance formulae of the hyperbolic comparison
quadrilateral with vertices corresponding to $\xi_n, \rho_n(0), h'_n$
and $h_n$, with one comparison angle $0$ and the three others at least
$\pi/2$ (see \cite[7.17]{Bea}), if $n$ (and hence $d(h_n ,C_\infty)$)
is big enough, then the distance between $\rho_n(0)$ and $h'_n$ is at
most $1$. Therefore $h_n$ is at distance at most a universal constant
$c_2$ from $\rho_n$.  Since $d(q_n,[h_n,\xi_n[)\leq c_1$, the point
$q_n$ is hence by convexity at distance at most $c_1+c_2$ from a point
$x_n$ of $\rho_n$, so that
$$
d(p_n,\rho_n)\leq d(p_n,x_n)\leq d(p_n,q_n)+ d(q_n,x_n)\leq
\Delta+c_1+c_2\;.
$$  
In particular, $d(\rho_n(0),x_n)$ tends to $+\infty$ as $n\ra+
\infty$.  By the distance formulae of the hyperbolic comparison
quadrilateral with vertices corresponding to $\rho_\xi(0),\rho_n(0),
x_n, p_n$ (see \cite[7.17]{Bea}), with two comparison angles at least
$\pi/2$ and the length of the opposite segment, that is $d(p_n,x_n)$,
bounded by a constant, if $n$ is big enough, then $\rho_\xi(0)$ is at
distance at most $1$ from $\rho_n(0)$.

Let $r_n$ be the class of $\ga_n$ in $\Ga_\infty\backslash \Ga/\Ga_0$,
so that $D(r_n)=d(h_n,C_\infty)$. Using the triangle inequality, we
have
$$ 
d(p_n,\rho_\xi(0))\ge d(q_n,C_\infty)-d(q_n,p_n)\geq
d(\ga_nC_0,C_\infty)-d(q_n,p_n)\geq D(r_n)-\Delta\;.
$$ 
Since $d(\rho_\xi(0),\rho_n(0))\leq 1$, there exists a universal
constant $c_3\geq 0$ such that some point on the geodesic line between
$\xi$ and $\xi_n$ is at distance at most $c_3$ from both a point
$p'_n$ of $\rho_\xi$ and a point $q'_n$ of $\rho_n$.  In particular,
$d(p'_n,q'_n)\leq 2\,c_3$. Since $p_n$ is at bounded distance from
$\rho_n$, there exists (by a geometric argument of geodesic triangles)
a constant $c_4\geq 0$ (which depends only on $\Delta$) such that
$$
d(p'_n,\rho_\xi(0))\geq d(p_n,\rho_\xi(0))-c_4\;.
$$
By the definition of the distance-like map $d_{C_\infty}$ and by the
triangle inequality, we have
$$
d_{C_\infty}(\xi,\xi_n)\leq 
e^{-\frac{1}{2}(d(\rho_\xi(0),\,p'_n)-c_3+d(\rho_n(0),\,q'_n)-c_3)}\;.
$$
By the triangle inequality, $d(\rho_n(0),q'_n)\geq d(\rho_\xi(0),p'_n)
- d(\rho_\xi(0),\rho_n(0)) - d(p'_n,q'_n)$. By the previous
inequalities, we hence have
$$
d_{C_\infty}(\xi,\xi_n)\leq e^{-d(\rho_\xi(0),\,p'_n)+c_3+
  \frac{d(p'_n,\,q'_n)}2+ \frac{d(\rho_\xi(0),\,\rho_n(0))}{2}} \leq
e^{\Delta+2c_3+c_4+\frac{1}{2}}\;e^{-D(r_n)}\;.
$$
The approximation constant $c(\xi)$ of $\xi$, for every $\xi$ in
$\Ga_\infty\bs(\Lambda_c\Ga-(\partial_\infty C_\infty
\cup\bigcup_{\ga\in \Ga}\ga\Lambda\Ga_0))$, is hence at most
$e^{\Delta+2c_3+c_4+\frac{1}{2}}$, which proves the result.  \cqfd

\medskip \rem Note that this upper bound depends only on
$(X,\Ga,C_\infty)$, but not on $C_0$.

\medskip
In special cases, the proof of Theorem \ref{theo:Dirichlet}
may be improved to give a simple explicit constant.

\bprop\label{prop:specialDirichlet} If $\Ga$ is geometrically finite,
if $\Delta$ is the diameter of the thick convex core of $\Ga\bs\Ga$,
if $C_\infty$ is a point in $X$ or a horoball in $X$, and if $C_0$ is
a geodesic line, then
$$
K_\D\le (1+\sqrt{2})e^{\Delta}.
$$
\eprop

The following consequence follows from Equation \eqref{eq:divspec}.

\bcoro \label{coro:bornsupunifsiralspec}
Let $M$ be a geometrically finite complete connected Riemannian
manifold with sectional curvature at most $-1$ and dimension at least
$2$, and let $\Delta$ be the diameter of the thick convex core of
$M$. Let $A_0$ be a closed geodesic in $M$, and let $A_\infty$ be
either a ball or a Margulis neighbourhood of a cusp of $M$. Then the
spiraling spectrum $\operatorname{Sp}_{A_\infty,A_0}(M)$ around $A_0$
is contained in $[0,(1+\sqrt{2})e^{\Delta}]$.  
\cqfd \ecoro

\noindent{\bf Proof of Proposition \ref{prop:specialDirichlet}. }
Assume first that $C_\infty=\{x_\infty\}$ with $x_\infty\in X$.  For
every $\xi$ belonging to $\Lambda_c\Ga-(\partial_\infty C_\infty
\cup\bigcup_{\ga\in \Ga}\ga\,\partial_\infty C_0)$, we define
$\rho_\xi,\ga_n,p_n,q_n,h_n$ as in the beginning of the proof of
Theorem \ref{theo:Dirichlet}, so that $\rho_\xi(0)=x_\infty$. Let
$\xi_n$ be an endpoint of the geodesic line $C_0$ such that $q_n\in
[h_n,\xi_n[$\,. Let $\rho_n$ be the geodesic ray from $x_\infty$ to
$\xi_n$. Let $r_n=[\ga_n]\in\Ga_\infty\bs\Ga/\Ga_0$, so that
$d(x_\infty,q_n)\geq D(r_n)$.  Let $z_n$ be the closest point on
$\rho_n$ to $q_n$, which satisfies
$$ 
d(z_n,q_n)\leq \delta\;,
$$ 
with $\delta=\log (1+\sqrt{2})$, by looking at the comparison triangle
of the geodesic triangle with vertices $x_\infty, \xi_n, h_n$ (see the
picture below).

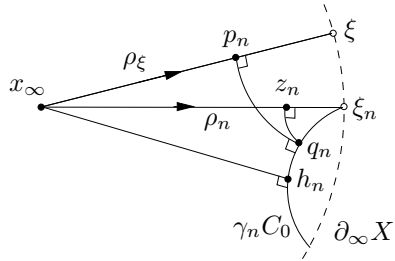
\begin{center}
\input{fig_dirichlet.pstex_t}
\end{center}

By the triangle inequality, for all $t$ big enough,
\begin{align*}
& d(\rho_n(t),x_\infty)+d(\rho_\xi(t),x_\infty)-
d(\rho_n(t),\rho_\xi(t))\\ 
& \geq \big(d(\rho_n(t),z_n)+d(z_n,x_\infty)\big)+
\big(d(\rho_\xi(t),p_n)+d(p_n,x_\infty)\big)-\\
&\quad\quad\big(d(\rho_n(t),z_n)+d(z_n,p_n)+d(p_n,\rho_\xi(t))\big)\\
& = d(x_\infty,p_n)+d(x_\infty,z_n)-d(p_n,z_n)\\ 
& \geq  d(x_\infty,q_n)-d(q_n,p_n)+d(x_\infty,q_n)-d(q_n,z_n)-
\big(d(p_n,q_n)+d(q_n,z_n)\big)\\ 
& \geq 2D(r_n)- 2\Delta-2\delta\;.
\end{align*}
Hence,
$$
d_{x_\infty}(\xi,\xi_n)=\lim_{t\ra+\infty}e^{-\frac12
\big(d(\rho_n(t),\,x_\infty)+d(\rho_\xi(t),\,x_\infty)-
d(\rho_n(t),\,\rho_\xi(t))\big)}\le
e^{-D(r_n)+\Delta+\delta},
$$
which proves that $c(\Ga_\infty\xi)\leq e^{\Delta+\delta}$. 

\medskip If $C_\infty$ is a horoball with point at infinity
$\xi_\infty$, the proof is similar, by replacing the geodesic rays
starting from $x_\infty$ by geodesic lines starting from $\xi_\infty$
and meeting $\partial C_\infty$ at time $0$, and using the fact
that $d(\rho_\xi(0),\rho_n(0))$ tends to $0$ as $n\ra+\infty$.  
\cqfd

\medskip Theorem \ref{theo:mainintroun} in the introduction follows
from Corollary \ref{coro:bornsupunifsiralspec}.

\subsection{On the closedness of the approximation spectra}
\label{subsec:closedness}

In this subsection, we prove that in the constant curvature manifold
case, the spiraling spectrum around a closed geodesic is closed.

\btheo \label{theo:closedspec} 
Let $\D=(X,\Ga,\Ga_0,C_\infty)$ be a quadruple of data such that $X$
is the real hyperbolic $n$-space, $\Ga$ is geometrically finite,
$C_0=\C\Ga_0$ is a geodesic line, and $C_\infty$ is a horoball. Then
$\spd{\D}$ is equal to the closure in $\RR$ of the set of the
approximation constants of the (orbits under $\Ga_\infty$ of the) fixed
points of the hyperbolic elements of $\Ga$ (that are not conjugated to
elements of $\Ga_0$).  
\etheo

In particular, the approximation spectrum $\spd{\D}$ is closed, the
Hurwitz constant of $\D$ is the maximum of $\spd{\D}$, and the
approximation constants of the (orbits under $\Ga_\infty$ of the)
hyperbolic fixed points of $\Ga$ are dense in $\spd{\D}$.

\medskip \rem The result is still true if $C_0$ is any totally
geodesic subspace of dimension at least $1$ and at most $n-1$; the
adaptation of the proof below is left to the reader.

\medskip The following consequence follows from Equation
\eqref{eq:divspec}, and proves the first claim in Theorem
\ref{theo:mainintrodeux} in the Introduction.

\bcoro Let $M$ be a geometrically finite complete connected Riemannian
manifold with constant sectional curvature $-1$ and dimension at least
$2$. Let $A_0$ be a closed geodesic in $M$, and let $A_\infty$ be a
Margulis neighbourhood of a cusp of $M$. Then the spiraling spectrum
$\operatorname{Sp}_{A_\infty,A_0}(M)$ around $A_0$ is closed, and is
equal to the closure of the set of spiraling constants of the geodesic
lines spiraling around closed geodesics distinct from $A_0$.  \cqfd
\ecoro

\noindent{\bf Proof of Theorem \ref{theo:closedspec}. } Let $\L_0$ be
the set of images under $\Ga$ of the two oriented geodesics defined by
$C_0$. Let $\tilde v$ be an element of $T^1X$, and $\tilde x$ be its
base point. For every $C$ in $\L_0$, define $p_C$ to be the point of
$C$ the closest to $\tilde x$, which depends continuously on $\tilde
v$; define $\theta_C$ to be the angle at $p_C$ between the parallel
transport along $[\tilde x,p_C]$ of $\tilde v$ and $C$, which depends
continuously on $\tilde v$.  Let
$$
\tilde f(\tilde v)=\inf_{C\in\L_0}\;\cosh d(\tilde x,C)-\cos \theta_C\;.
$$
Since $\L_0$ is locally finite and is preserved by $\Ga$, the lower
bound defining $\tilde f$ is locally a minimum over a finite set. 
Thus, the map $\tilde f:T^1X\ra \RR$ is continuous and invariant under
$\Ga$, and it defines a continuous map $f:\Ga\bs T^1X\ra \RR$.  As the
image of $C_0$ in $\Ga\bs X$ is compact and $\Ga\bs X$ is a proper
metric space, the distance to this image is a proper map on $\Ga\bs
X$. Hence $f$, which is at least $\cosh (0) -1= 0$, is
proper. Let $(\phi^t:T^1X\ra T^1X)_{t\in\RR}$ be the geodesic flow of
$X$, and denote again by $(\phi^t)_{t\in\RR}$ its quotient flow under
$\Ga$.

We will use the following result of F.~Maucourant \cite[Theo.~2
(2)]{Mau}, whose main tool is Anosov's closing lemma (and which builds on a
partial result of \cite{HPMZ}). The result extends to our orbifold
case.

\btheo \label{theo:maucourant} Let $V$ be a complete Riemannian
manifold with sectional curvature at most $-1$, let
$(\phi^t)_{t\in\RR}$ be its geodesic flow, and let $J_0$ be the subset
of $T^1V$ which consists of periodic unit tangent vectors. If
$f:T^1V\ra \RR$ is a proper continuous map, then
$$
\RR\cap\{\;\liminf_{t\ra+\infty} f(\phi^t v)\;:\;v\in T^1V\}=
\overline{\{\;\min_{t\in\RR} f(\phi^t v)\;:\;v\in J_0\}}
\;.\;\;\;\mbox{\cqfd}
$$  
\etheo

Assume that $X$ is the upper halfspace model of $\HH^n_\RR$, and that
$C_\infty$ is centered at $\infty$. By the assumptions on the data
$\D$, we are in the situation of Example 1 of Section
\ref{sec:framework}. In the following, we use the notation of that example.
Let $\xi$ be the endpoint of the geodesic line
defined by $\tilde v$, and note that this geodesic line is asymptotic
to the geodesic line from $\infty$ to $\xi$.  There are three cases to
consider:
\begin{enumerate}
\item If $\xi\in\Lambda_c\Ga-\R_{\Gamma_0}$, then it
follows from the definition of $\tilde f$ and from Equation
\eqref{eq:calccstapproxconstcurv} that $\liminf_{t\ra+\infty} \tilde
f(\phi^t \tilde v)=c(\Ga_\infty\xi)$.

\item If $\xi\in\R_{\Gamma_0}$, then $\liminf_{t\ra+\infty} \tilde f(\phi^t
\tilde v)=0$, and we have already seen (in Proposition
\ref{prop:hurwitzpositiv}) that $0$ belongs to $\spd{\D}$.

\item If $\xi\in \partial_\infty X-\Lambda_c\Ga$, then since $\Ga$ is
geometrically finite, either $\xi$ does not belong to the limit set,
or $\xi$ is a parabolic fixed point. In both cases, since $f$ is
proper, we have $\liminf_{t\ra+\infty} \tilde f(\phi^t \tilde
v)=+\infty$, which is not in $\RR$.
\end{enumerate}
These observations imply that
$\RR\cap\{\;\liminf_{t\ra+\infty} f(\phi^t v)\;:\;v\in \Ga\bs T^1X\}$
is contained in 
$\spd{\D}=\{c(\Ga_\infty\xi)\;:\;\xi\in\Lambda_c\Ga-\R_{\Gamma_0}\}$.

By considering a vertical unit tangent vector $\tilde v$ ending at a
given $\xi\in\Lambda_c\Ga-\R_{\Gamma_0}$, the opposite inclusion also
holds.  If $J'_0$ is the subset of vectors in $J_0$ that are not the
image in $\Ga\bs T^1X$ of unit tangent vectors to $C_0$, then the set $A$
of the approximation constants of the (orbits under $\Ga_\infty$ of
the) points of $\partial_\infty X$ fixed by hyperbolic elements not
conjugated to elements of $\Ga_0$, is equal to $\{\;\inf_{t\in\RR}
f(\phi^t v)\;:\;v\in J'_0\} $.  Furthermore, the approximation
constant of (the orbit under $\Ga_\infty$ of) a point of
$\partial_\infty X$ fixed by a hyperbolic element conjugated to an
element of $\Ga_0$, is equal to $0$. Hence, by Theorem
\ref{theo:maucourant}, we have $\spd{\D}=\overline{A\cup\{0\}}$. Since
$0$ is not isolated in $\spd{\D}$ (see Theorem
\ref{prop:hurwitzpositiv}), this implies that $\spd{\D}=\overline{A}$.
This proves Theorem \ref{theo:closedspec}. \cqfd

\subsection{Some upper bounds on the approximation spectra}
\label{subsec:upperbounds} 

We give estimates of the Hurwitz constants of data $\D=(X,\Ga,\Ga_0,
C_\infty)$ in a number of arithmetically defined cases. The estimates
are not likely to be very sharp, except for Proposition
\ref{prop:calchurwnbdor}. In the following five examples, $X$ is
$\HH_\RR^2$, $\HH_\RR^2$, $\HH_\RR^3$, $\HH_\RR^5$ or $\HH_\CC^2$
respectively. These dimensions are chosen with number theoretical
applications in mind, see Section \ref{sec:appdioapp}.  In the first
four examples, the convex set $C_\infty$ is the horoball centered at
infinity consisting of the points with Euclidean height at least $1$
in the upper half space model of $\hnr$.  The group $\Ga$ is specified
in each example, and the subgroup $\Ga_0$ is the stabilizer in $\Ga$
of any geodesic in $X$ whose quotient in $\Ga\backslash X$ is compact.

The following classical fact is used repeatedly in the examples: For
every $\alpha\leq \pi/2$, the distance $\ell$ in $\hdr$ between the
points of angle $\alpha$ and of angle $\pi/2$ with respect to the real
line $\partial_\infty\hdr-\{\infty\}$ on any Euclidean circle centered
at a point in $\partial_\infty\hdr-\{\infty\}$ is
\begin{equation} \label{eq:lengthcomput}
\ell=\log \cot \frac{\alpha}{2}=
\log \;\frac{\sqrt{1+\tan^2\alpha}\,+1}{\tan \alpha}\;.
\end{equation}
This equation is used to compute distances between points in
isometrically embedded copies of $\hdr$ in $\hnr$ and $\hdc$.

\medskip%\noindent
(1) Let $\Ga=\PSLZ$. It is well known that the hyperbolic triangle $F$
in $\hdr$ with vertices at $\infty,e^{i\frac\pi3}$ and
$e^{2i\frac\pi3}$ is a fundamental polygon for $\Ga=\PSLZ$. The
horoball $C_\infty$ covers the maximal Margulis neighbourhood $U$ of
the (only) cusp of $M= \Ga\backslash\hdr$.  The compact set
$K=\overline{F-C_\infty}$ covers the complement of $U$ in $M$.  The
symmetries imply that the diameter $\Delta$ of $M-U$ equals the
distance between the cone points of $M$ with angles $\pi$ and
$2\pi/3$. By Equation \eqref{eq:lengthcomput},
$$
\Delta=d(i,e^{i\frac\pi 3})=\frac 12 \log 3\approx 0.55\;.
$$ 
By Proposition \ref{prop:specialDirichlet}, we thus have $\spd{\D}
\subset [0,(1+\sqrt 2)\sqrt 3]\subset[0,4.19]$. Note that this upper
bound is uniform amongst the subgroups $\Ga_0$.

\medskip
Let us give an exact computation of the Hurwitz constant in a
particular case. Notice that the second assertion of the result below
shows a different behaviour than the classical Lagrange spectrum.

\bprop \label{prop:calchurwnbdor} Let $\Ga_0$ be the cyclic subgroup
of $\Ga=\PSLZ$ generated by $\ga_1=\pm\begin{pmatrix}2 &1\\1 & 1
\end{pmatrix}$, and let $\D=(\hdr,\Ga,\Ga_0,C_\infty)$. Then $K_\D=
1-1/\sqrt 5$, and $K_\D$ is not isolated in the approximation spectrum
$\spd{\D}$.  
\eprop

\dem The element $\ga_1$ is hyperbolic, and its translation axis $L_1$
is the geodesic line in $\hdr$ with endpoints at $(1\pm\sqrt 5)/2$
(see the picture below).  The translation length $\ell_1$ of $\ga_1$
satisfies $2\cosh(\ell_1/2) = 3$ (see \cite[page 173]{Bea}), and the
translates of $L_1$ intersect in pairs at the orbit of $i$, and form a
net (covering $\HH^2_\RR-\Ga C_\infty$) of equilateral triangles (the
images under $\Ga$ of the triangle with vertices $i,i+1, \frac{i+1}{2}
$) as in the figure below.  The edges of the triangles have length, by
Equation \eqref{eq:lengthcomput}, equal to $d(i,i+1)=\arcosh(3/2) =
\ell_1/2$, and the (interior) angles $\theta\in [0,\frac{\pi}{2}]$ of
the triangles satisfy $\cos\theta=3/5$. These facts are easily seen by
considering the $6$-fold cover of $M$ by the modular torus $M'$ which
is the quotient of $\hdr$ by the commutator subgroup of $\Ga=\PSLZ$
(see (2) below). Notice that $\ga_1$ and its translates by $z\mapsto
z\pm 1$ are lifts of the three shortest periodic geodesics of $M'$
which intersect at the three Weierstrass points of $M'$ (see for
instance \cite[Theo.~2]{Sch}).

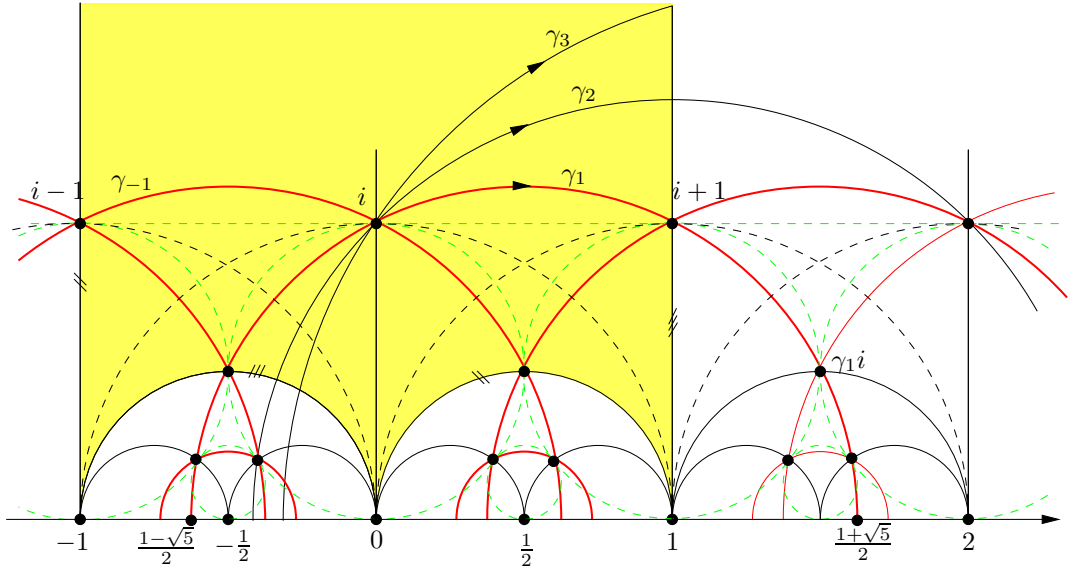
\begin{center}
\input{fig_modulartorus2.pstex_t}
\end{center}

Any geodesic line that connects $\infty$ with $\xi\in\RR-\QQ$
intersects infinitely many triangles. Let $T$ be one of the triangles,
and let $\ga$ be a geodesic line that intersects the interior of
$T$. The points of intersection of $\ga$ with $T$ and one of the
vertices of $T$ determine a triangle with angles $\theta,\phi_1,
\phi_2$. The supremum of $\min\{\cos \phi_1,\cos \phi_2\}$ over all
nondegenerate hyperbolic triangles with angles $\theta,\phi_1, \phi_2$
is obtained, by symmetry, when $\phi_1=\phi_2$ and when the triangles
become small, that is when they converge, after renormalization, to
the Euclidean triangle with angles $\theta,\phi_1=\phi_2$. Hence the
above supremum is $\cos(\frac{\pi-\theta}{2})=\sin\frac{\theta}{2}=
\frac{1} {\sqrt{5}}$, and it is not attained, by the Gauss-Bonnet
formula. Thus, the geodesic line from $\infty$ to $\xi$ cuts a sequence
of pairwise distinct $\Ga$-translates of $C_0$, the cosine of the
angle at each intersection point being at least $\frac{1}{\sqrt{5}}$.
By Equation \eqref{eq:calccstapproxconstcurv}, this implies that for
any $\xi\in\RR-\QQ$, we have $c(\xi)\le 1-1/\sqrt 5$.

For any nonzero integer $n$, consider the hyperbolic element  
$\ga_n=\pm\begin{pmatrix} n^2+1  &n\\n &1 \end{pmatrix}\in\Ga$.
The fixed points of $\ga_n$ are
    $\frac n2\pm\sqrt{(\frac n2)^2+1}$. Thus, the axis of $\ga_n$ is the intersection
    with the upper half plane of the Euclidean circle of center $n/2$ and
    radius  $\sqrt{(\frac n2)^2+1}$, which passes through the points $i$ and
  $n+i$.
The
translation distance of $\ga_n$, which is $2\,\arcosh(n^2/2+1)$ by
\cite[page 173]{Bea}, is twice the distance between the points $i$ and
  $n+i$, by Equation \eqref{eq:lengthcomput}. 
Thus, the
translation axis of $\ga_n$ intersects, at each $\Ga$-image of $i$ on
it, exactly two $\Ga$-translates of the axis of $\ga_1$, and always at
the same angle in absolute value, with alternating signs. As
$n\to\infty$, the smallest of the two positive angles approaches
(while strictly increasing) the angle $\theta'$ between the (oriented)
axis of $\ga_1$ and the (upward oriented) imaginary axis at $i$, which
satisfies $\cos\theta'=1/\sqrt 5$ . Thus, by Equation
\eqref{eq:calccstapproxconstcurv}, the approximation constants of the
lines from $\infty$ to the fixed points of $\ga_n$ converge to
$1-1/\sqrt 5$ (while being different).  The result follows.  
\cqfd

\medskip
(2) Let $\Ga$ be the commutator subgroup of $\PSLZ$. It is well known
(see for instance \cite{Sch}) that $\Ga$ is a torsion-free subgroup of
index $6$ in $\PSLZ$, and that the quotient $\Ga\bs\hdr$ is a
punctured torus, called {\em the modular torus}.

\medskip 
For every $k\in\NN$, let $H_k$ be the horoball centered at $k$ with
Euclidean height one. It is well known (see for instance
\cite{Coh,Sch}) that the modular torus is isometric to the quotient of
the ideal hyperbolic square $P$ with vertices $\infty, -1, 0, 1$ by
the gluing of the opposite faces of $P$, such that the horoball
$C_\infty$ maps by the two gluings to the horoballs $H_{-1}$ and
$H_1$. In particular, $C_\infty$ covers the maximal Margulis
neighbourhood $U$ of the cusp of $M=\Ga\backslash\hdr$.

Let $T'$ be the closure of the relatively compact component of
$\hdr-(C_\infty\cup H_{-1}\cup H_0)$. Then the closure of $M-U$ is the
union of the triangles with horocyclic sides $T'$ and $T'+1$, glued
along their vertices. The diameter of $T'$ for the induced distance of
$\hdr$ is, by a convexity argument, equal to
$d(i,i+1)=\arcosh(3/2)$. Any point of $T'+1$ is at distance at most
$d(i,e^{i\frac{\pi}{3}})=\frac{\log 3}{2}$ (by Equation
\eqref{eq:lengthcomput}) from one vertex of $T'+1$. Therefore
$$
\Delta\leq \frac{\log 3}{2}+\arcosh(3/2)\;.
$$
By Proposition \ref{prop:specialDirichlet}, we thus have
$\spd{\D}\subset[0,(1+\sqrt 2)e^\Delta]\subset[0,10.95]$.

\medskip%\noindent
(3) Let $m$ be a positive squarefree integer, and let $\Ga$ be the
{\it Bianchi group} ${\rm PSL}_2(\O_{-m})$, where $\O_{-m}$ is the
ring of integers of $\QQ(i\sqrt m)$.  All Bianchi groups contain the
transformation $z\mapsto z+1$, and thus, the interior of the horoball
$C_\infty$ is precisely invariant, by Shimizu's Lemma. Since
$\iota:z\mapsto -\frac{1}{z}$ also belongs to $\Ga$ and since the
horoballs $C_\infty$ and $\iota C_\infty$ are tangent, the horoball
$C_\infty$ covers the maximal Margulis neighbourhood of the cusp of
$M=\Ga\backslash\htr$. Fundamental domains for the Bianchi groups have
been determined in \cite{Bia,Swa,EGM} and we will use the tables of
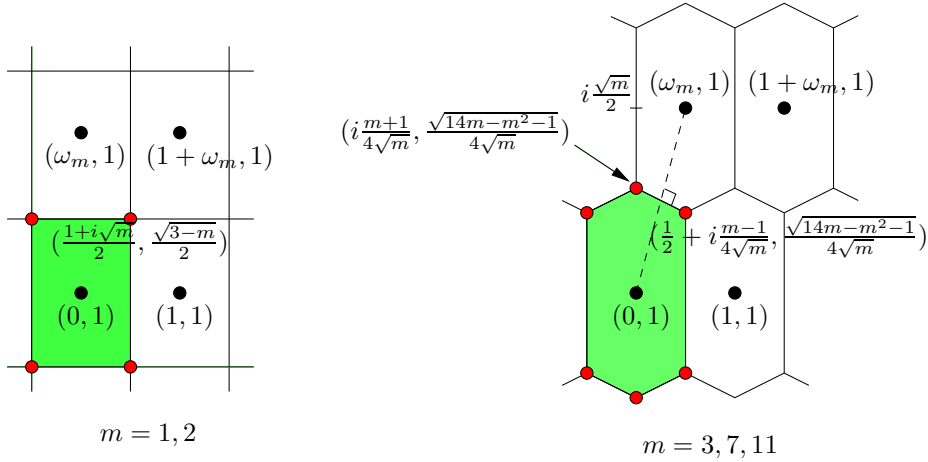
\cite[page 346]{Hat}.

\medskip
\begin{center}
\input{fig_forddomfond.pstex_t}
\end{center}

In the above picture, $\omega_m$ is equal to $i\sqrt{m}$ if
$m\equiv 1,2 \mod 4$ and is equal to $\frac{1+i\sqrt{m}}{2}$ if
$m\equiv 3 \mod 4$. The shaded area represents the vertical projection
to $\CC$ of the Ford fundamental domain at infinity $F_m$. The couples
are the coordinates in $\HH^3_\RR\subset\CC\times\RR$ of the finite
vertices of the polyhedron $F_m$ projecting to these points, and of
the center of the (unique) compact codimension $1$ face of $F_m$.

\smallskip
\noindent {\bf Cases $m=1,2$ : } A Ford fundamental domain of
$\Ga={\rm PSL}_2(\O_{-m})={\rm PSL}_2\big(\ZZ[i\sqrt{m}]\big)$ is
given by the polyhedron $F_m$ with  five vertices, one  at
$\infty$ and four finite ones at  $(\pm\frac
12\pm\frac {\sqrt{m}}{2} i,\frac {\sqrt {3-m}}{2})$. The diameter
$\Delta_m$ of the image of $K_m=\overline{F_m-C_\infty}$ in
$\Ga\bs\HH^3_\RR$ satisfies, by the symmetries and Equation
\eqref{eq:lengthcomput},
$$
\Delta_m\le d\Big(\big(\frac 12+\frac{\sqrt{m}}{2} i,
\frac {\sqrt {3-m}}{2}\big),(0,1)\Big)
=\log\frac{2+\sqrt{1+m}}{\sqrt{3-m}}\;,
$$
which is $\log(1+\sqrt 2)$ if $m=1$ and $\log(2+\sqrt
3)$ if $m=2$. Now, as in (1) above, $\spd{\D}$ is
contained in $[0,(1+\sqrt 2)^2]\subset[0,5.83]$ if $m=1$ and
$[0,(1+\sqrt 2)(2+\sqrt 3)]\subset[0,9.01]$ if $m=2$.

\medskip
\noindent {\bf Cases $m=3,7,11$ : }
A Ford fundamental domain of $\Ga={\rm PSL}_2(\O_{-m})={\rm
  PSL}_2\big(\ZZ[\frac{1+i\sqrt{m}}{2}]\big)$ is given by the
polyhedron $F_m$ with seven vertices at $\infty$, $\big(\pm
i\frac{m+1}{4\sqrt{m}},\frac{\sqrt{14\,m-m^2-1}}{4\sqrt{m}}\big)$ as
well as $(\pm\frac{1}{2}\pm
i\frac{m-1}{4\sqrt{m}},\frac{\sqrt{14\,m-m^2-1}}{4\sqrt{m}}\big)$. The
diameter $\Delta_m$ of the image of $K_m=\overline{F_m-C_\infty}$ in
$\Ga\bs\HH^3_\RR$ satisfies, by the symmetries and Equation
\eqref{eq:lengthcomput},
$$
\Delta_m\le d\Big(\big(i\frac{m+1}{4\sqrt{m}},
\frac{\sqrt{14\,m-m^2-1}}{4\sqrt{m}}\big),(0,1)\Big)
=\log\frac{4\sqrt{m}+ m + 1}{\sqrt{14m-m^2-1}}\;.
$$
Now, as in (1) above, $\spd{\D}$ is contained in $\big[0,(1+\sqrt
2)\frac{4\sqrt{m}+ m + 1}{\sqrt{14m-m^2-1}}\big]$ which is for
instance contained in $[0,4.664]$ if $m=3$.

\medskip (4) By the classification of the hyperbolic Coxeter simplices
(see for instance \cite[page 207]{VS}), there exists one, called $F$
thereafter, whose Coxeter diagram is
$$
\circ \!\!-\!\!\!-\!\!\!-\!\!    \circ \!\!-\!\!\!-\!\!\!-\!\!
\circ \!\!-\!\!\!-\!\!\!-\!\!    \circ \!\!-\!\!\!\overset{4}-
\!\!\!-\!\!\circ \!\!-\!\!\!-\!\!\!-\!\!\circ.
$$
Up to isometry of $\HH^5_\RR=\RR^4\times\;]0,+\infty[$, we may assume
that its ideal vertex is at infinity, and that the opposite face lies
on the Euclidean unit sphere centered at $0$.

Let $\Ga=\Ga_5$ be the group of isometries of $\HH^5_\RR$ generated by
the reflexions on the codimension-one faces of $F$.  The one-cusped
orbifold $\Ga_5\bs\HH_\RR^5$ is the minimal volume cusped hyperbolic
orbifold of dimension $5$, see \cite{Hil}.

The horoball $C_\infty$ is the maximal precisely invariant horoball
centered at $\infty$, see \cite[Prop.~5]{Hil}. It is easy to see (see for
instance \cite[page 216]{Hil}) that the vertical projection of $F$ in
$\RR^4$ (which is a Euclidean Coxeter simplex with Coxeter diagram $
\circ \!\!-\!\!\!-\!\!\!-\!\!  \circ \!\!-\!\!\!-\!\!\!-\!\!  \circ
\!\!-\!\!\!\overset{4}- \!\!\!-\!\!\circ \!\!-\!\!\!-\!\!\!-\!\!\circ
$ of type $\wt{F_4}$), has diameter $1/\sqrt 2$.  Thus the diameter of
$\overline{F-C_\infty}$ is at most $2d_{\HH^2_\RR}(i,\frac 1{\sqrt 2}+
i\frac 1{\sqrt 2})=2\log(1+\sqrt 2)$, and $\spd{\D}$ is contained in
$[0,(1+\sqrt 2)^3]\subset[0,14.08]$.

Let $\HH$ be the skew field of Hamilton's quaternions. The Hurwitz
ring $\O'$ consists of all quaternions in $\HH$ of the form $\frac
12(a_0+a_1i+a_2j+a_3k)$ such that the coefficients
$a_0,a_1,a_2,a_3\in\ZZ$ have equal parity. The {\it Hurwitz modular
  group} $\Ga={\rm PSL}_2(\O')$ (defined using the Dieudonn\'e
determinant) is, up to conjugation, the derived subgroup of $\Ga_5$,
which has index $4$ in $\Ga_5$, see \cite[page 186]{JW}. Since a
fundamental domain for $\Ga$ can be built as the connected union of
four copies of the fundamental domain $F$ of $\Ga_5$, the
approximation spectrum $\spd{\D}$ is contained in $[0,(1+\sqrt 2)^9]
\subset[0,2787]$, a very rough estimate.

\medskip (5) Before giving the fifth and last example, notice that the
model of $\hdc$ used therein will differ from the one used in Example
2 of Section \ref{sec:framework} to facilitate references to \cite{FP}
on which the example is based: We will use the Siegel domain of the
complex hyperbolic plane, whose underlying space is, as a subset of
the complex projective plane $\PP^2(\CC)$ with nonhomogeneous
coordinates,
$$
\HH^2_\CC=\big\{[W_0:W:1]\in
\PP^2(\CC)\;:\; 2\,{\rm Re}\;W_0 +|W|^2<0\big\}\;. 
$$ 
Consider the Hermitian form $q=Z_0\overline{Z_2} +Z_2\overline{Z_0}
+Z_1\overline{Z_1}$ on $\CC^3$, whose signature is $(1,2)$.  The {\it
  Eisenstein-Picard modular group} $\Ga={\rm PU}_q(\O_{-3})$ is the
projective unitary group of the form $q$ with coefficients in
$\O_{-3}$, acting projectively on $\HH^2_\CC$. Let
$$
C_\infty=\big\{[W_0:W:1]\in
\PP^2(\CC)\;:\; 2\,{\rm Re}\;W_0 +|W|^2\leq -2\big\}\;, 
$$
which is a horoball centered at $\infty=[-1:0:0]$. Let
$\omega=\frac{-1+i\sqrt 3}2$ and
$$
P=\begin{pmatrix} 1 & 1 & \omega\\0 & \omega & -\omega\\0 & 0 &
  1\end{pmatrix}, \;\;\;Q=\begin{pmatrix} 1 & 1 & \omega\\0 & -1 &
  1\\0 & 0 & 1\end{pmatrix}, \;\;\;R=\begin{pmatrix} 0 & 0 & 1\\0 & -1
  & 0\\1 & 0 & 0\end{pmatrix}\;,
$$ 
which define elements of $\Ga$. A fundamental domain $D$ for $\Ga$ is
constructed in \cite[Theo.~4.15]{FP}, as a simplex with one infinite
vertex at $\infty$, which is the geodesic cone with cone point
$\infty$ over a tetrahedron $T_0$ with four finite vertices
$$
z_0=[\bar\omega:0:1],\ \  z_1=[-1:-\omega:1],\ \ z_2=[-1:1:1],\ \  
z_3=[\omega:0:1].
$$ 
The images of the interior of $D$ by $P,Q,R$ are disjoint from $D$
since $D$ is a fundamental domain for $\Ga$. By \cite[Prop.~4.6]{FP},
the element $R$ maps $T_0$ to itself, and the geodesic cones with
vertex $\infty$ over the four faces of $T_0$ are paired, by $PQ^{-1}$
and $P$.

The horoball $C_\infty$, which is invariant by $P,Q$, and which meets
its image by $R$ only in $u_0=[-1:0:1]\in T_0$, is hence the maximal
precisely invariant horoball centered at $\infty$.  As $D$ is a cone
with vertex $\infty$ and $T_0$ intersects $C_\infty$ at $u_0$, the
diameter of the complement of the horoball $C_\infty$ in $D$ is
attained by two points of $T_0$.  Note that (see \cite{FP}, Definition
4.5 and the claim before Proposition 4.6 therein) the faces of $T_0$
are foliated by geodesic arcs between points of the edges, and that
the edges are geodesic arcs. Hence, by convexity of the distance map,
the maximal distance between two points of $T_0$ is attained by a pair
of vertices, that is by the maximal length of an edge of the
tetrahedron $T_0$.

The intersection of $\HH^2_\CC$ with the complex lines of equations
$W=0$ and $W_0=-1$ are totally geodesic, and are respectively a copy
of the constant curvature $-1$ real hyperbolic left halfplane and disc
of radius $\sqrt 2$ and center $0$ in $\CC$ (and $u_0$ corresponds to
the point $(-1,0)$ in this left halfplane and to the center of this
disc).  Hence, for $i=0$ and $i=3$, we have $d_{\HH^2_\CC}(z_i,u_0)=
\log(2+\sqrt{3})$ by Equation \eqref{eq:lengthcomput}. For $i=1$ and
$i=2$, we have $d_{\HH^2_\CC}(z_i,u_0)=\log(3+2\sqrt 2)$. Thus for
$0\leq i,j\leq 3$, by the triangle inequality, we have

$$
d_{\HH^2_\CC}(z_i,z_j)\leq
d_{\HH^2_\CC}(z_i,u_0)+d_{\HH^2_\CC}(z_j,u_0)\leq 2\log(3+2\sqrt{2})=
4\log(1+\sqrt{2})\;.
$$

Therefore, in the metric of $\HH^2_\CC$ with sectional curvature
between $-4$ and $-1$, we have the estimate
$$
\Delta\le \max_{0\leq i,j\leq 3}
d'_{\HH^2_\CC}(z_i,z_j)=2\log(1+\sqrt{2})\;,
$$
and a corresponding estimate on the approximation spectrum by
Proposition \ref{prop:specialDirichlet}
$$
\spd{\D}\subset[0,(1+\sqrt 2)^3]\subset[0,14.08]\;.
$$

%%%%%%%%%%%%%%%%%%%%%%%%
%% Modified 28.12.09/Frédéric
%%%%%%%%%%%%%%%%%%%%%%%%
\section{Hall rays in approximation spectra} 
\label{sec:hallray} 

In this section, for some quadruples of data $\D=(X,\Ga,\Ga_0,
C_\infty)$, we will prove that the approximation spectrum $\spd{\D}$
contains a segment $[0,c]$ for some $c>0$.

\medskip We start by recalling the following result from
\cite{PP2}. It says that given a family of almost disjoint
neighbourhoods of geodesic lines, there exists a geodesic ray or line,
with starting point (at infinity in the case of a line) any given
point outside these neighbourhoods, that has a prescribed penetration
in one given neighbourhood, and does not penetrate too much in the
neighbourhoods thereafter. We refer to Section \ref{sec:background}
for the definitions of the various penetration maps.

\btheo\label{theo:totgeoddimun} \cite[Theorem 5.9]{PP2}~ For every
$\epsilon>0$ and $\delta\geq0$, there exists a positive constant
$h'_1$ such that the following holds. Let $X$ be a complete simply
connected Riemannian manifold with sectional curvature at most $-1$
and dimension at least $3$. Let $(L_n)_{n\in\NN}$ be a family of
geodesic lines in $X$, such that $\diam(\N_\epsilon L_n\cap\N_\epsilon
L_m)\leq\delta$ for all $n\neq m$ in $\NN$. For every $\xi\in(X\cup
\partial_\infty X)-(\N_\epsilon L_0 \cup \partial_\infty L_0)$, let
$f_0 : T^1_\xi X \ra [0,+\infty[$ be either $f_0=\ftp_{L_0}$, or
$f_0=\ell_{\N_\epsilon L_0}$ if $X$ has constant curvature, or
$f_0=\cp_{L_0}$ (in which case $\xi\in\partial_\infty X -
\partial_\infty L_0$) if the metric spheres for the Hamenst\"adt
distances (on $\partial_\infty X-\{\xi'\}$ for any $\xi'\in 
\partial_\infty X$) are topological spheres. Let $h\geq h'_1$.

Then there exists a geodesic ray or line $\rho$ starting from $\xi$
and entering $\N_\epsilon L_0$ at time $0$ with $f_0(\rho) = h$, such
that $\ell_{\N_\epsilon L_n} (\rho)\le h'_1$ for every $n\neq 0$ such
that $\rho(]\delta,+\infty[)$ meets $\N_\epsilon L_n$. \cqfd \etheo

Note that the condition on the metric spheres of the Hamenst\"adt
distance being topological spheres is satisfied by all negatively
curved symmetric spaces.

The following result is an analog of Theorem 5.13 of \cite{PP2}, where
we considered cusp excursions. It has as hypothesis the conclusion of
the previous theorem. It says that if a given family of almost
disjoint neighbourhoods of geodesic lines is rich enough, then we can
find a geodesic line which has a prescribed upper asymptotic
penetration in these neighbourhoods.

We first define what we mean precisely by this.  Let $X$ be a proper
${\rm CAT}(-1)$ space and let $\xi\in X\cup \partial _\infty X$. Let
$\epsilon>0$, $\delta,\kappa\geq 0$.  Let $(C_\alpha)_{\alpha\in\A}$
be a family of convex subsets of $X$ such that $\diam(\N_\epsilon
C_\alpha\cap\N_\epsilon C_\beta)\leq\delta$ for all $\alpha\neq \beta$
in $\A$.  For each $\alpha\in\A$ such that $\xi\notin C_\alpha \cup
\partial_\infty C_\alpha$, let $f_\alpha:T^1_\xi X\ra [0,+\infty]$ be
a map such that $||f_\alpha-\ell_{\N_\epsilon C_\alpha}||_{\infty}\leq
\kappa$.  These assumptions guarantee that for every $\rho\in T^1_\xi
X$, the set $\E_\rho$ of times $t\geq 0$ such that $\rho$ enters in
some $C_\alpha$ at time $t$ with $f_\alpha(\rho)>\delta+\kappa$ is
discrete in $[0,+\infty[$, and that such an $\alpha$ is then unique,
denoted by $\alpha_t$.  Hence $\E_\rho=(t_i)_{i\in\N}$ for some
initial segment $\N$ in $\NN$, with $t_i<t_{i+1}$ for $i,i+1$ in $\N$.
With $a_i(\rho)=f_{\alpha_{t_i}}(\rho)$, the (finite or infinite)
sequence $\big(a_i(\rho)\big)_{i\in\N}$ will be called the {\em
penetration sequence} of $\rho$ with respect to $(\N_\epsilon
C_\alpha,f_\alpha)_ {\alpha\in\A}$ (and $\delta,\kappa$).  We will be
interested in the possible values of $\limsup_{i\to+\infty}
\;a_i(\rho)$, when $\N=\NN$.

\btheo \label{theo:limsuppenseq} Let $\epsilon>0$ and
$\delta,\nu,\nu'\geq 0$.  Let $X$ be a proper ${\rm CAT}(-1)$ space,
with $\partial_\infty X$ infinite, and let $\xi\in
X\cup \partial_\infty X$.  Let $(L_\alpha)_{\alpha\in\A}$ be a family
of geodesic lines in $X$, such that $\diam(\N_\epsilon
L_\alpha\cap\N_\epsilon L_\beta)\leq\delta$ for all $\alpha\neq \beta$
in $\A$.  For every $\alpha\in\A$ such that $\xi\notin \N_\epsilon
L_\alpha \cup \partial_\infty L_\alpha$, let $f_\alpha$ be either
$\ell_{\N_\epsilon L_\alpha}$ or $\ftp_{L_\alpha}$ or
$\cp_{L_\alpha}$, and in this last case, assume that $\xi
\in \partial_\infty X$. Let $\kappa$ be the upper bound of the
$||f_\alpha-\ell_{\N_\epsilon L_\alpha}||_\infty$ for all $\alpha$ in
$\A$ such that $\xi\notin \N_\epsilon L_\alpha \cup \partial_\infty
L_\alpha$.  Assume that $\bigcup_{\alpha\in\A} \partial_\infty
L_\alpha$ is dense in $\partial_\infty X$. Assume that for every
$h\geq \nu$ and $\alpha\in\A$ such that $\xi\notin \N_\epsilon
L_\alpha \cup \partial_\infty L_\alpha$, there exists a geodesic ray
or line $\rho$ starting from $\xi$ and entering $\N_\epsilon L_\alpha$
at time $t=0$ with $f_\alpha(\rho) = h$, and with $f_\beta(\rho)\leq
\nu'$ for every $\beta$ in $\A-\{\alpha\}$ such that
$\rho(]\delta,+\infty[)$ meets $\N_\epsilon L_\beta$. Let
$\big(a_i(\rho')\big)_{n\in\N}$ be the penetration sequence of a
geodesic ray or line $\rho'$ with respect to $(\N_\epsilon L_\alpha,
f_\alpha)_ {\alpha\in\A}$ (and $\delta,\kappa$).

Then, there exists $h_*=h_*(\epsilon,\delta,\kappa,\nu,\nu')>0$ such
that for every $h\geq h_*$, there exists a geodesic ray or line $\rho$
starting from $\xi$ such that
$$
\limsup_{i\to+\infty} \;a_i(\rho)=h\;.
$$
\etheo

\dem We start by recalling two lemmas from \cite{PP2}, which explain
the relative penetration behaviour of a pair of geodesic lines in the
$\epsilon$-neighbourhoods of geodesic lines.

\blemm \cite[Lemma 2.3]{PP2}~ \label{lem:c1}%
Let $C$ be a convex subset in $X$, let $\epsilon>0$ and let $\xi_0\in
(X\cup\partial_\infty X) -(\N_\epsilon C\cup\partial_\infty C)$. If
two geodesic rays or lines $\rho,\rho'$ which start from $\xi_0$
intersect $\N_\epsilon C$, then the first intersection points $x,x'$
of $\rho,\rho'$ respectively with $\N_\epsilon C$ are at a distance at
most $c'_1(\epsilon)= 2\,\arsinh(\coth\epsilon) $.  \cqfd 
\elemm

\blemm\cite[Lemmas 2.5 and 2.6]{PP2}~ \label{lem:h0} %
For every $\epsilon,\eta>0$, there exist (explicit) constants
$c'_2(\epsilon),c'_3(\epsilon)>0$ and $c(\epsilon,\eta)>0$ such that
the following holds. Let $X$ be a ${\rm CAT}(-1)$ space, $C$ a convex
subset in $X$, $\xi_0\in X\cup\partial_\infty X$, and $\rho,\rho'$ two
geodesic rays or lines starting from $\xi_0$. If $\rho$ enters
$\N_\epsilon C$ at a point $x\in X$ and exits $\N_\epsilon C$ at a
point $y\in X$ such that $d(x,y)\geq c(\epsilon,\eta)$ and
$d(y,\rho')\leq \eta$, then $\rho'$ enters $\N_\epsilon C$ at a point
$x'\in X$ such that $d(x,x')\leq c'_2(\epsilon)\,d(x,\rho')$ and exits
$\N_\epsilon C$ at a point $y'\in X\cup\partial_\infty X$ such that
$$
d(y,y')\leq c'_3(\epsilon)\,d(y,\rho')\;\;{\rm or} \;\; 
d(x',y')> d(x,y)\;. \;\;\mbox{\cqfd}
$$ 
\elemm

Let $X, (L_\alpha, f_\alpha)_{\alpha\in \A}, \xi, c,c',\kappa$ be as
in the statement of Theorem \ref{theo:limsuppenseq}. Note that by
the equations \eqref{eq:normftpl} and \eqref{eq:normcrpftp}, we have
$$
\kappa \leq 2\,c'_1(\epsilon)+2\epsilon+4\log(1+\sqrt{2})\;.
$$ 
In particular, $\kappa$ is finite. We start the proof of this theorem by
defining the constants that will be used therein. Let

$$
c_*= \kappa+ \max\big\{2(c'_1(\epsilon)+
\delta),\; c(\epsilon,\delta+c'_1(\epsilon)),\;
c'_1(\epsilon)c'_2(\epsilon)+c'_3(\epsilon)(c'_1(\epsilon)+\delta)+ 
\nu'+\kappa\big\}\;,
$$
where the positive constants $c'_i(\cdot)$ for $i=1,2,3$ and
$c(\cdot,\cdot)$ are defined in the lemmas \ref{lem:c1} and
\ref{lem:h0}.  Note that $c_*>\max\{\kappa+2\delta,\nu'\}$, since
$c'_1(\epsilon)\geq 1$ for all $\epsilon>0$. Let
\begin{equation}\label{hstar}
h_*=h_*(\epsilon,\delta,\kappa,\nu,\nu')=\max\{c_*,\nu\}\;.
\end{equation}

\medskip
Let $h\geq h_*$, and let $\alpha_0\in\A$ be such that $\xi\notin
\N_\epsilon L_{\alpha_0}\cup \partial_\infty L_{\alpha_0}$. The
existence of such an index follows from the assumptions: Indeed, as
$\partial_\infty X$ is 
(Hausdorff and) infinite, and by the density of $\bigcup_{\alpha\in\A}
\partial_\infty L_\alpha$, the set $\A$ is infinite; note that
$\partial_\infty L_\alpha\cap \partial_\infty L_\beta$ is empty if
$\alpha\neq \beta$, otherwise, as geodesic rays converging to the same
point at infinity become exponentially close, we would have
$\diam(\N_\epsilon L_\alpha\cap\N_\epsilon
L_\beta)=+\infty$; hence $\xi$ belongs to at most one $\partial_\infty
L_{\alpha}$ if $\xi\in\partial_\infty X$; if $\xi\in X$, then $\xi$
belongs to at most finitely many $\N_\epsilon L_{\alpha}$ for
$\alpha\in\A$, as $X$ is proper and $\diam(\N_\epsilon
L_\alpha\cap\N_\epsilon L_\beta)\leq \delta$ if $\alpha\neq \beta$.

As $h\geq h_*\geq \nu$, there exists, by the assumptions of Theorem
\ref{theo:limsuppenseq}, a geodesic ray or line $\rho_0$ starting from
$\xi$, entering $\N_\epsilon L_{\alpha_0}$ at time $t=0$, such that
$f_{\alpha_0}(\rho_0)=h$, and $f_\alpha(\rho_0)\le \nu'$ for every
$\alpha\neq \alpha_0$ such that $\rho_0(]\delta,+\infty[)$ meets
$\N_\epsilon L_\alpha$.

If a geodesic ray or line $\rho$ starting from $\xi$ meets
$\N_\epsilon L_\alpha$ such that $\xi\notin \N_\epsilon L_\alpha
\cup \partial_\infty L_\alpha$, let $t^-_\alpha(\rho)$ and
$t^+_\alpha(\rho)$ be the entrance and exit times.

We construct, by induction, sequences $(\rho_k)_{k\in\NN}$ of geodesic
rays or lines starting from $\xi$, $(\alpha_k)_{k\in\NN}$ of elements
of $\A$, and $(t_k)_{k\in\NN-\{0\}}$ of elements in $[0,+\infty[$
converging to $+\infty$, such that for every $k\in\NN$,
\begin{enumerate}
\item %
  $\rho_k$ enters the interior of $\N_\epsilon L_{\alpha_0}$ at time
  $0$, with $d(\rho_k(0),\rho_{k-1}(0))\leq \frac1{2^k}$ if $k\geq 1$;
\item %
  $\rho_k$ enters $\N_\epsilon L_{\alpha_k}$, $\xi\notin \N_\epsilon L_{\alpha_k}
\cup \partial_\infty L_{\alpha_k}$ and
  $f_{\alpha_k}(\rho_k)=h$;
\item %
  if $0\leq j\leq k-1$, then $\rho_k(]0,+\infty[)$ enters the interior
    of $\N_\epsilon L_{\alpha_j}$ before entering $\N_\epsilon
    L_{\alpha_k}$ with $t^-_{\alpha_j}(\rho_k)< t_k=
    t^+_{\alpha_{k-1}}(\rho_k)<t^+_{\alpha_{k}}(\rho_k)$;
\item %
  if $k\geq 1$, then for every $\alpha$ such that
  $\rho_k(]0,+\infty[)$ meets $\N_\epsilon L_\alpha$, we have 
\begin{itemize}
\item[$\bullet$] $\big| f_{\alpha}(\rho_k)- f_{\alpha} (\rho_{k-1})
  \big| < \frac1{2^k}$ if $t^-_\alpha(\rho_k)< t_k$, 
\item[$\bullet$] $f_{\alpha} (\rho_k)\leq c_*$ if $\alpha\neq\alpha_k$
  and $t_k\leq t^-_\alpha (\rho_k)\leq t^-_{\alpha_k}(\rho_k)+\delta$, 
\item[$\bullet$] $f_{\alpha} (\rho_k) \leq \nu'$ if
  $t^-_\alpha(\rho_k)> t^-_{\alpha_k}(\rho_k)+\delta$.
\end{itemize}
\end{enumerate}

Let us first prove that the existence of such sequences implies
Theorem \ref{theo:limsuppenseq}.  By the assertion (1), the sequence
$\big(\rho_k(0)\big)_{k\in\NN}$ stays at bounded distance from
$\rho_0(0)$, by a geometric series argument. Hence as $X$ is proper,
up to extracting a subsequence, the sequence $(\rho_k)_{k\in\NN}$
converges to a geodesic ray or line $\rho_\infty$ starting from $\xi$,
entering in $\N_\epsilon L_{\alpha_0}$ at time $t=0$, by the
continuity of the entering point in the interior of the
$\epsilon$-neighbourhood of a convex subset of $X$ (see for instance
\cite[Lemma 3.1]{PP2}). Let us prove that
$$
\limsup_{i\to+\infty}\;a_i(\rho_\infty)=h\;.
$$

The lower bound $\limsup_{i\to+\infty} \;a_i(\rho_\infty)\geq h$ is
immediate by a semicontinuity argument. Indeed, for every $k>i$ in
$\NN$, we have by the assertions (2), (3) and (4),
$$
\big| f_{\alpha_i}(\rho_k)-h\big|=
\big |f_{\alpha_i}(\rho_k)-f_{\alpha_i}(\rho_i)\big|
\leq \sum_{j=i}^{k-1}
\big|f_{\alpha_i}(\rho_{j+1})-f_{\alpha_i}(\rho_j)\big| 
\leq \sum_{j=i}^{k-1}\frac1{2^{j+1}}\leq \frac1{2^i}\;.
$$
Hence by the continuity of $f_{\alpha_i}$ (see Section
\ref{sec:background}), we have the inequality $f_{\alpha_i}
(\rho_\infty) \geq h-\frac1{2^i}$, whose right side converges to $h$
as $i$ tends to $+\infty$, which proves the lower bound, by the
definition of $\kappa$ and of the penetration sequence, as $h\geq
h_*\geq c_*>\delta+\kappa$.

\medskip To prove the upper bound $\limsup_{i\to+\infty}
\;a_i(\rho_\infty) \leq h$, assume by contradiction that there exists
$\eta >0$ such that for every $\lambda>0$, there exists $\alpha=
\alpha(\lambda)\in\A$ such that $\rho_\infty$ enters $\N_\epsilon
L_\alpha$ with $f_{\alpha}(\rho_\infty) \geq h+ \eta$ and
$t^-_\alpha(\rho_\infty) > \lambda+2\,c'_1(\epsilon)$, where
$c'_1(\epsilon)$ has been defined in Lemma \ref{lem:c1}. Take
$$
\lambda_0=\max\Big\{t_{i+1}:\frac 1{2^i}\ge\frac\eta 2\Big\}
$$
and $\alpha= \alpha(\lambda_0)$.

By continuity of $f_{\alpha}$, if $k$ is big enough, we have
$f_{\alpha}(\rho_k)\geq h+ \frac{\eta}{2}>h_*$. In particular,
$\alpha\neq \alpha_k$ by the assertion (2). Since
$$
h_*\geq c_*\geq \kappa\geq 
|f_{\alpha}(\rho_k)-\ell_{\N_\epsilon L_\alpha}(\rho_k)|\;,
$$ 
the geodesic $\rho_k$ meets $\N_\epsilon L_\alpha$. The entry time of
$\rho_k$ in $\N_\epsilon L_\alpha$ is positive, as $d(\rho_k(0),
\rho_\infty(0))\leq c'_1(\epsilon)$ and the entrance points of
$\rho_k$ and $\rho_\infty$ in $\N_\epsilon L_\alpha$ are at distance
at most $c'_1(\epsilon)$, both by Lemma \ref{lem:c1}, and as the
entrance time of $\rho_\infty$ in $\N_\epsilon L_\alpha$ is bigger
than $2\,c'_1(\epsilon)$. Hence, since $\nu'\leq c_*\leq h_*$ by the
definitions of $c_*$ and of $h_*$, we have $t^-_\alpha(\rho_k)< t_k$,
otherwise, by the assertion (4), $f_{\alpha}(\rho_k)\leq \max\{c_*,
\nu'\} = c_*\leq h_*$, a contradiction. Let $i\leq k-1$ be the minimum
element of $\NN$ such that for $j=i,\dots, k-1$, the geodesic
$\rho_{j+1}$ meets $\N_\epsilon L_\alpha$ at a positive time with
$t^-_\alpha(\rho_{j+1}) < t_{j+1}$.  By the triangle inequality, we
have
$$
\big |t^-_\alpha(\rho_{i+1})-t^-_\alpha(\rho_\infty) \big|\leq
d\big(\rho_{i+1}(t^-_\alpha(\rho_{i+1})),
\rho_\infty(t^-_\alpha(\rho_\infty))\big)+
d\big(\rho_{i+1}(0),\rho_\infty(0)\big)\leq 2\,c'_1(\epsilon)\;,
$$
by applying twice Lemma \ref{lem:c1}. Hence, by the definition of $i$
and of $\alpha$,
$$
t_{i+1}> t^-_\alpha(\rho_{i+1})\geq t^-_\alpha(\rho_\infty)
-2\,c'_1(\epsilon)>
\lambda_0+2\,c'_1(\epsilon)-2\,c'_1(\epsilon)=\lambda_0.
$$  
By the definition of $\lambda_0$, we hence have $\frac{1}{2^{i}}<
\frac{\eta}{2}$.  By the definition of $i$ and by the assertion
(4), we have
\begin{align*}
f_{\alpha}(\rho_i)=\; & f_{\alpha}(\rho_k)+\sum_{j=i}^{k-1}
\big(f_{\alpha}(\rho_{j})-f_{\alpha}(\rho_{j+1})\big)
\geq h+ \frac{\eta}{2} -
\sum_{j=i}^{k-1} \frac{1}{2^{j+1}}\\ \geq \;&
h+ \frac{\eta}{2} - \frac{1}{2^{i}}\geq
h\geq h_*\;,
\end{align*}
and in particular by the same argument as for $\rho_k$ above, $\rho_i$
enters $\N_\epsilon L_\alpha$ at a positive time and
$t^-_\alpha(\rho_i)< t_i$.  This contradicts the minimality of $i$.
This completes the proof of Theorem \ref{theo:limsuppenseq}, assuming
the existence of the sequences with the properties (1)--(4).

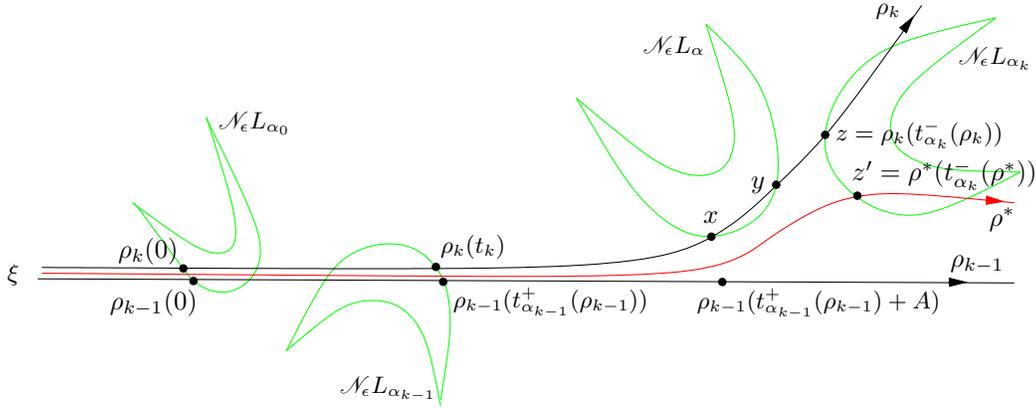
\begin{center} 
\input{fig_limsup_geod.pstex_t}
\end{center}

Let us now construct the sequences $(\rho_k)_{k\in\NN}$,
$(\alpha_k)_{k\in\NN}$ and $(t_k)_{k\in\NN-\{0\}}$. We already have
defined $\rho_0$ and $\alpha_0$, and they satisfy the properties
(1)--(4). Let $k\geq 1$, and assume that $\rho_{k-1}$, $\alpha_{k-1}$,
as well as $t_{k-1}$ if $k\geq 2$, have been constructed.

For every small $\mu>0$ and big $A>0$ (to be precised later on),
consider the set $\E'=\E'_k$ of couples $(\alpha,\rho)$ where $\alpha
\in\A$ and $\rho$ is a geodesic ray or line, starting from $\xi$,
entering $\N_\epsilon L_{\alpha_0}$ at time $t=0$, which is
$\mu$-close to $\rho_{k-1}$ on $[0,t^+_{\alpha_{k-1}}(\rho_{k-1})+A]$,
and which enters $\N_\epsilon L_\alpha$ with $t^-_{\alpha}(\rho)\geq
t^+_{\alpha_{k-1}}(\rho_{k-1})+\frac{A}{2}$. This set $\E'$ is not
empty, as $\bigcup_{\alpha\in\A} \partial_\infty L_\alpha$ is dense in
$\partial_\infty X$, and as the assumption that $\diam(\N_\epsilon
L_\alpha\cap \N_\epsilon L_\beta)\leq\delta$ for all $\alpha\neq
\beta$ in $\A$ implies that a compact subset of $X$ meets only
finitely many $\N_\epsilon L_{\alpha}$ for $\alpha\in\A$ (we may even
find such a couple $(\alpha,\rho)$ with $\rho(+\infty)
\in \partial_\infty L_\alpha$, which implies that a subray of $\rho$
is contained in $\N_\epsilon L_{\alpha}$).  Let $(\alpha_k,\rho^*)$ be
an element of $\E'$ with $t^-_{\alpha_k}(\rho^*)$ minimal, which
exists since the family $(\N_\epsilon L_\alpha)_{\alpha\in\A}$ is
locally finite and by a continuity argument (when $\A$ is given the
discrete topology, the subset $\E'$ is closed).  Note that $\xi\notin
\N_\epsilon L_{\alpha_k}\cup\partial_\infty L_{\alpha_k}$, as
$t^-_{\alpha_k}(\rho^*)>0$.

By the last hypothesis of Theorem \ref{theo:limsuppenseq}, let
$\rho_k$ be a geodesic ray or line starting from $\xi$ with
$f_{\alpha_k}(\rho_k) =h$ and $f_{\alpha}(\rho_k)\le \nu'$ for every
$\alpha\in\A$ such that $\rho_k (]t^-_{\alpha_k}(\rho_k) +\delta,
+\infty[)$ enters $H_\alpha$. In particular, this proves the assertion
(2) at rank $k$, since
\begin{equation}\label{eq:probdelta}
h\geq h_*\geq c_*> \kappa\geq |f_{\alpha_k}(\rho_k)-
\ell_{\N_\epsilon L_{\alpha_k}} (\rho_k)|\;.
\end{equation}
Let $z=\rho_k (t_{\alpha_k}^-(\rho_k))$ and $z'=\rho^*
(t_{\alpha_k}^-(\rho^*))$ be the entering points of $\rho_k$ and
$\rho^*$ in $\N_\epsilon L_{\alpha_k}$. By Lemma \ref{lem:c1}, we have
$d(z,z')\leq c'_1(\epsilon)$. Hence, by hyperbolicity and as
$(\alpha_k,\rho^*)$ is in $\E'$, if $A$ is big
enough, then $\rho_k$ is $(2\mu)$-close to $\rho_{k-1}$ between $\xi$
and $\rho_{k-1}(t^+_{\alpha_{k-1}} (\rho_{k-1})+1)$. In particular, if
$\mu$ is small enough, and using properties (1) and (3) at rank
$k-1$, we have the following properties.
\begin{enumerate}
\item[$\bullet$] The geodesic ray or line $\rho_k$ enters the interior
  of $\N_\epsilon L_{\alpha_0}$, at a time that we may assume to be
  $0$, with $d(\rho_k(0),\rho_{k-1}(0))\leq \frac1{2^k}$ (this proves
  the assertion (1) at rank $k$).
\item[$\bullet$] For $0\leq j\leq k-1$, the geodesic ray or line
  $\rho_{k-1}$ meets the interior of $\N_\epsilon L_{\alpha_{j}}$ at a
  time strictly between $0$ and $t^+_{\alpha_{k-1}}(\rho_{k-1})$, by
  the inductive assertions (3) if $k\neq 1$ and $j\leq k-2$, or (1) if
  $k=1$ or (2) if $j=k-1$ (by Equation \eqref{eq:probdelta} where $k$
  has been replaced by $k-1$). Hence the geodesic ray $\rho_k$ also meets
  the interior of $\N_\epsilon L_{\alpha_j}$ at a time strictly
  between $0$ and $t^+_{\alpha_{k-1}}(\rho_{k-1})$. This allows, in
  particular, to define $t_k= t^+_{\alpha_{k-1}}(\rho_k)$, and proves
  the assertion (3) at rank $k$.
\item[$\bullet$] For every $\alpha$ such that $\rho_k(]0,+\infty[)$
  meets $\N_\epsilon L_\alpha$ and $t^-_\alpha(\rho_k)< t_k$, we may
  assume, by the continuity of $f_\alpha$, up to taking $\mu$ small
  enough, that $\big| f_{\alpha}(\rho_k)-f_{\alpha}(\rho_{k-1})\big| <
  \frac1{2^k}$.
\end{enumerate}
Hence (using also the construction of $\rho_k$), to prove the
assertion (4) at rank $k$, we consider $\alpha\in\A-\{\alpha_k\}$ such
that $\rho_k$ meets $\N_\epsilon L_\alpha$ with $t_k \leq
t^-_\alpha(\rho_k) \leq t^-_{\alpha_k}(\rho_k)+\delta$, and we prove
that $f_{\alpha} (\rho_k)\leq c_*$.

\medskip Assume by absurd that $f_{\alpha} (\rho_k)> c_*$. In
particular, $\ell_{\N_\epsilon L_\alpha}(\rho_k)\geq c_*-\kappa>0$ (by
the definition of $c_*$), so that $\rho_k$ enters $\N_\epsilon
L_\alpha$. Let $x=\rho_k(t^-_{\alpha}(\rho_k))$ be the entering point
of $\rho_k$ in $\N_\epsilon L_\alpha$. Note that
$$
\ell_{\N_\epsilon L_{\alpha_k}}(\rho_k)\geq f_{\alpha_k}(\rho_k)-\kappa=
h-\kappa\geq h_*-\kappa\geq c_*-\kappa\;,
$$ 
by the definition of $h_*$. If $t^-_\alpha(\rho_k) \geq
t^-_{\alpha_k}(\rho_k)$, then, since $t^-_\alpha(\rho_k)\leq
t^-_{\alpha_k}(\rho_k)+\delta$ and
$$
\min\{\ell_{\N_\epsilon L_\alpha}(\rho_k),
\ell_{\N_\epsilon L_{\alpha_k}}(\rho_k)\} \geq c_*-\kappa>2\delta
$$ 
by the definition of $c_*$, this would imply that the intersection
$\N_\epsilon L_\alpha\cap\N_\epsilon L_{\alpha_k}$ has diameter bigger
than $\delta$, which contradicts $\alpha\neq \alpha_k$. Hence
$t^-_\alpha(\rho_k) < t^-_{\alpha_k}(\rho_k)$ (which implies that
$x\in[\xi,z]$) and we have $t^+_\alpha(\rho_k)\leq t^-_{\alpha_k}
(\rho_k) +\delta$, again since $\diam(\N_\epsilon L_\alpha\cap
\N_\epsilon L_{\alpha_k})\leq\delta$ and $\alpha\neq\alpha_k$.  In
particular, $y=\rho_k(t^+_{\alpha}(\rho_k))$ is a point in $X$.

We want to apply Lemma \ref{lem:h0} with $\eta=\delta+c'_1(\epsilon)$,
$\rho=\rho_k$, $\rho'=\rho^*$, $C=L_\alpha$ and $\xi_0=\xi$. We first
check the hypotheses of this lemma.

We do have that $\rho_k$ enters $\N_\epsilon L_\alpha$ at $x$ and
exits it at $y$,  and
\begin{equation}\label{eq:minodxy}
d(x,y)\geq f_{\alpha}(\rho_k)-\kappa>
c_*-\kappa\geq c(\epsilon,\eta)\;,
\end{equation}
by the definition of $c_*$.  As
$$
d(\rho_k(t^-_{\alpha_k}(\rho_k)+\delta),\rho^*)\leq \delta+
d(\rho_k(t^-_{\alpha_k}(\rho_k)),\rho^*)\leq \delta+d(z,z')\leq \eta\;,
$$
and by convexity, we have $d(y,\rho^*)\leq \eta$. 

Hence we may indeed apply Lemma \ref{lem:h0}, and the geodesic
$\rho^*$ enters $\N_\epsilon L_\alpha$ at a point $x'$ such that
\begin{equation}\label{eq:dxx}
d(x,x')\leq c'_2(\epsilon)\,d(x,\rho^*)\leq 
c'_2(\epsilon)\,d(z,\rho^*)\leq c'_2(\epsilon)\,c'_1(\epsilon) \;,
\end{equation}
where the middle inequality holds by convexity, and the last one since
$d(z,z')\leq c'_1(\epsilon)$.

Furthermore, the geodesic $\rho^*$ exits $\N_\epsilon L_\alpha$ at a
point $y'$ (possibly at infinity) and, by the alternative at the end
of Lemma \ref{lem:h0} and the equations \eqref{eq:minodxy} and
\eqref{eq:dxx},
$$
d(x',y')> d(x,y) > c_*-\kappa\geq \nu'+\kappa
$$ 
or
\begin{align*}
d(x',y')&\geq d(x,y)-d(x,x')-d(y,y')\geq 
d(x,y)-d(x,x')-c'_3(\epsilon)\,d(y,\rho^*)\\ &
> (c_*-\kappa) - c'_2(\epsilon)\,c'_1(\epsilon)- c'_3(\epsilon)\eta
\geq \nu'+\kappa\;,
\end{align*}
by the definition of $c_*$. In both cases, $d(x',y')> \nu'+\kappa$.

\medskip Let us prove that 
\begin{equation}
\label{eq:ineqcontradicpart}
t^-_{\alpha_k}(\rho^*)> t^-_{\alpha}(\rho^*)\;.
\end{equation} 
Otherwise, the point $z'$ belongs to $[\xi,x']$. Hence, with $x'',z''$
the closest points to $x',z'$ respectively on $\rho_k$, we have
$z''\in[\xi,x'']$. Note that $x''\in[\xi,y]$ since $d(x,x'')\leq
d(x,x')\leq c'_1(\epsilon)$ and $d(x,y)>c_*-\kappa>c'_1(\epsilon)$ by
Equation \eqref{eq:minodxy} and the definition of $c_*$.

Respectively by Equation \eqref{eq:minodxy}, by the triangle
inequality, since $z''\in[\xi,x'']$ and $x''\in[\xi,y]$, since
$t^-_\alpha(\rho_k)\leq t^-_{\alpha_k}(\rho_k)+\delta$, since closest
point maps do not increase distances, and by Lemma \ref{lem:c1},
\begin{align*}
c_*-\kappa&\leq d(x,y)\leq d(x,x'')+d(x'',y)\leq d(x,x'')+d(z'',y)\leq 
d(x,x'')+d(z'',z)+\delta\\ &
\leq 
d(x,x')+d(z',z)+\delta\leq 2\,c'_1(\epsilon)+\delta\;,
\end{align*}
which contradicts the definition of $c_*$.

\medskip 
Now, recall the constants $\mu>0$ and $A\geq 0$ introduced in
the definition of $\rho^*$. It follows from Equation
\eqref{eq:ineqcontradicpart}, and from the minimality assumption in
the definition of $\alpha_k$, that we have
$$
t^-_{\alpha}(\rho^*)<t^+_{\alpha_{k-1}}(\rho_{k-1})+\frac{A}{2}\;.
$$
Assume that $\mu$ is small enough and that $A$ is big enough.  Since
$d(x',y')>\nu'+\kappa$, and as $\rho^*$ is $\mu$-close to $\rho_{k-1}$
on $[0, t^+_{\alpha_{k-1}}(\rho_{k-1})+A]$, this implies that
$\rho_{k-1}$ enters $\N_\epsilon L_\alpha$ at a point $x^\sharp$ close
to $x'$, and exits at a point $y^\sharp$ (possibly at infinity)
such that $d(x^\sharp,y^\sharp)> \nu'+\kappa$. Hence
$$
f_\alpha(\rho_{k-1})\geq d(x^\sharp,y^\sharp)-\kappa> \nu'\;.
$$
This implies that $t^-_\alpha(\rho_{k-1})\leq
t^-_{\alpha_{k-1}}(\rho_{k-1})+\delta$, otherwise we have in
particular that $\alpha\neq\alpha_{k-1}$ and by the assertion (4) at
rank $k-1$ if $k\geq 2$ or by the construction of $\rho_0$ if $k=1$,
we would have $f_\alpha(\rho_{k-1})\leq \nu'$. Hence
\begin{align*}
t^+_{\alpha_{k-1}}(\rho_{k-1})-t^-_\alpha(\rho_{k-1})& \geq
t^+_{\alpha_{k-1}}(\rho_{k-1})-t^-_{\alpha_{k-1}}(\rho_{k-1})-\delta
=\ell_{\N_\epsilon L_{\alpha_{k-1}}}(\rho_{k-1})-\delta\\ &
\geq f_{\alpha_{k-1}}(\rho_{k-1})-\kappa-\delta=h-\delta-\kappa\geq 
h_*-\delta-\kappa\geq 
c_*-\delta-\kappa\;.
\end{align*}
That is, $\rho_{k-1}$ enters in $\N_\epsilon L_\alpha$ well before
exiting $\N_\epsilon L_{\alpha_{k-1}}$, the amount of time being at
least the constant $c_*-\delta-\kappa$ (which is positive by the
definition of $c_*$). But since the entering points in $\N_\epsilon
L_\alpha$, as well as the exiting points out of $\N_\epsilon
L_{\alpha_{k-1}}$, of the geodesic rays or lines $\rho_{k-1}$, $\rho*$ and
$\rho_k$ are very 
close, this contradicts the fact that $t^-_{\alpha}(\rho_k)\geq
t_k=t^+_{\alpha_{k-1}}(\rho_k)$.

This proves the result.  
\cqfd

\bcoro 
\label{coro:corohalgene}
Let $X$ be a complete simply connected Riemannian manifold with
sectional curvature at most $-1$ and dimension at least $3$, such that
the metric spheres for the Hamenst\"adt distances (on $\partial_\infty
X-\{\xi'\}$ for any $\xi'\in\partial_\infty X$) are topological spheres.
Let $\Ga$ be a discrete group of isometries of $X$ with finite
covolume, and let $\ga_0$ be a hyperbolic element of $\Ga$.  Let
$\xi_0\in \partial_\infty X$ be a parabolic fixed point, and $H_0$ be
a horosphere centered at $\xi_0$. For every $\xi\in \partial_\infty
X$ which is not a fixed point of a conjugate of $\ga_0$ or a parabolic
fixed point, define
$$
c'(\xi)=\liminf 
\frac{d_{\xi_0,H_0}(\xi,\ga_-)}{d_{\xi_0,H_0}(\ga_+,\ga_-)}\;,
$$
where the lower limit is taken over the conjugates $\ga$ of $\ga_0$ or
its inverse, with fixed points $\ga_-,\ga_+$ and
$d_{\xi_0,H_0}(\ga_+,\ga_-)$ tending to $0$.

Then the subset of $\RR$ consisting of the $c'(\xi)$ for $\xi\in
\partial_\infty X$ which is neither a fixed point of a conjugate of
$\ga_0$ nor a parabolic fixed point, contains a segment $[0,c]$ for
some $c>0$. 
\ecoro

\dem We will apply Theorem \ref{theo:limsuppenseq} with
$(L_\alpha)_{\alpha\in\A}$ the family of translation axes of the
conjugates of the element $\ga_0$ (where each line appears exactly
once), with $\xi=\xi_0$ and with $f_\alpha=\cp_{L_\alpha}$ for every
$\alpha$ in $\A$. Let
$\kappa=2\,c'_1(\epsilon)+2\epsilon+4\log(1+\sqrt{2})$, which
satisfies $\|f_\alpha-\ell_{N_\epsilon L_\alpha}\|_\infty\leq\kappa$
by Section \ref{sec:background}.

For some positive $\epsilon$ and $\delta$, this family satisfies the
assumption that $\diam\big(\N_\epsilon L_\alpha\cap\N_\epsilon
L_\beta\big)\leq\delta$ for all $\alpha\neq \beta$ in $\A$. Otherwise,
there would exist a sequence $(\ga_n)_{n\in\NN}$ in $\Ga-\Ga_0$,
where $\Ga_0$ is the stabilizer in $\Ga$ of the translation axis $L_0$
of $\ga_0$, such that $\diam\big(\N_\epsilon L_0\cap\ga_n\N_\epsilon
L_{0}\big)$ converges to $+\infty$. Up to multiplying $\ga_n$
on the right and on the left by a power of $\ga_0$ or its inverse, the
element $\ga_n$ moves a point of $L_0$ less than a constant. Hence
$\ga_n$ stays in a compact subset of the isometry group of $X$. By
discreteness, up to extracting a subsequence, $\ga_n$ does not
depend on $n$. But then $L_0$ and $\ga_1L_{0}$ are two distinct
translation axes that meet at least in one point at infinity, which
contradicts the discreteness of $\Ga$.

As $\Ga$ has finite covolume, the set of fixed points of the
conjugates of $\ga_0$ is dense in $\partial_\infty X$, hence
$\bigcup_{\alpha\in\A} \partial_\infty L_\alpha$ is dense in
$\partial_\infty X$.  The last hypothesis of Theorem
\ref{theo:limsuppenseq} holds true by Theorem \ref{theo:totgeoddimun},
with $\nu=h'_1$ and $\nu'=h'_1+\kappa$ (by definition of $\kappa$). By
Theorem \ref{theo:limsuppenseq}, there exists $h_*>0$ such that for
every $h\geq h_*$, there exists a geodesic line $\rho$ starting from
$\xi_0$ such that $\limsup_{i\to+\infty} \;a_i(\rho)=h$ where
$a_i(\rho)$ is the penetration sequence of $\rho$ with respect to
$(\N_\epsilon L_\alpha, f_\alpha)_ {\alpha\in\A}$ (and
$\delta,\kappa$).

For every $i$ in $\NN$, let $\alpha_i\in\A$ be the unique element such
that $a_i(\rho)=\cp_{L_{\alpha_i}}(\rho)$. For every $\alpha\in\A$
such that $\rho$ meets $\N_\epsilon L_\alpha$ at a time big enough
with $\cp_{L_{\alpha}}(\rho)>0$, let $L_{\alpha,\pm}$ be the two
endpoints of $L_{\alpha}$ such that, by the definition of the
crossratio penetration map and by Equation \eqref{eq:formcrosham},
$$
\cp_{L_{\alpha}}(\rho)=[\xi_0,L_{\alpha,-},\rho(+\infty),L_{\alpha,+}]
=\log\frac{d_{\xi_0,H_0}(L_{\alpha,+},L_{\alpha,-})}
{d_{\xi_0,H_0}(\rho(+\infty),L_{\alpha,-})}\;.
$$
Only finitely many $L_\alpha$'s meet a given compact subset of $X$.
Thus, for every subsequence $(i_k)_{k\in \NN}$ such that the sequence
$\big(\cp_{L_{\alpha_{i_k}}}(\rho)\big)_{k\in\NN}$ is bounded, since
the entrance time of $\rho$ in $\N_\epsilon L_{\alpha_{i_k}}$ tends to
$+\infty$, the distance $d_{\xi_0,H_0}(L_{\alpha_{i_k},+},
L_{\alpha_{i_k},-})$ tends to $0$ as $k\ra+\infty$. Also note that, by
definition of the penetration sequence, if $\alpha\in\A$ does not
belong to $\{\alpha_i\;:\;i\in \NN\}$, then either $\rho$ does not
meet $\N_\epsilon L_\alpha$ at a positive time, or
$\cp_{L_{\alpha}}(\rho)\leq \delta+\kappa< h_*\leq h$, see Equation
\eqref{hstar}.

Finally, if $(\ga_k)_{k\in\NN}$ is a sequence of conjugates of $\ga_0$
or its inverse with fixed points $\ga_{k,-},\ga_{k,+}$ and with
$d_{\xi_0,H_0}(\ga_{k,-},\ga_{k,+})$ tending to $0$, such that
 the sequence $\frac{d_{\xi_0,H_0}(\rho(+\infty),\ga_{k,-})}
{d_{\xi_0,H_0}(\ga_{k,-},\ga_{k,+})}$ is bounded from above, then
$\ga_{k,-}$ tends to $\rho(+\infty)$. Hence, if
$$
\liminf_{k\ra+\infty}\;\frac{d_{\xi_0,H_0}(\rho(+\infty),\ga_{k,-})}
{d_{\xi_0,H_0}(\ga_{k,-},\ga_{k,+})}\leq e^{-h},
$$ 
then for every $\epsilon\in\;]0,h_*-\delta-\kappa[\,$, for $k$ big
enough, there exists $\alpha\in\A$ such that $L_{\alpha}$ is the
translation axis of $\ga_k$, and $\cp_{L_{\alpha}}(\rho)\geq
h-\epsilon\geq h_*-\epsilon> \delta+\kappa$ so that $\rho$ meets
$\N_\epsilon L_{\alpha}$, at a positive time. In particular $\alpha$
belongs to $\{\alpha_i\;:\;i\in \NN\}$. Therefore, we have
$$
\liminf 
\frac{d_{\xi_0,H_0}(\rho(+\infty),\ga_-)}{d_{\xi_0,H_0}(\ga_+,\ga_-)}
=e^{-h}\;,
$$ 
where the lower limit is taken as in the statement of the corollary.
This proves the result, with $c=e^{-h_*(\epsilon,\delta,\kappa, h'_1,
  h'_1+\kappa)}$.  \cqfd

\medskip Specializing the above Corollary \ref{coro:corohalgene} to
the particular cases of the real or complex hyperbolic space (see the
examples at the end of Section \ref{sec:framework}), we have the
following applications.

\bcoro \label{coro:hallrayconstcurv} Let $n\geq 3$, let $\Ga$ be a
discrete group of isometries of $X=\HH^n_\RR$ with finite covolume,
and let $\Ga_0$ be the stabilizer in $\Ga$ of the translation axis of
a hyperbolic element of $\Ga$.  Let $C_\infty$ be a precisely
invariant horoball centered at a parabolic fixed point of $\Ga$, and
$\D=(X,\Ga,\Ga_0,C_\infty)$.
Then $\spd{\D}$ contains a segment $[0,c]$ for some $c>0$.  \cqfd
\ecoro

\medskip 
By the last equality in the proof of Corollary \ref{coro:corohalgene}
(and since the constant $h'_1$ appearing in Theorem
\ref{theo:totgeoddimun} is explicited in \cite{PP2}), if one wants in
particular situations to be able to give an explicit (lower bound on
the) constant $c$ appearing in Corollary \ref{coro:hallrayconstcurv}
(which is the same as in Corollary \ref{coro:corohalgene}), one only
needs to find explicit $\epsilon,\delta$ such that
$\diam\big(\N_\epsilon L_\alpha\cap\N_\epsilon L_\beta\big)\leq\delta$
for all $\alpha\neq \beta$ in $\A$. We give such a computation in the
following remark.

\brema \label{rem:computinters} In the real hyperbolic upper halfplane
$\HH^2_\RR$, consider the geodesic line $L$ with endpoints
$\frac{1+\sqrt{5}}{2}$ and $\frac{1-\sqrt{5}}{2}$. Let
$(L_\alpha)_{\alpha\in\A}$ be the family of the images of $L$ by ${\rm
PSL}_2(\ZZ)$, acting by homographies on $\HH^2_\RR$, modulo the (global)
stabilizer of $L$. Let $\epsilon=\frac{\log 5}{2}$ and $\delta=2\log
(2+\sqrt{5})$. Then $\diam\big(\N_\epsilon L_\alpha\cap\N_\epsilon
L_\beta\big)\leq\delta$ for all $\alpha\neq \beta$ in $\A$.  \erema

\begin{center}
\input{fig_delteps.pstex_t}
\end{center}

\dem Call {\it Weierstrass points} the points of $\HH^2_\RR$ in the
image of $i$ by ${\rm PSL}_2(\ZZ)$.  Note that, by an easy
computation, the $\epsilon$-neighbourhood of $L$ is the tubular
neighbourhood of $L$ of biggest radius such that the Weierstrass
points in its interior lie on $L$. As seen by considering the
fundamental domain for the integer horizontal translations in
$\HH^2_\RR$, between the geodesic lines with endpoints $0,\infty$ and
$1,\infty$ respectively, if two images of $L$ by elements of ${\rm
  PSL}_2(\ZZ)$ are disjoint, then the interiors of their
$\epsilon$-neighbourhoods do not meet. If two images of $L$ by
elements of ${\rm PSL}_2(\ZZ)$ are distinct but meet, to prove that
the diameter of the intersection of their $\epsilon$-neighbourhoods is
at most $\delta$, we may assume that these images are $L$ and the
image $L'$ of $L$ by the translation by $-1$.

Recall that the boundary of the $\epsilon$-neighbourhood of the
   geodesic line carried by the Euclidean circle of center $x\in\RR$ and
   Euclidean radius $r$ is the union of two arcs of circles between $x+r$
   and $x-r$ that are invariant by reflection in the vertical line
   through $x$. Since the upper arc of circle $\partial_+ \N_\epsilon L'$
   of $\partial \N_\epsilon L'$ is tangent to the vertical line through
   the Weierstrass point $i+1$, its Euclidean center is the point
   $-\frac{1}{2}+i$, and its radius is $\frac 32$.
   Let $L''$ be the geodesic line with endpoints
   $-1,1$, which is a bisectrix of $L$ and $L'$. Let $u$ and $u'$ be
   the intersection points of $L''$ with $\partial_- \N_\epsilon L$ and
   $\partial_- \N_\epsilon L'$ respectively. By considering the Euclidean
   quadrangle with vertices at $u, i-\frac{1}{2}, -\frac{1}{2},0$, we
   easily compute that the Euclidean height of $u$ is
   $\frac{1}{\sqrt{5}}$. Using Formula \eqref{eq:lengthcomput}, we hence
   have $d(i,u)= \log(2+\sqrt{5})$. The distance between $i$ and the
   intersection point of $\partial_+ \N_\epsilon L'$ and $\partial_+
   \N_\epsilon L$ is $\log(1+\sqrt{2})<\log(2+\sqrt{5})$. Hence by
   convexity and symmetry arguments, the intersection $\N_\epsilon L'\cap
   \N_\epsilon L$ is contained in the hyperbolic ball of center $i$ and radius
   $\log(2+\sqrt{5})$. Therefore, the diameter of this intersection is
   $d(u,u')=2
   \log(2+\sqrt{5})$. The result follows.
   \cqfd

\medskip
The following result is a consequence of Corollary
\ref{coro:hallrayconstcurv} and Equation \eqref{eq:divspec}, and
proves the second claim of Theorem \ref{theo:mainintrodeux} in the
Introduction.

\bcoro Let $M$ be a geometrically finite complete connected Riemannian
manifold  of constant sectional curvature $-1$ and of dimension at least
$3$. Let $A_0$ be a closed geodesic in $M$, and let $A_\infty$ be a
Margulis neighbourhood of a cusp of M. Then the spiraling spectrum
$\operatorname{Sp}_{A_\infty,A_0}(M)$ around $A_0$
contains a segment $[0,c]$ for some $c>0$.  \cqfd 
\ecoro

\medskip
Using Equation \eqref{eq:egalhamcyg}, we obtain the following result in
the complex hyperbolic case.

\bcoro \label{cor:hallcomphyp} Let $n\geq 2$, let $\Ga$ be a discrete
group of isometries of the Siegel domain model of $\HH^n_\CC$ with
finite covolume, and let $\ga_0$ be a hyperbolic element of $\Ga$.
Assume that the point $\infty$ is a parabolic fixed point of
$\Ga$. For every $\xi\in \partial_\infty \HH^n_\CC$ which is neither a
fixed point of a conjugate of $\ga_0$ nor a parabolic fixed point,
define
$$
c'(\xi)=\liminf 
\frac{d_{\rm Cyg}(\xi,\ga_-)}{d_{\rm Cyg}(\ga_+,\ga_-)}\;,
$$
where the lower limit is taken over the conjugates $\ga$ of $\ga_0$ or
its inverse, with fixed points $\ga_-,\ga_+$ and $d_{\rm Cyg}
(\ga_+,\ga_-)$ tending to $0$.

Then the subset of $\RR$ consisting of the $c'(\xi)$ for $\xi\in
\partial_\infty \HH^n_\CC$ which is neither a fixed point of a
conjugate of $\ga_0$ nor a parabolic fixed point, contains a segment
$[0,c]$ for some $c>0$. \cqfd
\ecoro

\section{Applications to Diophantine approximation} 
\label{sec:appdioapp} 

In this section, we apply the results of Sections \ref{sec:bounded}
and \ref{sec:hallray} to study the Diophantine approximation by 
quadratic irrational elements in $\RR$, $\CC$ and the Heisenberg group.

In order to obtain the Khinchin-type results in Subsection
\ref{subsec:appdioappRC} and in Subsection \ref{subsec:appdioappHeis},
we will apply the following result, which follows as a slight
extension of a particular case from \cite[Theorem 4.6]{HPpre}.  We
refer for instance to \cite{HPpre} for the general definitions of the
critical exponent $\delta=\delta_\Ga\in[0,+\infty]$ and of the
Patterson-Sullivan measure $\mu_{\xi_\infty,H_\infty}$ associated to a
horosphere $H_\infty$ with point at infinity $\xi_\infty$, for a
nonelementary discrete group of isometries $\Ga$ of a complete simply
connected Riemannian manifold $X$ with sectional curvature at most
$-1$. In this paper, we will only be interested in the particular
cases explained after the statement.

\btheo\label{theo:HPspir} \cite{HPpre} 
Let $X$ be a complete simply connected Riemannian manifold with
sectional curvature at most $-1$ and dimension at least $2$; let $\Ga$
be a discrete group of isometries of $X$ with finite covolume and
critical exponent $\delta$; let $\ga_0$ be a hyperbolic element of
$\Ga$ and $\R_{\Ga_0}$ be the set of points in $\partial_\infty X$
fixed by some conjugate of $\ga_0$ in $\Ga$; let $\xi_\infty$ be a
parabolic fixed point of $\Ga$ and $H_\infty$ be a horosphere centered
at $\xi_\infty$; and let $f:[0,+\infty[\,\ra\,]0,+\infty[$ be a slowly
varying map (as defined in Section \ref{subsec:nontrivial}).

If $\int_{1}^{+\infty} f(t)^{\delta}\;dt$ converges (resp.~diverges),
then $\mu_{\xi_\infty,H_\infty}$-almost no (resp.~every) point of
$\partial_\infty X - \{\xi_\infty\}$ belongs to infinitely many balls
of center $r$ and radius $f(D(r))e^{-D(r)}$ for the Hamenst\"adt
distance $d_{\xi_\infty,H_\infty}$, where $r$ ranges over
$\R_{\Ga_0}$.  \cqfd 
\etheo

In our applications in Section \ref{subsec:appdioappRC}
(resp.~\ref{subsec:appdioappHeis}), $X$ is the upper halfspace model
of $\HH^n_\RR$ (resp.~the Siegel domain model of $\HH^n_\CC$ as in
Example 1 of Section \ref{sec:framework}) and $H_\infty$ is the
horosphere, centered at $\xi_\infty=\infty$, of the points at
Euclidean height $1$ (resp.~$H_\infty= \{(w_0,w)\in \HH^n_\CC \;:\;
2\operatorname{Re} w_0 - |w|^2=2\}$). In this situation, since $\Ga$
has finite covolume,
\begin{itemize}
\item[$\bullet$] the critical exponent is $\delta=n-1$
  (resp.~$\delta=2n$, see for instance \cite[\S\ 6]{CI}),
\item[$\bullet$] the Hamenst\"adt distance $d_{\xi_\infty,H_\infty}$
  on $\partial_\infty X-\{\infty\}$ is the Euclidean distance, see
  \cite[\S\ 2.1]{HPMZ} (resp.~a multiple of the Cygan distance, see
  \cite[\S\ 3.11]{HPsurv}), and
\item[$\bullet$] the measure $\mu_{\xi_\infty,H_\infty}$ on
  $\partial_\infty X-\{\infty\}$ is the Hausdorff measure of
  $d_{\xi_\infty,H_\infty}$, which is the Lebesgue measure (resp.~is in
  the same measure class as the Hausdorff measure of the Cygan
  distance, see \cite{CI}).
\end{itemize}

\subsection{Approximation in $\RR$ and $\CC$ by  
  quadratic irrational elements}
\label{subsec:appdioappRC}

Let $K$ be either the field $\QQ$ or an imaginary quadratic extension
of $\QQ$, and correspondingly, let $\wh K$ be either $\RR$ or
$\CC$. Let $\O_K$ be the ring of integers of $K$. Denote by $\Kq$ the
set of quadratic irrational  elements in $\widehat K$ over $K$. For
every $\alpha\in \Kq$, let $\alpha^\sigma$ be its Galois conjugate
over $K$.

The group ${\rm PGL}_2(\wh K)$ acts on $\PP^1(K)=\wh K\cup\{\infty\}$
by homographies, and its subgroup ${\rm PGL}_2(\O_K)$ preserves $K$
and $\Kq$. Note that, for every $\alpha\in\Kq$ and every $\ga\in {\rm
  PGL}_2(\O_K)$, we have $(\ga\cdot\alpha)^\sigma= \ga\cdot
(\alpha^\sigma)$.

Let us fix a finite index subgroup $\Ga$ of ${\rm PSL}_2(\O_K)$. An
orbit of $\Ga$ in $\Kq$ will be called a {\it congruence class} in
$\Kq$ under $\Ga$. We are interested in Section
\ref{subsec:appdioappRC} in the approximation of elements of $\wh K$
by elements in the union of a fixed congruence class and of its
Galois conjugate.

For every $\alpha\in \Kq$, let 
$$
\E_{\alpha,\Ga}= \Ga\cdot \{\alpha,\alpha^\sigma\}\;,
$$
endowed with its Fr\'echet filter, and let
$$
h(\alpha)=\frac{2}{|\alpha-\alpha^\sigma|}\;.
$$ 
Clearly, $h(\alpha)$ belongs to $]0,+\infty[$ (as $\alpha\neq
     \alpha^\sigma$), and $h(\alpha^\sigma)=h(\alpha)$.  We will see
     in the proof of Theorem \ref{theo:rc} that points $r\in
     \E_{\alpha,\Ga}$ exit every finite subset of $\E_{\alpha,\Ga}$ if
     and only if $h(r)$ tends to $+\infty$.  Define the {\it
       quadratic Lagrange spectrum} relative to $(\alpha,\Ga)$ by
$$
\operatorname{Sp}_{\alpha,\Ga}=
\big\{\;c_{\alpha,\Ga}(\xi)=\liminf_{r\in \E_{\alpha,\Ga}}\;h(r)\,|\xi-r|\;:\;
\xi\in\wh  K -( K \cup\E_{\alpha,\Ga})\big\}\;.
$$

The following result is very classical, its proof (given for the sake
of completeness) was indicated to us by Y.~Benoist.

\blemm \label{lem:benoist} 
Let $\alpha\in\widehat K$. Then $\alpha$ is   quadratic irrational over
$K$ if and only if there exists a hyperbolic element $\ga$ in ${\rm
  PSL}_2(\O_K )$ having $\alpha$ as a fixed point, the other one then
being $\alpha^\sigma$.  
\elemm

\dem Let $\delta=\operatorname{dim}_{\RR} \widehat K$.  If
$\ga=\begin{pmatrix} a & b \\ c & d \end{pmatrix}$ is a hyperbolic
element in ${\rm SL}_2(\O_K )$, then its two fixed points (in
$\partial_\infty\HH^{\delta+1}_\RR=\wh K\cup\{\infty\}$) are distinct
solutions of the quadratic equation $ax+b=x(cx+d)$ with coefficients
in $\O_K $, and in particular they are  quadratic irrational and Galois
conjugated by $\sigma$.

Conversely, let $\alpha\in \Kq$. We refer for instance to \cite{Bor}
for general information on linear algebraic groups. Let
$$
T(\wh K)= \{\ga\in \operatorname{SL}_2(\wh K)\;:\;
\ga\cdot \alpha=\alpha,\;\;\ga\cdot \alpha^\sigma=\alpha^\sigma\}\;.
$$ Since $\alpha$ and $\alpha^\sigma$ are two distinct points in the
boundary of $\HH^{\delta+1}_\RR$, the subgroup $T(\wh K)$ is the set
of $\wh K$-points of an algebraic torus $T$ in $\operatorname{SL}_2$.
This torus $T$ is defined by a set of polynomial equations with
coefficients in $K(\alpha)$, this set being invariant by the Galois
group of $K(\alpha)$ over $K$. Hence $T$ is defined over $K$. Notice
that $T$ does not split over $K$, as the eigenvectors of an element of
$T$ in the affine plane are $(\alpha,1)$ and $(\alpha^\sigma,1)$, that
are not multiples of an element with coordinates in $K$. Hence, by
Borel-Harish-Chandra's theorem (see for instance
\cite[Theo.~12.3]{BHC}), the subgroup $T(\O_K )$ is a lattice in
$T(\wh K)$, and in particular is not trivial. That is, there exists an
element in $\operatorname{SL}_2(\O_K )$ having $\alpha$ (and
$\alpha^\sigma$) as fixed point. The result follows.  \cqfd

\medskip
\noindent{\bf Remark 1. } Let $K =\QQ$ and $\Ga=\PSLZ$. Note that the
Golden Ratio $\phi=\frac{1+\sqrt{5}}{2}$ is in the same congruence
class under $\Ga$ as its Galois conjugate $\frac{1-
  \sqrt{5}}{2}=-1/\phi$, hence $\E_{\phi,\Ga}= \Ga\cdot\phi$. 
On the other hand, $\frac{1+\sqrt{3}}{2}$ and $\frac{1-\sqrt{3}}{2}$ are Galois
conjugate, but are not in the same congruence class under $\Ga$, and
this second example is more typical.  Many papers have given necessary
and sufficient condition for when a quadratic irrational element is
in the same orbit under $\Ga$ as its Galois conjugate, see for
instance \cite{Sar,Lon,PR,Bur} and also \cite[Prop.~5.3]{BPP}.

\medskip
\noindent{\bf Remark 2. } Let $K =\QQ$. Let us give another expression
of the approximation constants $c_{\alpha,\Ga}(x)$.

For every $x\in \PP^1(\RR)$, let $\Ga_x$ be the stabilizer of $x$ in
$\Ga$. For every $\alpha\in \Qq$, endow the infinite set $\Ga/
\Ga_\alpha$ with its Fr\'echet filter. Denote by $N(\alpha)= \alpha
\alpha^\sigma$ the {\it norm} of an element $\alpha\in \Qq$. For every
element $\ga$ in $\Ga$, let $\ga=\pm\begin{pmatrix} a(\ga)& b(\ga) \\
  c(\ga)& d(\ga)\end{pmatrix}$.

\bprop Let $\alpha\in \Qq$ and let $\Ga$ be a finite index subgroup of
$\PSLZ$. For every $x\in\RR- (\QQ\cup \E_{\alpha,\Ga})$, we have
$$ 
c_{\alpha,\Ga}(x)=h(\alpha)\;\;\liminf_{(\ga,\varepsilon)
  \in(\Ga/\Ga_\alpha)\times\{1,\sigma\}}\;
\big|N\big(\alpha \,c(\ga)+d(\ga)\big)\big|\;\big|x-\ga
\cdot\alpha^\varepsilon\big| \;.
$$
\eprop

For instance, $\E_{\phi,\Ga}$ is the set of real numbers whose
continued fraction expansion is eventually constant equal to $1$, and
$\operatorname{Sp}_{\phi,\Ga}$ is equal to
$$ 
\Big\{\frac{2}{\sqrt{5}} \;
\liminf_{a,b,c,d\in\ZZ,\,ad-bc=1,\,d^2+dc-c^2\ra+\infty} \;
\big|d^2+dc-c^2\big|\;\Big|x-\frac{a\phi+b}{c\phi+d}\Big|\;:\;x\in
\RR-(\QQ\cup\E_{\phi,\Ga})\Big\}\;.
$$

\dem An easy computation shows that $|\ga\cdot\alpha - \ga\cdot
\alpha^\sigma| = |\alpha-\alpha^\sigma|/ N(c(\ga)\alpha+d(\ga))$,
hence
$$
h(\ga\cdot\alpha)= N(c(\ga)\alpha+d(\ga))\; h(\alpha)\;.
$$
The map $(\Ga/\Ga_\alpha)\times\{1,\sigma\}\ra\E_{\alpha,\Ga}$ defined
by $(\ga,\epsilon)\mapsto \ga\alpha^{\epsilon}$ is a bijection if
$\alpha$ and $\alpha^\sigma$ are not in the same congruence class,
and is a $2$-to-$1$ map otherwise. The result follows. 
\cqfd

\bigskip
\noindent{\bf Remark 3. }  
Assume again that $K=\QQ$. The quantity $h(\alpha)$ behaves in a very
different way from the naive height $H(\alpha)$ of $\alpha$ (defined
in the introduction).  Clearly, $h(\alpha+n)=h(\alpha)$ for every $n$
in $\NN$, but $H(\alpha+n)\to\infty$ as $n\to\infty$. Hence there
exists a sequence $(\alpha_i)_{i\in\NN}$ in $\Qq$ such that the ratio
$\frac{h(\alpha_i)}{H(\alpha_i)}$ tends to $0$ when $i\ra+\infty$. But
this ratio cannot tend to $+\infty$, as for every $\alpha\in\Qq$, we
have
$$ 
\frac{h(\alpha)}{H(\alpha)}\leq 2\;.
$$
Indeed, let $aX^2+bX+c$ be a minimal polynomial of $\alpha$ over
$\ZZ$, so that the naive height of $\alpha$ is
$H(\alpha)=\max\{|a|,|b|,|c|\}$. Then
$$
h(\alpha)=\frac{2}{\sqrt{(b/a)^2-4(c/a)}}= 
\frac{2|a|}{\sqrt{b^2-4ac}} \leq 2|a|\leq 2\;H(\alpha)\;,
$$ 
There are many possibilities for the relative behaviour of $h(\alpha)$
and $H(\alpha)$.  In particular, there exist sequences
$(\alpha_i)_{i\in\NN}$ in $\Qq$ and constants $c,c',c''>0$ such that
$h(\alpha_i)$ is equivalent as $i$ tends to $+\infty$ either to $c\;
H(\alpha_i)^{-\frac{1}{2}}$ (take $c=1$ and $\alpha_i=\sqrt{p_i}$ for
$p_i$ the $i$-th prime number), or to $c'\;H(\alpha_i)^{\frac{1}{2}}$
(take $c'=1$ and $\alpha_i=1/\sqrt{p_i}$ for $p_i$ the $i$-th prime
number) or to $c''\;H(\alpha_i)$ (take $c''=2/\sqrt{5}$ and $\alpha_i=
\frac{3 +2i+\sqrt{5}}{1+3i+i^2}$ for every $i\in\NN$).  This
difference is good to bear in mind when comparing our results with for
example the results of \cite{DS1,Spr,Bug} cited in the
Introduction. We refer to \cite[Lem.~5.2]{BPP} for a treatment of the
algebraic number theory aspects of $h(\alpha)$.

\btheo\label{theo:rc} Let $K =\QQ$ or $K=\QQ(i\sqrt{m})$ where $m$ is
a squarefree positive integer, let $\wh K =\RR$ or $\wh K =\CC$
respectively, and let $\delta=\operatorname{dim}_{\RR} \widehat
K$. Let $\Ga$ be a finite index subgroup of ${\rm PSL}_2(\O_K)$.
\begin{enumerate}
\item For every $\alpha_0\in \Kq$, the quadratic Lagrange spectrum
  $\operatorname{Sp}_{\alpha_0,\Ga}$ is closed, and equal to the
  closure of the set of the approximation constants
  $c_{\alpha_0,\Ga}(x)$ for $x$ a quadratic irrational over $K$, not
  in $\E_{\alpha_0,\Ga}$.
\item There exists $C\geq 0$ such that for every $\alpha_0\in \Kq$,
$$
\max\operatorname{Sp}_{\alpha_0,\Ga}\leq C\;.
$$ 
\item If $\widehat K =\CC$, then for every $\alpha_0\in \Kq$, there
exists $c>0$ such that ${\rm Sp}_{\alpha_0,\Ga}$ contains $[0,c]$.
\item Let $\alpha_0\in \Kq$ and $\varphi:\;]0,+\infty[\;\ra\;]0,+\infty[$
    be a map such that $t\mapsto \varphi(e^t)$ is slowly varying. If
    the integral $\int_{1}^{+\infty} \varphi(t)^\delta/t\;dt$ diverges
    (resp.~converges), then for Lebesgue almost every $x\in\widehat K
    $,
$$
\liminf_{r\in\E_{\alpha_0,\Ga}} \;
\frac{h(r)}{\varphi(h(r))}\;|x-r|=0\;
({\rm resp.} =+\infty)\;.
$$
\end{enumerate}
\etheo

\dem Consider the data $\D=(X,\Ga,\Ga_0,C_\infty)$, where $X$ is the
upper halfspace model of $\HH_\RR^{\delta+1}$, $\Ga$ is as in the
statement, $\Ga_0$ is the stabilizer in $\Ga$ of the translation axis
of a hyperbolic element of $\Ga$ one of whose fixed points is
$\alpha_0$, and $C_\infty$ is the set of points in $X$ with Euclidean
height at least $1$.  Note that $\Ga_0$ is nontrivial: By Lemma
\ref{lem:benoist}, the point $\alpha_0$ is a fixed point of a
hyperbolic element of ${\rm PSL}_2(\O_K)$, hence of a hyperbolic
element of $\Ga$, since $\Ga$ is a finite index subgroup of ${\rm
  PSL}_2(\O_K)$.

The data $\D$ satisfies the general assumptions of Example 1 of
Section \ref{sec:framework}: Since $z\mapsto z+1$ belongs to ${\rm
  PSL}_2(\O_K)$, the point $\infty$ is fixed by a parabolic element of
${\rm PSL}_2(\O_K)$, hence of $\Ga$; moreover, $C_\infty$ is precisely
invariant under the stabilizer of $\infty$ in ${\rm PSL}_2(\O_K)$,
hence in $\Ga$, by Shimizu's lemma; furthermore, the quotient of $X$
by ${\rm PSL}_2(\O_K)$, hence by $\Ga$, has finite volume.

Using the notations of  Example 1 of Section \ref{sec:framework}, we easily
check that 
$\R_{\Ga_0}=\E_{\alpha_0,\Ga}$. Furthermore, let $r\in \E_{\alpha_0,\Ga}$
and let $\ga_r\in\Ga$ be an element such that $\ga_r\ga_0\ga_r^{-1}$
fixes $r$. Recall that the other fixed point of a hyperbolic element
fixing $r$ is the Galois conjugate of $r$ over $K$, by Lemma
\ref{lem:benoist}. Hence, by Equation \eqref{eq:calcDr}, we have
\begin{equation}
D([\ga_r])=\log h(r)\;.\label{eqn:depth}
\end{equation}
Therefore, by Lemma \ref{lem:depthgrows}, points $r\in
\E_{\alpha_0,\Ga}$ exit every finite subset of $\E_{\alpha_0,\Ga}$ if
and only if $h(r)$ tends to $+\infty$.  The set of parabolic fixed
points of $\Ga$ is equal to the set of parabolic fixed point of ${\rm
  PSL}_2(\O_K )$, as $\Ga$ has finite index, hence it is equal to $ K
\cup\{\infty\}$.  Therefore
$$
\wh K -( K \cup\E_{\alpha_0,\Ga}) =\Lambda_c\Ga- \R_{\Ga_0}\;.
$$
For every $\xi$ in this set, we have
$$
c(\Ga_\infty\xi)= \liminf_{r\in \E_{\alpha_0,\Ga}}\;h(r)|\xi-r|
$$ 
by the first equality in Equation
\eqref{eq:calccstapproxconstcurv}. Hence, the quadratic Lagrange
spectrum $\operatorname{Sp}_{\alpha_0,\Ga}$ coincides with the
approximation spectrum $\spd{\D}$.

\medskip
We can now conclude that the assertions (1), (2) and (3) follow,
respectively, from Theorem \ref{theo:closedspec}, Proposition
\ref{prop:specialDirichlet}, and Corollary
\ref{coro:hallrayconstcurv}.  In particular, in (2) we get an upper
bound on $\operatorname{Sp}_{\alpha_0,\Ga}$ which depends only on $\Ga$,
and not on the (congruence class under $\Ga$) of $\alpha_0$.

\medskip To prove the assertion (4), define $f:t\mapsto \varphi(e^t)$,
which is slowly varying by the assumptions of (4). By an easy change
of variable, the integral $\int_{1}^{+\infty} f(t)^{\delta}\,dt$
diverges if and only if $\int_{1}^{+\infty} \varphi(t)^\delta/t \;dt$
diverges.  Hence by Theorem \ref{theo:HPspir}, by the comments
following it, and by Equation \eqref{eqn:depth}, if
$\int_{1}^{+\infty} \varphi(t)^\delta/t \;dt$ diverges, then for
almost every (for the Lebesgue measure) point $x$ in $\wh K$,
$$
\liminf_{r\in\E_{\alpha_0,\Ga}} 
\;\frac{h(r)}{\varphi(h(r))}\;|x-r| \leq 1.
$$ 
Replacing $\varphi$ by $\frac{1}{k}\varphi$ and letting $k\in\NN$ go
to $+\infty$, this proves the divergence part of the assertion (4) in
Theorem \ref{theo:rc}. The convergence part follows similarly.
\cqfd

\medskip \rem Replacing in the above proof, as in \cite{PP2,PP3,PP4}, 
$X$ by $\HH^5_\RR$, $\delta$ by $4$ and $\Ga$ by the image in the
isometry group of $X$ of a finite index subgroup of ${\rm SL}_2(\O')$
where $\O'= \ZZ(1+i+j+k)/2+\ZZ i+\ZZ j+\ZZ k$ is the Hurwitz order in
the usual Hamilton quaternion algebra $A$ over $\QQ$ (using
Dieudonn\'e's determinant), we could get similar Diophantine
approximation results of points in $A(\RR)$ by points in quadratic
extensions of $A(\QQ)$.

\medskip
To prove Theorem \ref{theo:goldratintro} of the Introduction, apply
Theorem \ref{theo:rc} with $K=\QQ$ (so that $\delta=1)$, $\Ga=\PSLZ$,
and $\varphi:t\mapsto t\psi(2/t)$, and notice that $\log(\varphi(e^{t}))$
is Lipschitz, and that $\int_1^{+\infty}\varphi(t)/t \,dt$ converges
if and only if $\int_0^{1}\psi(t)/t^2 \,dt$ converges.

\subsection{Approximation in the Heisenberg
  group}\label{subsec:appdioappHeis}

Let $m$ be a squarefree positive integer, let $K$ be the number field
$\QQ(i\sqrt{m})$, let $\O_K$ be its ring of integers, and let $\Kq$ be
the set of elements of $\CC$ which are  quadratic irrational over $K$.

Let $n\ge 2$, and let $(w',w)\mapsto w'\cdot \overline{w} =
\sum_{i=1}^{n-1} w'_i\overline{w_i}$ be the usual Hermitian scalar product on
$\CC^{n-1}$.  Consider the real Lie group
$$
{\rm Heis}_{2n-1}(\RR)=\{(w_0,w)\in\CC\times\CC^{n-1}\;:\;
2\operatorname{Re} w_0 - w\cdot \overline{w}=0\}
$$
whose law is
$$
(w_0,w)\cdot(w'_0,w')=(w_0+w_0'+w'\cdot w,w+w')\;.
$$ 
Endow it with the Cygan distance (see for instance \cite{Gol}),
which is the unique left-invariant distance such that
$$
d_{\rm Cyg}((w_0,w),(0,0))= \sqrt{2|w_0|},
$$ 
as well as with the {\it Cygan measure} (which is the Hausdorff
measure of the Cygan distance). The Lie group ${\rm Heis}_{2n-1}(\RR)$
is isomorphic to the standard $(2n-1)$-dimensional Heisenberg group,
that is for $n=2$ to the Lie group $\Big\{\left(\begin{array}{ccc} 1 &
    x & z \\ 0 & 1 & y \\ 0 & 0 & 1\end{array}\right)\;:\; x,y
,z\in\RR\Big\}$.

We are interested in the Diophantine approximation of the points of
${\rm Heis}_{2n-1}(\RR)$ by elements all of whose coordinates are
rational or quadratic over $K$, that is by elements of ${\rm Heis}_{2n-1}(\RR)\cap
{(K\cup\Kq)}^{n}$ (or nice subsets of it).

\medskip Let us identify ${\rm Heis}_{2n-1}(\RR)$ with its image in
the projective space $\PP_n(\CC)$ of $\CC^{n+1}$ by the map
$(w_0,w)\mapsto [w_0:w:1]$. Let $q$ be the Hermitian form of signature
$(1,n)$ on $\CC^{n+1}$ with coordinates $(z_0,z,z_n)$, defined by
$$
q=-z_0\overline{z_n}-z_n\overline{z_0}+z\cdot\overline{z}\;.
$$ 
The induced action on $\PP_n(\CC)$ of the special unitary group ${\rm
  SU}_q$ of $q$ preserves (the image in $\PP_n(\CC)$ of) ${\rm
  Heis}_{2n-1} (\RR)\cup\{\infty\}$. Let ${\rm SU}_q(\O_K)$ be the
arithmetic subgroup of ${\rm SU}_q$ which consists of matrices with
coefficients in $\O_K$. Note that the action of ${\rm SU}_q(\O_K)$ on
${\rm Heis}_{2n-1} (\RR)\cup\{\infty\}$ preserves both the set 
%(which is the set of rational points
%of a quadric defined over $\QQ$ which is not useful to
%explicit here, see for instance \cite[Sect.~6.3]{PP2})
$$
{\rm Heis}_{2n-1}(\QQ)\cup\{\infty\}=
\big({\rm Heis}_{2n-1}(\RR)\cap K^{n}\big)\cup\{\infty\}
$$ 
and the set $\big({\rm Heis}_{2n-1}(\RR)\cap K(\alpha)^{n}\big) \cup
\{\infty\}$ for every $\alpha\in\Kq$.

Let 
$$
\alpha_0=\frac{i}{2}\big(\sqrt{m+4}-\sqrt{m}\,\big)\;,
$$ 
which is an element of $\Kq$, since it is a root of the quadratic
polynomial $X^2+i\sqrt{m}\,X+1$ whose coefficients are in $K$, and it
does not belong to $K$. (We could have taken many other examples, but
$\alpha_0$ is one of the simplest ones.) The Galois conjugate of
$\alpha_0$ over $K$ is $\alpha_0^\sigma =
\frac{i}{2} (-\sqrt{m+4}-\sqrt{m})$.

Let $\Ga$ be a finite index subgroup of ${\rm SU}_q(\O_K)$, 
and let $\E'_{\alpha_0,\Ga}=\Ga\cdot\{(\alpha_0,0),(\alpha_0^\sigma,0)\}$.
For every $r\in\E'_{\alpha_0,\Ga}$, let $r^\sigma$ be the componentwise
Galois conjugate of $r$, and let
$$
h'(r)=\frac{1}{d_{\rm Cyg}(r,r^\sigma)}\;.
$$

Endow $\E'_{\alpha_0,\Ga}$ with its Fr\'echet filter. We will see in the
proof below that points $r\in\E'_{\alpha_0,\Ga}$ tend to infinity in
$\E'_{\alpha_0,\Ga}$ if and only if $h'(r)$ tends to $+\infty$.

In order to understand the Diophantine approximation of elements $\xi$
of ${\rm Heis}_{2n-1}(\RR)$ by elements in the subset $\E'_{\alpha_0,
  \Ga}$ of the set ${\rm Heis}_{2n-1}(\RR)\cap K(\alpha_0)^{n}$, we
introduce the {\it approximation constant} of $\xi$, defined by
$$
c'(\xi)=\liminf_{r\in \E'_{\alpha_0,\Ga}}\;h'(r)\,d_{\rm Cyg}(\xi,r)\;,
$$
and the {\it quadratic Heisenberg-Lagrange spectrum} with respect to
$(\alpha_0,\Ga)$
$$
\operatorname{Sp}'_{\alpha_0,\Ga}=\Big\{c'(\xi)\;:\;\xi\in {\rm
  Heis}_{2n-1}(\RR)-\big({\rm Heis}_{2n-1} (\QQ)\cup
\E'_{\alpha_0,\Ga}\big)\Big\}\;.
$$

\btheo \label{theo:applicapproxheis}
If $\Ga$ is  a finite index subgroup of ${\rm SU}_q(\O_K)$, then
\begin{enumerate}
\item[(1)] the Heisenberg-Lagrange spectrum
  $\operatorname{Sp}'_{\alpha_0,\Ga}$ is bounded;
\item[(2)] there exists $c>0$ such that $\operatorname{Sp}'_{\alpha_0,\Ga}$
contains an interval $[0,c[\;$;
\item[(3)] let $\varphi:\;]0,+\infty[\;\ra\;]0,+\infty[$ be a map such
  that $t\mapsto \varphi(e^t)$ is slowly varying, if the integral
  $\int_{1}^{+\infty} \varphi^{2n}(t)/t\;dt$ diverges
  (resp.~converges), then for Cygan almost every $x\in{\rm
    Heis}_{2n-1}(\RR)$,
$$
\liminf_{r\in\E'_{\alpha_0,\Ga}} \;
\frac{h'(r)}{\varphi(h'(r))}\;d_{\rm Cyg}(x,r)=0\;
({\rm resp.} =+\infty)\;.
$$
\end{enumerate}
\etheo

\dem Let $X$ be the Siegel domain model of the complex hyperbolic
$n$-space $\HH^n_\CC$, as in Example 2 of Section \ref{sec:framework},
identified, as well as its boundary at infinity, with its image in
$\PP_n(\CC)$ by $(w_0,w)\mapsto [w_0:w:1]$. The group $\Ga\subset
\operatorname{SL}_{n+1}(\CC)$ acts (with finite kernel) on
$X\subset\PP_n(\CC)$ as a discrete group of isometries of $X$ with
finite covolume, by the restriction of the projective action.  Note
that $\infty$ (corresponding to $[1:0:0]$) is a parabolic fixed point
of ${\rm SU}_q(\O_K)$, hence of $\Ga$.  Let
$$
{\ga_0}=
\left(\begin{array}{ccc} m+1 & 0 & -i\sqrt{m} \\ 0 & I & 0\\
    i\sqrt{m} & 0 & 1\end{array}\right)\;.
$$
It is easy to check that ${\ga_0}$ is a hyperbolic element of ${\rm
  SU}_q(\O_K)$ whose fixed points in $\partial_\infty X$ are exactly
$(\alpha_0,0)$ and $(\alpha_0^\sigma,0)$. Since $\Ga$ has finite index
in ${\rm SU}_q(\O_K)$, there exists $k\in \NN-\{0\}$ such that
${\ga_0}^k\in\Ga$.  Let $\Ga_0$ be the stabilizer in $\Ga$ of the
translation axis of $\ga_0$.

The horoball $\H_{2\sqrt{m}}$ is precisely invariant under the
stabilizer of $\infty$ in ${\rm SU}_q(\O_K)$, hence under $\Ga_\infty$
in $\Ga$, by the Kamiya-Parker inequality, see for instance
\cite[Lemma 6.4]{PP2}.  Thus, the data $\D=(X,\Ga,\Ga_0,
\H_{2\sqrt{m}})$ satisfy the conditions of Example 2 of Section
\ref{sec:framework}.  Using the notation of this example,
$\E'_{\alpha_0,\Ga}$ is the set $\R_{\Ga_0}$ of fixed points of the
conjugates of $\ga_0^k$ in $\Ga$. The set of parabolic fixed points of
$\Ga$ is equal to the set of parabolic fixed points of ${\rm SU}_q(\O_K)$,
which is exactly ${\rm Heis}_{2n-1} (\QQ)\cup\{\infty\}$ (see for
instance \cite[Sect.~6.3, Exam.~2]{PP2}). In particular,
$$
\Lambda_c\Ga-\R_{\Ga_0}={\rm Heis}_{2n-1}(\RR)-\big({\rm Heis}_{2n-1} 
(\QQ)\cup \E'_{\alpha_0,\Ga}\big)\;.
$$

For every $r$ in $\E'_{\alpha_0,\Ga}$, let $\ga_r\in\Ga$ be such that
$r$ is fixed by $\ga_r\ga_0^k\ga_r^{-1}$.  Note that if $r= \ga_r
(\alpha_0^\epsilon,0)$ for some $\epsilon\in\{1,\sigma\}$, then the
other fixed point of $\ga_r\ga_0^k\ga_r^{-1}$ is $\ga_r
((\alpha_0^{\epsilon})^\sigma,0)$, which is equal to $r^\sigma$ since
${\rm SU}_q$ acts projectively on ${\rm Heis}_{2n-1}(\RR) \cup
\{\infty\}$. Since the Cygan distance and the modified Cygan distance
are equivalent (see \cite [Sect.~6.1]{PP2}), there exists a constant
$c_1>0$ such that for every $r$ in $\E'_{\alpha_0,\Ga}$, with
$D([\ga_r])$ computed in Lemma \ref{lem:calcDcomphyp}, we have
\begin{equation}\label{eq:approxDhprim}
|\;D([\ga_r])-\log h'(r)\;|\leq c_1\;.
\end{equation}
Hence, it follows from Lemma \ref{lem:depthgrows} that points
$r\in\E'_{\alpha_0,\Ga}$ tend to infinity in $\E'_{\alpha_0,\Ga}$ if
and only if $h'(r)$ tends to $+\infty$.

Again, since the Cygan distance and the modified Cygan distance are
equivalent, there exists a constant $c_2>0$ such that for every $\xi\in
\Lambda_c\Ga-\R_{\Ga_0}$, we have
$$
\liminf_{r\in \E'_{\alpha_0,\Ga}}\; \frac{d_{\rm Cyg}(\xi,r)}{d_{\rm
    Cyg}(r,r^\sigma)}\leq c_2\;\liminf_{r\in \R_{\Ga_0}} \;\sqrt{2}\,
\frac{d'_{\rm Cyg}(r,r^\sigma)\;d_{\rm Cyg}(\xi,r)}{d_{\rm
    Cyg}(r,r^\sigma)^2}\;.
$$
The assertion (1) of Theorem \ref{theo:applicapproxheis} then follows
from Equation \eqref{eq:spdCygan} and from Proposition
\ref{prop:specialDirichlet}, and 
the assertion (2) follows from Corollary \ref{cor:hallcomphyp}.

\medskip To prove the assertion (3) of Theorem
\ref{theo:applicapproxheis}, consider the map $f:t\mapsto
\varphi(e^t)$, which is slowly varying.  In particular, by Equation
\eqref{eq:approxDhprim}, there exists a constant $c_3>0$ such that for
every $r\in\E'_{\alpha_0,\Ga}$, we have
$$
\frac{1}{c_3}\;\frac{\varphi(h'(r))}{h'(r)}\;\leq \;
f\big(D([\ga_r])\big)\;e^{-D([\ga_r])}
\;\leq\; c_3\;\frac{\varphi(h'(r))}{h'(r)}\;.
$$ Hence the assertion (3) of Theorem \ref{theo:applicapproxheis}
follows from Theorem \ref{theo:HPspir}, as in the proof of the
assertion (4) in Theorem \ref{theo:rc}. \cqfd

\noindent {\small 
\begin{tabular}{l}  
Department of Mathematics and Statistics, P.O. Box 35\\ 
40014 University of Jyv\"askyl\"a, FINLAND\\ 
{\it e-mail: parkkone@maths.jyu.fi} 
\end{tabular} 
\\
\smallskip
 \mbox{} 
\\ 
\begin{tabular}{l}  
D\'epartement de Math\'ematique et Applications, UMR 8553 CNRS\\ 
\'Ecole Normale Sup\'erieure, 45 rue d'Ulm\\ 
75230 PARIS Cedex 05, FRANCE\\ 
{\it e-mail: Frederic.Paulin@ens.fr} 
\end{tabular} 
}

\end{document}